%% file: main.tex
\documentclass[12pt, letter]{amsart}
\usepackage{graphicx} %

\DeclareRobustCommand{\SkipTocEntry}[5]{}
\usepackage{etoolbox}

\makeatletter
\let\old@tocline\@tocline
\let\section@tocline\@tocline
\newcommand{\subsection@dotsep}{4.5}
\newcommand{\subsubsection@dotsep}{4.5}
\patchcmd{\@tocline}
  {\hfil}
  {\nobreak
     \leaders\hbox{$\m@th
        \mkern \subsection@dotsep mu\hbox{.}\mkern \subsection@dotsep mu$}\hfill
     \nobreak}{}{}
\let\subsection@tocline\@tocline
\let\@tocline\old@tocline

\patchcmd{\@tocline}
  {\hfil}
  {\nobreak
     \leaders\hbox{$\m@th
        \mkern \subsubsection@dotsep mu\hbox{.}\mkern \subsubsection@dotsep mu$}\hfill
     \nobreak}{}{}
\let\subsubsection@tocline\@tocline
\let\@tocline\old@tocline

\let\old@l@subsection\l@subsection
\let\old@l@subsubsection\l@subsubsection

\def\@tocwriteb#1#2#3{%
  \begingroup
    \@xp\def\csname #2@tocline\endcsname##1##2##3##4##5##6{%
      \ifnum##1>\c@tocdepth
      \else \sbox\z@{##5\let\indentlabel\@tochangmeasure##6}\fi}%
    \csname l@#2\endcsname{#1{\csname#2name\endcsname}{\@secnumber}{}}%
  \endgroup
  \addcontentsline{toc}{#2}%
    {\protect#1{\csname#2name\endcsname}{\@secnumber}{#3}}}%

\newlength{\@tocsectionindent}
\newlength{\@tocsubsectionindent}
\newlength{\@tocsubsubsectionindent}
\newlength{\@tocsectionnumwidth}
\newlength{\@tocsubsectionnumwidth}
\newlength{\@tocsubsubsectionnumwidth}
\newcommand{\settocsectionnumwidth}[1]{\setlength{\@tocsectionnumwidth}{#1}}
\newcommand{\settocsubsectionnumwidth}[1]{\setlength{\@tocsubsectionnumwidth}{#1}}
\newcommand{\settocsubsubsectionnumwidth}[1]{\setlength{\@tocsubsubsectionnumwidth}{#1}}
\newcommand{\settocsectionindent}[1]{\setlength{\@tocsectionindent}{#1}}
\newcommand{\settocsubsectionindent}[1]{\setlength{\@tocsubsectionindent}{#1}}
\newcommand{\settocsubsubsectionindent}[1]{\setlength{\@tocsubsubsectionindent}{#1}}

\renewcommand{\l@section}{\section@tocline{1}{\@tocsectionvskip}{\@tocsectionindent}{}{\@tocsectionformat}}%
\renewcommand{\l@subsection}{\subsection@tocline{1}{\@tocsubsectionvskip}{\@tocsubsectionindent}{}{\@tocsubsectionformat}}%
\renewcommand{\l@subsubsection}{\subsubsection@tocline{1}{\@tocsubsubsectionvskip}{\@tocsubsubsectionindent}{}{\@tocsubsubsectionformat}}%
\newcommand{\@tocsectionformat}{}
\newcommand{\@tocsubsectionformat}{}
\newcommand{\@tocsubsubsectionformat}{}
\expandafter\def\csname toc@1format\endcsname{\@tocsectionformat}
\expandafter\def\csname toc@2format\endcsname{\@tocsubsectionformat}
\expandafter\def\csname toc@3format\endcsname{\@tocsubsubsectionformat}
\newcommand{\settocsectionformat}[1]{\renewcommand{\@tocsectionformat}{#1}}
\newcommand{\settocsubsectionformat}[1]{\renewcommand{\@tocsubsectionformat}{#1}}
\newcommand{\settocsubsubsectionformat}[1]{\renewcommand{\@tocsubsubsectionformat}{#1}}
\newlength{\@tocsectionvskip}
\newcommand{\settocsectionvskip}[1]{\setlength{\@tocsectionvskip}{#1}}
\newlength{\@tocsubsectionvskip}
\newcommand{\settocsubsectionvskip}[1]{\setlength{\@tocsubsectionvskip}{#1}}
\newlength{\@tocsubsubsectionvskip}
\newcommand{\settocsubsubsectionvskip}[1]{\setlength{\@tocsubsubsectionvskip}{#1}}

\patchcmd{\tocsection}{\indentlabel}{\makebox[\@tocsectionnumwidth][l]}{}{}
\patchcmd{\tocsubsection}{\indentlabel}{\makebox[\@tocsubsectionnumwidth][l]}{}{}
\patchcmd{\tocsubsubsection}{\indentlabel}{\makebox[\@tocsubsubsectionnumwidth][l]}{}{}

\newcommand{\@sectypepnumformat}{}
\renewcommand{\contentsline}[1]{%
  \expandafter\let\expandafter\@sectypepnumformat\csname @toc#1pnumformat\endcsname%
  \csname l@#1\endcsname}
\newcommand{\@tocsectionpnumformat}{}
\newcommand{\@tocsubsectionpnumformat}{}
\newcommand{\@tocsubsubsectionpnumformat}{}
\newcommand{\setsectionpnumformat}[1]{\renewcommand{\@tocsectionpnumformat}{#1}}
\newcommand{\setsubsectionpnumformat}[1]{\renewcommand{\@tocsubsectionpnumformat}{#1}}
\newcommand{\setsubsubsectionpnumformat}[1]{\renewcommand{\@tocsubsubsectionpnumformat}{#1}}
\renewcommand{\@tocpagenum}[1]{%
  \hfill {\mdseries\@sectypepnumformat #1}}

\let\oldappendix\appendix
\renewcommand{\appendix}{%
  \leavevmode\oldappendix%
  \addtocontents{toc}{%
    \protect\settowidth{\protect\@tocsectionnumwidth}{\protect\@tocsectionformat\sectionname\space}%
    \protect\addtolength{\protect\@tocsectionnumwidth}{2em}}%
}
\makeatother

\makeatletter
\settocsectionnumwidth{2em}
\settocsubsectionnumwidth{2.5em}
\settocsubsubsectionnumwidth{3em}
\settocsectionindent{1pc}%
\settocsubsectionindent{\dimexpr\@tocsectionindent+\@tocsectionnumwidth}%
\settocsubsubsectionindent{\dimexpr\@tocsubsectionindent+\@tocsubsectionnumwidth}%
\makeatother

\settocsectionvskip{5pt}
\settocsubsectionvskip{0pt}
\settocsubsubsectionvskip{0pt}

\settocsectionformat{\scshape}%
\settocsubsectionformat{\mdseries}
\settocsubsubsectionformat{\mdseries}
\setsectionpnumformat{\scshape}%
\setsubsectionpnumformat{\mdseries}
\setsubsubsectionpnumformat{\mdseries}

\let\oldtableofcontents\tableofcontents
\renewcommand{\tableofcontents}{%
  \vspace*{-\linespacing}%
  \oldtableofcontents}

\usepackage[dvipsnames]{xcolor}
\usepackage[draft]{say}
\definecolor{MHcol}{RGB}{100, 180, 255}

\definecolor{egyptianblue}{rgb}{0.06, 0.2, 0.65}
 \newcommand{\sayVS}[1]{\say[VS]{\color{egyptianblue}{\bf VS:}\;#1}}

\usepackage{amsthm}
\usepackage{amssymb}
\usepackage{amsmath}
\usepackage{parskip}
\usepackage{tikz}
\usepackage{tikz-cd}
\usepackage{bbm}
\usepackage{enumitem}
\usepackage{xcolor}
\usepackage{import}
\usepackage{yhmath}
\usepackage{aliascnt}
\usepackage{standalone}
\usepackage{caption}
\usepackage{subcaption}
\usepackage{mathtools}

\usepackage{mathrsfs}

\usepackage{comment}

\usepackage[margin=1in,marginparwidth=0.8in, marginparsep=0.1in]{geometry}

\usepackage[bookmarks=true, bookmarksopen=true,%
bookmarksdepth=3,bookmarksopenlevel=2,%
colorlinks=true,%
linkcolor=blue,%
citecolor=blue,%
filecolor=blue,%
menucolor=blue,%
urlcolor=blue]{hyperref}

\usepackage{mathtools}

\input{macros}

\ifdefined\DeclareUnicodeCharacter
  \DeclareUnicodeCharacter{FF1A}{:}
\fi

\DeclareMathOperator{\Ind}{Ind}

\DeclareMathOperator{\Fuk}{Fuk}

\newtheorem{dummy}{dummy}[section]
\newtheorem{lemma}[dummy]{Lemma}
\newtheorem{theorem}[dummy]{Theorem}

\newtheorem{corollary}[dummy]{Corollary}
\newtheorem{proposition}[dummy]{Proposition}

\theoremstyle{definition}
\newtheorem{definition}[dummy]{Definition}

\newtheorem{example}[dummy]{Example}
\newtheorem{remark}[dummy]{Remark}

\newtheorem{setup}[dummy]{Setup}

\begin{document}

\title{On sinks and mirrors}

\author{Mingyuan Hu} \author {Vivek Shende}

\begin{abstract}
    Given a full rank sink recurrent cluster $\mathcal{A}$ variety, we prove that its wrapped Fukaya category is equivalent to the category of coherent sheaves on a stacky elaboration of the dual Fock-Goncharov cluster $\mathcal{X}$ variety. 
    
    In particular,  braid varieties for a semisimple group are mirror to the braid varieties for the Langlands dual group.  Another example: for compact surfaces with no interior marked points, a least one marked point on each boundary component, and at least two marked points in total, the Fock-Goncharov moduli for SL(2) and PGL(2)  are homologically mirror. 
    The results were obtained in communion with ChatGPT Astra.
\end{abstract}

\maketitle
\setlength{\parskip}{0pt}

\renewcommand\contentsname{\vspace*{-30pt}}
\tableofcontents

\thispagestyle{empty}

\setlength{\parskip}{.5\baselineskip plus 2pt}

\begin{flushright}
\small\itshape
He that breaks a thing  to find out what it is has left the path of wisdom.\\[2pt]
\upshape --- Gandalf, in J.\,R.\,R.~Tolkien, \emph{The Fellowship of the Ring}.
\end{flushright}

\section{Introduction}

The Fock-Goncharov pairs of ``Langlands dual'' cluster varieties -- e.g., moduli of decorated local systems for Langlands dual groups -- have long been understood to be related by various forms of mirror symmetry \cite{FG09, FG06, GHK, GHKK, Goncharov-Shen, arguz-bousseau}, and been conjectured to be homologically mirror.  It is also  known that this conjectural homological mirror symmetry is compatible with the Lagrangian Floer theory between the compact real tori in the cluster charts \cite{PT, GLe}.  However, such tori cannot probe the complement of the cluster charts; addressing this ``deep locus'' has been an essential challenge in mirror symmetry for cluster varieties \cite{Castronovo-Gorsky-Simental-Speyer}.  

In this article, we establish homological mirror symmetry for all full rank sink-recurrent cluster varieties.  This class of cluster varieties was introduced in \cite{GLSS}, and shown to include various kinds of braid varieties, open Richardson varieties, open positroid varieties, etc.  In Appendix \ref{app:surfaces}, we add another example: for compact surfaces with no interior marked points, at least one marked point on each boundary, and at least two boundary marked points in total,  the Fock-Goncharov \cite{FG06} cluster structure for decorated $\mathrm{SL}_2$ local systems is also sink-recurrent.  (This condition on boundary points amounts to asking that the local systems have only irregular singularities, with total Katz invariant at least 1.) 

In particular, braid varieties for a semisimple group are homologically mirror to the braid varieties of the Langlands dual group.  Such varieties include the moduli of Stokes data for connections on $\mathbb{P}^1$ which have no singularities away from $\infty$, and there have a sufficiently irregular singularity; thus we show these are homologically mirror to the corresponding moduli for the Langlands dual group.  
We similarly obtain homological mirror symmetry between the aforementioned $\mathrm{SL}_2$ and $\mathrm{PGL}_2$ Fock-Goncharov moduli spaces.

We recall the basic setup of cluster algebra \cite{FZ-I}. 
Fix the following data: 
\[
 \widetilde B=\binom BC,\qquad B\in\operatorname{Mat}_{m\times m}(\bZ),
 \quad C\in\operatorname{Mat}_{f\times m}(\bZ)
\]
and assume there is a positive rational diagonal matrix $D$
satisfying
\begin{equation}\label{eq:symmetrizer}
 BD=-DB^{\mathsf T}.
\end{equation}

In the cluster algebra literature, $B$ is the {\em exchange matrix}, $\widetilde B$ is the {\em extended exchange matrix}, and the data is also referred to as a {\em skew-symmetrizable seed}. We write $M$ for the set of row indices of $B$, and $F$ for the set of row indices of $C$, so that $|M|=m$ and $|F|=f$.  
When \(B\) is skew-symmetric, it is often recorded by the quiver  where the number of 
arrows from $i\to j$ is $B_{ij}$ if this is positive, and zero otherwise.

Given $k \in M$, the {\em mutation} $\mu_k \widetilde B$ is another matrix of the same size, given by the formula: 
\[
 (\mu_k\widetilde B)_{ij}=\begin{cases}
 -\widetilde B_{ij},&i=k\text{ or }j=k,\\
 \widetilde B_{ij}+[\widetilde B_{ik}]_+\widetilde B_{kj}+\widetilde B_{ik}[-\widetilde B_{kj}]_+,&i,j\ne k.
 \end{cases}
\]

We write $\sC(\widetilde B) := \bC[z_i]_{i\in M} [z_j^{\pm}]_{j \in F}$, and $\sC^*(\widetilde B) := \bC[z_i^{\pm}]_{i\in M \sqcup F}$.  These $z_i$ are called {\em cluster variables}, and said to be {\em mutable} if $i\in M$, and {\em frozen} otherwise. 
One embeds $\sC(\mu_k \widetilde B) = \bC[z'_i]_{i\in M} [{z'}_j^{\pm}]_{j \in F}$ in  $\sC^*(\widetilde B)$ by taking $z_i' = z_i$ for $i \ne k$, and imposing the {\em exchange relation}: 
\begin{equation}\label{eq:exchange}
 z_kz_k'=\prod_i z_i^{[\widetilde B_{ik}]_+}
              +\prod_i z_i^{[-\widetilde B_{ik}]_+},
 \qquad [a]_+=\max(a,0).
\end{equation}
The {\em cluster algebra}  $\sA(\widetilde B)$ is by definition the span of the $\sC(\widetilde B') \subset \bC(z_i:i\in M\sqcup F)$ over all iterated mutations $\widetilde B'$ of $\widetilde B$, and the {\em upper cluster algebra} is the corresponding intersection 
$\sU(\widetilde B) := \bigcap_{\widetilde B'} \sC^*(\widetilde B') \subset \bC(z_i:i\in M\sqcup F)$. 
A fundamental result of the theory is that  the cluster variables associated to any iterated mutation are in fact Laurent polynomials in the original variables, i.e.
$\sC(\widetilde B') \subset \sC^*(\widetilde B)$, from which it follows that $\sA(\widetilde B) \subset \sU(\widetilde B)$ \cite{FZ-I}.

For the class of {\em locally acyclic} cluster algebras \cite{M}, $\sA(\widetilde B) = \sU(\widetilde B)$ and this algebra is finitely generated \cite[Theorem~4.1 and 4.2]{M}.  Here we will be interested in the  class of {\em sink-recurrent} cluster varieties, which are known to be locally acyclic 
(\cite[Proposition~5.5]{GLSS} for quivers and  \cite[Proposition~3.7]{GLSB} in general).  Let us recall the definition.

For a mutatable vertex \(k\in M\), let \(N(k) = \{j \in M \mid B_{jk} \neq 0 \}\).  We say \(k\in M\) is a \emph{sink} if \(B_{jk} \ge 0\) for all \(j\in M\). 

\begin{definition}[{\cite[Definition~5.4]{GLSS}}]
\label{def:sink-recurrence}
Sink-recurrent skew-symmetrizable
 matrices are the smallest mutation-closed class containing the zero matrix, and such that if $k \in M$ is a sink, and $B|_{M \setminus \{k\}}$ and $B|_{M \setminus (\{k\} \sqcup N(k))}$ are sink-recurrent, then $B$ is also sink-recurrent. 
\end{definition}

For an extended exchange matrix $\widetilde B= (B;C)$, we write $\cA(\widetilde B)=\Spec\sA(\widetilde B)$. 

We recall the `principal coefficient' matrix and associated structures: 
\[
B_{\mathrm{prin}}=\binom B{I_m},\qquad
P^\vee=\cA((-B^{\mathsf T})_{\mathrm{prin}})\times(\bC^*)^f,
\qquad H=(\bC^*)^m
\]
The torus $H$ acts on the initial coordinates $(X;Z;U)$ of $P^\vee$
with weight matrix
\begin{equation}
\label{eq:weights}
 \begin{pmatrix}
     I_m & B^{\mathsf T} & C^{\mathsf T}
 \end{pmatrix}.
\end{equation}
This weight matrix annihilates the exchange matrix. 
Every exchange relation is therefore homogeneous, and the action extends
to the  cluster algebra. Following \cite{FG09, FG06, GHK, GHKK}, we write:
\begin{equation}
\label{eq:canonical-mirror}
 \cX(\widetilde B^\vee):=[P^\vee/H].
\end{equation}
Some of the cited  authors have taken the coarse moduli space; we take the stack, as did \cite{Goncharov-Shen}.  These agree on the locus of cluster tori. 

There is also a description of \(\cX(\widetilde B^\vee)\) as a finite quotient, which however depends on additional choices.
Choose an integral matrix with nonvanishing determinant:
\begin{equation}\label{intro:algcompletion}
 T=\begin{pmatrix}-B^{\mathsf T}&-C^{\mathsf T}\\
 C'&T_{\rm ff}\end{pmatrix} \in \Mat_{m+f}(\bZ),
 \qquad
 \widetilde B_T^\vee=\binom{-B^{\mathsf T}}{C'}.
\end{equation}
Then \(G_T \coloneqq \Hom(\coker(T), \bC^*)\) is a finite group, and
 Proposition~\ref{prop:finite-atlas} shows:
\begin{equation}
\label{eq:intro_DM_mirror}
    \cX(\widetilde B^\vee)\cong[\cA(\widetilde B_T^\vee)/G_T].
\end{equation}

Under certain circumstances,  $\mathcal{A}(\widetilde B)$ and $\cX(\widetilde B^\vee)$ can be identified.  We consider the following classes of possible $\widetilde{B}$, at least the first two of which are standard \cite{GLSB}: 
\begin{enumerate}
    \item Full rank (FR):
    \(
        \operatorname{rank}_{\bQ}\widetilde B = m.
    \)
    \item Really full rank (RFR):
    \(
        \widetilde B^{\mathsf T}\bZ^{m+f}=\bZ^m.
    \)
    Equivalently, $\widetilde B \bZ^m$ is saturated.

    \item Unimodular: the mutable columns extend to a matrix
    \(E\in GL_{m+f}(\bZ) \),
    \(    E_{M\sqcup F,M}=\widetilde B, \)
    and there exists a positive rational diagonal $(m+f) \times (m+f)$ matrix
    $\widetilde D$ with
    \[
        \widetilde D|_M=D,\qquad
        E\widetilde D=-\widetilde D E^{\mathsf T}.
    \]
\end{enumerate}
Apparently, \( (3) \Rightarrow (2) \Rightarrow (1)\).

In the locally acyclic case, \eqref{eq:intro_DM_mirror}  shows that \(\cX(\widetilde B^\vee)\) is a smooth affine variety if \(\widetilde B\) is RFR, and is isomorphic to \(\cA(\widetilde B)\) if there is a matrix 
\[
    T =
    \begin{pmatrix}
        B&-C^{\mathsf T}\\
        C &T_{\mathrm{ff}}
    \end{pmatrix}
    \in GL_{m+f}(\bZ).
    \]
 Note that we do not ask \(T_\mathrm{ff}\) to be skew-symmetric. 

The space $\cX$ also carries tori charts, which however transform by a dual mutation formula.

For a seed $\widetilde B'$ obtained by mutating $\widetilde B$ once, we write $\cX^{\mathrm{tori}}_{\widetilde B^\vee, (\widetilde B')^\vee}$ for the corresponding union of $\cX$ cluster charts.  We denote the restriction maps:
\begin{equation} \label{eq: one mutation X comparison}
\Perf_{\C[x_1^{\pm1}, \ldots, x_{m+f}^{\pm1}]} \leftarrow \Coh(\cX^{\mathrm{tori}}_{\widetilde{B}^\vee, (\widetilde{B}')^\vee}) \rightarrow \Perf_{\C[x_1'^{\pm1}, \ldots, x_{m+f}'^{\pm1}]} 
\end{equation}

Given a parameterized torus $\T = (S^1)^{m+f}$, we fix $\Loc(\T) = \Perf_{\C[x_1^{\pm1}, \ldots, x_{m+f}^{\pm1}]}$ by identifying $x_i$ with the monodromy of the $i$'th $S^1$.  
Given two parameterized exact Lagrangian tori $\T, \T'$ in a Weinstein symplectic manifold $Y$, the Viterbo restriction functors give maps 
\begin{equation} \label{eq: lagrangian cluster comparison}
\Perf_{\C[x_1^{\pm1}, \ldots, x_{m+f}^{\pm1}]} \leftarrow \Fuk(Y) \rightarrow \Perf_{\C[x_1'^{\pm1}, \ldots, x_{m+f}'^{\pm1}]} \end{equation} 

\begin{definition}
    Let $\widetilde{B}$ be a skew-symmetrizable seed and let $Y$ be a  Liouville symplectic manifold.  A   Lagrangian cluster ensemble for $\widetilde{B}$ on $Y$ consists of the data of a collection of parameterized exact Lagrangian tori $T_{\widetilde{B}'}: (S^1)^{m+f} \to Y$, indexed by iterated mutations $\widetilde{B}'$ of $\widetilde{B}$, such that for every pair of adjacent cluster charts $\widetilde{B}', \widetilde{B}''$,  there is a map $\Fuk(Y) \to \Coh(\cX^{\mathrm{tori}}_{(\widetilde{B}')^\vee, (\widetilde{B}'')^\vee})$ intertwining \eqref{eq: lagrangian cluster comparison} and \eqref{eq: one mutation X comparison}:
    \[
 \begin{tikzcd}[column sep=large]
\Perf_{\C[x_1'^{\pm1}, \ldots, x_{m+f}'^{\pm1}]} \arrow[d, equal] & \Fuk(Y) \arrow[l] \arrow[r] \arrow[d] & \Perf_{\C[x_1''^{\pm1}, \ldots, x_{m+f}''^{\pm1}]} \arrow[d, equal] \\
\Perf_{\C[x_1'^{\pm1}, \ldots, x_{m+f}'^{\pm1}]}  & \Coh(\cX^{tori}_{(\widetilde{B}')^\vee, (\widetilde{B}'')^\vee}) \arrow[l] \arrow[r] & \Perf_{\C[x_1''^{\pm1}, \ldots, x_{m+f}''^{\pm1}]} 
 \end{tikzcd}
\]
Note that by Zariski descent, each such map is unique if it exists. 
\end{definition}

It follows from the results of \cite{PT, GLe} that for a Liouville structure on $Y$ for which the $(S^1)^{m+f}$ in the cluster charts are exact, these tori will form a Lagrangian cluster ensemble. We  construct such a Liouville structure in Proposition \ref{prop:lambda_ad} below.

When $\cX(\widetilde B^\vee)$  is a smooth affine scheme, its derived category is generated by the structure sheaf, so there should be some corresponding Lagrangian which generates $\Fuk(\cA(\widetilde B))$.  There is a known candidate: since cluster transformations are positive, the positive real locus of  various cluster tori coincide; we denote said locus as $L_+$.  It is evidently totally real.

\begin{theorem} \label{thm:main} For any full rank sink-recurrent skew-symmetrizable \(\widetilde B\), there is a Weinstein structure on $\cA(\widetilde B)$ for which $L_+$ is a strictly exact Lagrangian, and for which the central compact torus in each cluster chart is exact.   There is  an equivalence 
$\Fuk( \cA(\widetilde B) ) \simeq \Coh( \cX (\widetilde B^\vee))$ carrying $L_+ \mapsto \mathcal{O}$ and exhibiting that said tori form a Lagrangian cluster ensemble.  

When $\cX (\widetilde B^\vee)$ is a scheme (then necessarily affine and smooth), we also have a mirror symmetry 
$\Fuk(\cX (\widetilde B^\vee)) \cong \Coh(\cA( \widetilde{B}))$.
\end{theorem}

The proof consists of the following ideas.  We first prove the result in the unimodular case, and then use the fact that, for $\Gamma$ finitely generated abelian, $\Gamma$-covers are mirror to $\Hom(\Gamma, \mathbb{G}_m)$ quotients to obtain the general case, similarly to  \cite{Sh21}.

On the $A$-side, we realize the space as a parameterized handle attachment to a smaller sink  recurrent cluster variety $\cA(\widetilde B_{[k]})$, along  $T^* S^1$ times another smaller sink recurrent cluster variety $\cA(\widetilde B_{\Sigma})$.  By sectorial descent \cite{GPS2}, this reduces us to computing the partially wrapped Fukaya category of the former, stopped at the latter.  We use a new descent result, given in Theorem \ref{thm: SOV descent}, to compute this as a pullback of the un-stopped category for $\cA(\widetilde B_{[k]})$, and a product of a standard stopped category $((T^* S^1)^2, T^*S^1)$ with the category for $\cA(\widetilde B_{\Sigma})$.  
The mirror categories have corresponding mirror decompositions; the parameterized handle attachment is mirror to a certain affine blowup (as in \cite{GLe}), and the aforementioned new descent is mirror to a Zariski cover.  We use the known mirror relationship between tori to rigidify the diagrams to make matching them easy. 

In the main we work with a specially designed primitive $\lambda_{\mathrm{ad}}$ for the imaginary part of the holomorphic symplectic form. However, in Appendix \ref{app:affine}, we show using \cite{HS26} that the resulting Fukaya category agrees with that for the standard primitive of the K\"ahler form.

We now mention some examples to which this result applies. 

After \cite{CGGLSS, GLSS, GLSB}, it is known that braid varieties are sink-recurrent cluster varieties.  There are various closely related versions of braid varieties;  we use that of   \cite{CGGLSS}.  Let us recall the definition.
Fix a connected semisimple algebraic group $G$, with Borel subgroup $B_+$, and its Langlands dual $G^\vee$.
For a positive braid word $\beta=\sigma_{i_1}\cdots\sigma_{i_l}$, recall the Demazure product is defined to be \(\delta (\beta) \coloneqq s_{i_1} \star \cdots \star s_{i_l}\), where
\[
w\star s_i=
\begin{cases}
ws_i,&\ell(ws_i)>\ell(w)\\
w,&\ell(ws_i)<\ell(w)
\end{cases}
\]
and \(s_i \in N_G(T)/T\) is the reflection corresponding to the generator of braid group \(\sigma_i\).  
The braid variety is
\[
 X_G(\beta)=\left\{(F_0,\ldots,F_l)\in(G/B_+)^{l+1}\ \middle|\
 \begin{array}{l}F_0=B_+,\quad F_l= \delta(\beta)B_+,\\ F_{j-1} \xrightarrow{s_{i_j}}F_j  \quad(1\le j\le l)\end{array}\right\}.
\]
As indicated by the notation, $X_G(\beta)$ depends only on $\beta$  and not its expression as a word. The flag variety, and so also these braid varieties, depend only on the Lie algebra of $G$. We identify the labeled Coxeter systems of $G$ and $G^\vee$.

\begin{corollary}\label{cor:braid}
$\Fuk\bigl(X_G(\beta)\bigr)\simeq \Coh\bigl(X_{G^\vee}(\beta)\bigr).$
\end{corollary}

\begin{proof}
By~\cite[Sections~4.8--4.9 and~6.3]{GLSB}, braid varieties are sink-recurrent and really full rank.
By~\cite[Proposition~6.14]{CGGLSS}, the cluster seed $\widetilde B_G(\beta)$ for $X_G(\beta)$ is dual to $\widetilde B_{G^\vee}(\beta)$ for $X_{G^\vee}(\beta)$.
Write $\varepsilon$ for the full exchange matrix in the convention of~\cite{CGGLSS}; in our notations above, 
$
 \widetilde B_G(\beta)=(\varepsilon_{M,M\sqcup F})^{\T}$ and $\widetilde B_{G^\vee}(\beta)=-\varepsilon_{M\sqcup F,M}$.
By~\cite[Equation~(45), Lemmas~8.2 and~8.4]{CGGLSS}, there is $p\in GL_{M\sqcup F}(\bZ)$ with $p_{ij}=\varepsilon_{ij}$ unless both indices are frozen.
Take $T=-p$ in~\eqref{intro:algcompletion}. Then 
\[
 T_{M,M\sqcup F}=-\widetilde B_G(\beta)^{\T},\qquad
 \widetilde B_T^\vee=\widetilde B_{G^\vee}(\beta),\qquad G_T=1.
\]
Proposition~\ref{prop:finite-atlas} then identifies $\cX(\widetilde B_G^\vee(\beta)) = X_{G^\vee}(\beta)$, and we conclude by  Theorem~\ref{thm:main}.
\end{proof}

\begin{example}
    The positroid strata of  Grassmannians have long been known to be expressed in terms of (type A) braids \cite{KLS}, though in a slightly different formalism (which was also shown to be sink-recurrent in \cite{GLSB}).  They were shown to be braid varieties in the sense above in \cite{CGGS}.  It follows that these are self-mirror. 
\end{example}

\begin{example}
$G$-local systems on $\mathbb{P}^1$ with one irregular singularity (and no regular singularities) may be described in terms of their Stokes data, which (essentially by definition) is a cyclic sequence $(F_0, F_1, \ldots, F_k)$ of elements in $G/B_+$, with adjacent terms in prescribed positions characterized by the formal type of the irregular singularity; see e.g. discussions in \cite{boalch-yamakawa}.  Writing these adjacency conditions as $\sigma_0, \ldots, \sigma_k$, one can see readily that if $\sigma_0 \ldots \sigma_k = w_0 \beta$, where $w_0$ is the longest word, $\beta$ is positive, and $\delta(\beta) = w_0$, then this parameter space is just $X_G(\beta)$.  

It is expected that mirror symmetry for character varieties is in some sense a microlocalization of Betti geometric Langlands \cite{benzvi-nadler} to the locus of stable Higgs bundles (see e.g. \cite{shende-hitchinfibers, nadler-shende-higgs}); correspondingly, we expect that Corollary \ref{cor:braid} should go a long way towards establishing some version of Betti geometric Langlands  for $\mathbb{P}^1$ with one irregular puncture.  Other recent results on mirror symmetry in this context include \cite{jin-yun, BHOS}. 
\end{example}

\begin{example}
    The special case where $\beta$ is a power of the half-twist gives \cite[Conjecture~1.5, equation~(35)]{Goncharov-Shen} for $n \ge 3$. 
\end{example}
 
A different class of examples comes from the Fock-Goncharov \cite{FG06} moduli of decorated local systems for marked surfaces.  
For a semisimple simply connected group $G$ and marked surface $S$, we write $\mathcal{A}_{G, S}^\times$ for the moduli of $G$-local systems with markings in $G/U$, and adjacent boundary markings transverse.  For the Langlands dual (adjoint type) group $G^\vee$, we write $\mathcal{P}^\times_{G^\vee, S}$ for the variant of the $\mathcal{X}$ space introduced in \cite{GS-quantum}, which  instead parameterizes local systems with transverse $G^\vee/B^\vee$ markings and ``pinnings''.  It is shown in \cite{GS-quantum} (see also \cite{Goncharov-Shen, Shen-Poisson, IO-Wilson, IOS}) that the 
$\mathcal{A}_{G, S}^\times$ and $\mathcal{P}^\times_{G^\vee, S}$ are the $\mathcal{A}$ and $\mathcal{X}$ varieties for the Langlands dual seeds originally proposed in \cite{FG06}.  So long as there are boundary markings and no interior markings, $\mathcal{P}^\times_{G^\vee, S}$ is non-stacky.   

We show in Proposition \ref{prop:surface} that when \(S\) has at least one marked point on each boundary component, and at least two marked points in total, the seed for \(\cA_{SL_2, S}\) is sink-recurrent. 

We conclude: 

\begin{corollary}
When $S$ is a compact connected surface, with no internal markings, all boundaries marked, and with at least two total markings,  
$\Fuk\bigl(\mathcal A^\times_{\mathrm{SL}_2,S})
 \simeq \Coh\bigl(\mathcal \cP^\times_{\mathrm{PGL}_2,S}\bigr)$ and $\Fuk\bigl(\mathcal \cP^\times_{\mathrm{PGL}_2,S})
 \simeq \Coh\bigl(\mathcal A^\times_{\mathrm{SL}_2,S}\bigr).$ 
\end{corollary}

{\bf AI Statement.} 
All results of this article are obtained through discussion with ChatGPT-Astra.  Here we give a short history. 
After two days of discussion, ChatGPT-Astra was first able to compute skeleta and microsheaves first for $A_{2n}$ cluster varieties, then eventually the acyclic case (Appendix~\ref{app:acyclic-skeleton}).
We asked it to compute the skeleton for more general cases. It was not able to do this, but in the process, it found the sink cover~\eqref{eq:sink_cover} and Theorem~\ref{thm: SOV descent} in special cases (in a somewhat obfuscated form). 
We investigated this gluing formula, and formulated and proved the results first in the unimodular setting (with the assistance of Astra).
Astra originally could not handle the full rank case, but we had the idea of using the methods of \cite{Sh21} to descend from the unimodular case.  We also extracted and significantly generalized many intermediate steps into their present forms, especially Theorem \ref{thm: SOV descent}.

{\bf Acknowledgements.} 
After obtaining the results, we learned that a collaboration of Linhui Shen, Zhe Sun, Daping Weng, and Peng Zhou \cite{SSWZ}, and related but separate work of Josephine Hlavinka \cite{Hlavinka-mirror}, have pursued the acyclic case by heuristically determining  skeleta and computing microsheaves.    We similarly learned that  Hlavinka had previously  proposed that sink-recurrent cluster varieties may be a good place to study mirror symmetry.
We are supported by Villum Fonden Villum Investigator 37814.  

\section{Sylvan-Orlov-Viterbo descent}

Here we prove the following new descent result for partially wrapped Fukaya categories.

\begin{theorem}\label{thm: SOV descent}
Let $O\subset V$ be Weinstein domains 
whose contact boundaries coincide in 
the neighborhood of a Weinstein stop $F \subset \partial_\infty O \cap \partial_\infty V$.  Consider the commutative diagram formed by the Orlov functors, Viterbo restriction and stop removal: 
\begin{equation}\label{eq:s4-viterbo-pullback}
 \begin{tikzcd}[column sep=large]
 \Fuk(F) \arrow[r, "J"] \arrow [dr, "J_0"]
 & \Fuk(V,F)\arrow[r]\arrow[d,"q"']
   &\Fuk(V)\arrow[d]\\ & 
 \Fuk(O,F)\arrow[r]&\Fuk(O)
 \end{tikzcd}
\end{equation}
Assume $F$ is swappable, and the monodromy for the resulting spherical twist on $\Fuk(V, F)$ preserves the kernel of $q$, and
$J_0 c =0 \implies J c = 0$.  Then the square
is a pullback. 
\end{theorem}

In fact, once we use the standard properties of the Orlov and Viterbo functors \cite{GPS2, Sy}, Theorem \ref{thm: SOV descent} is really a categorical fact about spherical functors, which we prove first in Corollary \ref{cor:pullback from spherical}.

A sufficient condition for both swappability and the triviality of monodromy is the existence of a global circle action, positive for the stop: 

\begin{lemma}
\label{lem:circle-twist}
Let $V$ be a Liouville domain, and let  $F \subset \partial_{\infty} V$ be a Liouville hypersurface.   Suppose given a strictly exact circle action denoted for 
$\alpha \in \bR/2\pi \bZ$ by 
$\rho_\alpha: V \to V$, and assumed equipped with a lift to grading/orientation data.  Assume  $\rho_\alpha(F) \cap F = \emptyset$ for $\alpha \in (0, 2\pi)$, and that 
$\alpha\mapsto\rho_\alpha(F)$ is a positive isotopy.

Assume $\Fuk(V, F) \to \Fuk(V)$ is the quotient by the image of the  Orlov functor $J:\Fuk(F)\to\Fuk(V,F)$, e.g. $F$ Weinstein.  
Then $J$ is spherical, and the spherical twist is the identity. 
\end{lemma}
\begin{proof}
The loop $\alpha\mapsto\rho_\alpha(F)$ is a stop swap in the sense of Sylvan, hence $J$ is spherical by \cite[Thm. 1.3]{Sy}, with (inverse) twist given by transport around this loop.
On the other hand, we may use the circle action to trivialize the entire family of sectors, $
 \rho_\alpha:(V,F)\xrightarrow{\sim}(V, \rho_\alpha F)$.
It follows that transport around the loop is given by the identity functor. 
\end{proof}

\subsection{A pullback criterion}

Throughout this section, all categories are stable presentable $\infty$-categories.

\begin{lemma}\label{lem:magic-pullback}
Let $q:\cP\to\cQ$ be a localization, and assume the fully faithful right adjoint $q^R$ preserves sequential colimits (e.g. $\cP$ compactly generated and $q$ preserves compact objects).
Let $T: \cP \to \cP$ and $T_0: \cQ \to \cQ$ be autoequivalences, equipped with $qT\simeq T_0q$ and natural transformations
\[
 f:\id_{\cP}\longrightarrow T,\qquad
 f_0:\id_{\cQ}\longrightarrow T_0.
\]
Assume $qf\simeq f_0q$, $Tf\simeq fT$, and $T_0f_0\simeq f_0T_0$, coherently with iteration.
Write $\eta:\id_{\cP}\to q^Rq$ for the adjunction unit.
Suppose the following square of endofunctors induces Cartesian squares when applied to objects: 
\[
\begin{tikzcd}[column sep=large]
 \id_{\cP}\ar[r,"f"]\ar[d,"\eta"']
 &T\ar[d,"\eta T"]\\
 q^Rq\ar[r,"q^Rf_0q"']&q^RT_0q.
\end{tikzcd}
\]
Let $\cP[f^{-1}]\subset\cP$ be the full subcategory of objects $X$ for which $f_X:X\to TX$ is an equivalence.
Define $\cQ[f_0^{-1}]\subset\cQ$ similarly.
Then their localization functors $L,L_0$ give a pullback square
\[
\begin{tikzcd}[column sep=large]
 \cP\ar[r,"L"]\ar[d,"q"']
 &\cP[f^{-1}]\ar[d,"q"]\\
 \cQ\ar[r,"L_0"']&\cQ[f_0^{-1}].
\end{tikzcd}
\]
\end{lemma}

\begin{proof}
Note that the compatibility assumptions between $T, T_0, q$ and between $f, f_0, q$ amount to canonical isomorphisms $T_0 \simeq q T q^R$ and $f_0 \simeq q f q^R$. 

Let us recall the telescope formula for the localizations:
\[
\begin{aligned}
 LX&=\colim\bigl(X\xrightarrow{f_X}TX
       \xrightarrow{T(f_X)}T^2X\longrightarrow\cdots\bigr),\\
 L_0Y&=\colim\bigl(Y\xrightarrow{(f_0)_Y}T_0Y
       \xrightarrow{T_0((f_0)_Y)}T_0^2Y\longrightarrow\cdots\bigr).
\end{aligned}
\]
Here, 
the maps $X\to LX$ and $Y\to L_0Y$ come from the first terms.
To see these implement the localization, observe the relation $Tf\simeq fT$ identifies $f_{LX}:LX\to TLX$ with the map to the shifted telescope, hence an equivalence.
If $f_Z$ is invertible, precomposition induces equivalences
$\Hom_{\cP}(T^{n+1}X,Z)\to\Hom_{\cP}(T^nX,Z)$ by naturality of $f$ and invertibility of $T$.
Consequently, $
 \Hom_{\cP}(LX,Z)
 \simeq\lim_n\Hom_{\cP}(T^nX,Z)
 \simeq\Hom_{\cP}(X,Z)$, so  $L$ is left adjoint to the inclusion $\cP[f^{-1}]\subset\cP$; the same argument applies to $L_0$.

Taking right adjoints of $qT\simeq T_0q$ gives
$T^{-1}q^R\simeq q^RT_0^{-1}$, hence $Tq^R\simeq q^RT_0$.
This identification carries $fq^R$ to $q^Rf_0$.
Since $q$ preserves colimits and $q^R$ preserves sequential colimits, the telescope formulas give
\[
 qL\simeq L_0q,\qquad Lq^R\simeq q^RL_0,
\]
compatibly with the localization units.

For $X\in\cP$, set $K_X=\fib(X\xrightarrow{\eta_X}q^RqX)$.
The vertical fibers of the assumed square are $K_X$ and $TK_X$.
Thus $f_{K_X}$ is invertible and $K_X\simeq LK_X$.
Applying $L$ to the fiber sequence defining $K_X$ shows that the unit maps give
\[
 X\xrightarrow{\sim}LX\times_{q^RL_0qX}q^RqX.
\]

An object on the RHS is a triple $(a,b,\alpha)$ with
$a\in\cP[f^{-1}]$, $b\in\cQ$, and $\alpha:qa\xrightarrow{\sim}L_0b$.
Send it to
 $a\times_{q^RL_0b}q^Rb$,
where the arrows are $a\xrightarrow{\eta_a}q^Rqa\xrightarrow{q^R\alpha}q^RL_0b$ and $q^R(b\to L_0b)$.
Applying $q$ recovers $b$, since $qq^R\simeq\id$ and $\alpha$ is invertible.
Applying $L$ recovers $a$, since $Lq^Rb\to Lq^RL_0b$ is identified with the equivalence $q^RL_0b\to q^RL_0^2b$.
These identifications recover $\alpha$ as well.
Together with the preceding formula for $X$, they give mutually inverse functors.
\end{proof}

\begin{remark}
    The sequential-colimit assumption in Lemma \ref{lem:magic-pullback} can be weakened to requiring that the canonical comparison $Lq^R\to q^RL_0$ be invertible.
\end{remark}

Such morphisms of endofunctors arise naturally in the setting of spherical functors.  We recall this notion now. 

For the remainder of this subsection we fix morphisms of stable categories 
$\cC \xrightarrow{J} \cP \xrightarrow{q} \cQ$ 
with the following two properties: 
\begin{itemize}
    \item The map $q$ is a localization.  
    We write $q^R$  for its fully faithful right adjoint, 
    \(\eta:\id_{\cP}\to q^R q\) for its unit, and $
  \cK:=\ker q$.  
  
  We recall
    the counit \(qq^R\to\id_{\cQ}\) is an equivalence,  
 the essential image of
\(q^R\) is the full subcategory 
$\cK^\perp
  :=
  \{X\in\cP:\Hom_{\cP}(K,X)\simeq0\text{ for every }K\in\cK\}$, and
for each \(X\in\cP\), the cone of \(\eta_X:X\to q^RqX\) lies in \(\cK\).
\item $J$ is spherical.  In one of the equivalent formulations from \cite[Theorem~1.1]{AL}, this means: 
\begin{itemize}
    \item there are left and right adjoints  $J^L\dashv J\dashv J^R$ 
    \item the twist $M$ and cotwist $\Phi$ defined by the following exact triangles 
\begin{equation}\label{eq:twist}
JJ^R\longrightarrow\id_{\cP}\longrightarrow M,
\end{equation}
\[
  \id_{\cC}\longrightarrow J^RJ\longrightarrow\Phi.
\]
are invertible.
\item the composition  
\begin{equation}\label{eq:spherical-comparison}
  J^R\longrightarrow J^RJJ^L\longrightarrow\Phi J^L.
\end{equation}
is an equivalence. 
\end{itemize}

\end{itemize}

We will be interested in the properties of $J_0 := q \circ J$.  We note first $J_0$ has  a right adjoint   $J_0^R := J^R q^R$, and define $M_0: \cQ \to \cQ$  by the exact triangle
\[
J_0 J_0^R \longrightarrow\id_{\cQ} \longrightarrow M_0.
\]

\begin{lemma} 
    If \(M(\cK)\subset\cK\), then 
$qM \simeq M_0 q$. 
\end{lemma}
\begin{proof}
If \(M(\cK)\subset\cK\), then  \(qM\) annihilates \(\cK\), so it factors through \(q\) as some $M_0' q$ by
the universal property of the localization.  Then $M_0' qq^R = q M q^R$, and $q \, \circ$ \eqref{eq:twist} $\circ \, q^R$ gives $M_0' \simeq M_0$. 
\end{proof}

\begin{lemma}\label{lem:formal-descent}
Assume $M(\cK) \subset \cK$. Then
 the following conditions are equivalent:
\begin{enumerate}[label=\textup{(\roman*)}]
\item \(M(\cK)=\cK\).
\item \(M_0\) is an equivalence.
\item The map  \(M q^R\to q^R M_0\), induced by the
localization unit, is an equivalence.
\end{enumerate}
\end{lemma}

\begin{proof}

If (i) holds, both \(M\) and \(M^{-1}\) descend to the quotient.
Their descended functors are inverse to one another, proving (ii).
Conversely, if (ii) holds, then
$
  qMX\simeq0
  \Longleftrightarrow
  M_0 q X\simeq0
  \Longleftrightarrow
  qX\simeq0.
$
Thus \(M^{-1}(\cK)=\cK\), equivalently \(M(\cK)=\cK\), proving (i). 

Assuming (ii), 
one can take right adjoints of
\(qM\simeq M_0 q\) to obtain
\(M^{-1} q^R\simeq q^R M_0^{-1}\); rearranging gives  (iii). 

Finally, assume (iii). Since \(q^R\) and \(M\) are fully faithful, the
equivalence \(q^R M_0 \simeq Mq^R\) implies that \(M_0\)
is fully faithful. Since \(qM\) is
essentially surjective and \(qM\simeq M_0 q\), we see that $M_0$ is also essentially surjective.  This proves  (ii).
\end{proof}

\begin{lemma}\label{lem:local} 
The following conditions are equivalent:
\begin{enumerate}[label=\textup{(\roman*)}]
\item Composing $J$ with the unit for the localization gives an equivalence  \(J\to q^R qJ = q^R J_0\).
\item \(J^L|_{\cK}=0\).
\item \(J^R|_{\cK}=0\).
\item For every $K \in \cK$ and $c \in \cC$, we have $\Hom_{\cP}(K, Jc) = 0$. 
\item For every $K \in \cK$ and $c \in \cC$, we have $\Hom_{\cP}(Jc, K) = 0$. 
\end{enumerate}
Under these conditions, $M$ is the identity on $\cK$. 
\end{lemma}

\begin{proof}
The equivalence of (i) and (iv) is a standard property of localization.  
For \(K\in\cK\) and \(c\in\cC\), adjunction gives
$  \Hom_{\cP}(K,Jc)\simeq\Hom_{\cC}(J^L K,c)$.
Consequently, every \(Jc\) is local if and only if \(J^L K=0\) for all
\(K\in\cK\). This proves the equivalence of (iv) and (ii).  A similar argument shows equivalence of (v) and (iii). 
The equivalence of (ii) and (iii) follows from
\(J^R\simeq\Phi J^L\) and invertibility of \(\Phi\).

The fact that (iii) implies $M$ is the identity on $\cK$ is immediate from \eqref{eq:twist}. 
\end{proof}

\begin{lemma} \label{lem: categorical vcs criterion}
If $M(\cK) \subset \cK$ and $J_0 c = 0 \implies Jc = 0$, the equivalent conditions of Lemma \ref{lem:local} hold. 
\end{lemma}
\begin{proof}
Applying \(q\) to \eqref{eq:twist} at \(K\) gives an exact triangle $
  J_0 J^R K\longrightarrow qK\longrightarrow qMK. $
  The second term vanishes by definition, and the last because $M(\cK) \subset \cK$,  so \(J_0 J^R K=0\).
By assumption, then $JJ^RK=0$;
adjunction then gives
$\Hom_{\cC}(J^RK,J^RK)\simeq\Hom_{\cP}(JJ^RK,K)\simeq0$,
so $J^R K=0$.
\end{proof}

\begin{theorem}\label{thm: descending spherical functor}
Suppose  the equivalent conditions of Lemma \ref{lem:local} hold.
Then the functor \(J_0:\cC\to\cQ\) has left and right adjoints $J_0^L := J^L  q^R$ and $J_0^R := J^R q^R$.
In addition, $q$ adjunction morphisms induce isomorphisms $J^L \xrightarrow{\sim}  J_0^L q$ and $J^R \xrightarrow{\sim} J_0^R q$. 
Finally, $J_0$ is spherical, with
twist $M_0$, and cotwist $\Phi$.  \end{theorem}

\begin{proof}
For \(y\in\cQ\) and \(c\in\cC\), using Lemma \ref{lem:local} (i)  gives
\[
  \Hom_{\cC}(J^L q^R y,c)
  \simeq\Hom_{\cP}(q^R y,Jc)
  \simeq\Hom_{\cP}(q^R y,q^R J_0c)
  \simeq\Hom_{\cQ}(y,J_0c),
\]
so $J_0^L:= J^L q^R$ is left adjoint to $J_0$.  

For every \(X\in\cP\), the cone of \(X\to q^R qX\) belongs to \(\cK\).
Applying  Lemma \ref{lem:local} (ii) and (iii) gives:
\[
  J^L \xrightarrow{\sim}J^L q^R q=J_0^L q,
  \qquad
  J^R \xrightarrow{\sim}J^R q^R q=J_0^R q.
\]

Per Lemma \ref{lem:local}, $M$ is the identity on $\cK$, hence by Lemma \ref{lem:formal-descent}, the twist $M_0$ is an equivalence.

The  cotwist $\Phi_0$ of \(J_0\) is
defined by the exact triangle 
$
  \id_{\cC}\longrightarrow J_0^RJ_0\longrightarrow\Phi_0
$.
Since \(J\simeq q^R J_0\), we have 
$
  J_0^R J_0=J^R q^R J_0\simeq J^R J
$, and so 
$\Phi_0\simeq\Phi$; in particular, $\Phi_0$ is an equivalence.

Finally, composing $q^R$ with 
\eqref{eq:spherical-comparison} gives  
$
  J_0^R =J^R q^R \xrightarrow{\sim}\Phi J^L q^R \simeq \Phi_0 J_0^L.
$
Thus $J_0$ is spherical. 
\end{proof}

\begin{corollary} \label{cor:pullback from spherical} 
Suppose $\mathcal{P}$ is compactly generated, $q$ preserves compact objects.   Assume $M(\cK)\subseteq\cK$ and  $J_0 c = 0 \implies Jc =0$.   Then there is a pullback diagram
    \[
\begin{tikzcd}
 \cP
 \ar[r]\ar[d]
 &
 \cP / \operatorname{Image(J)}
 \ar[d]
 \\
 \cQ
 \ar[r]
 &
 \cQ / \operatorname{Image(J_0)}.
\end{tikzcd}
\]

\end{corollary}
\begin{proof}
We consider the natural morphism of twist triangles:  
    \[
\begin{tikzcd}
 J J^R
 \ar[r]\ar[d]
 &
 \mathrm{id}_{\cP}
 \ar[d] \ar[r, "f"] & T \ar[d] 
 \\
 q^R J_0 J_0^R q \ar[r] & q^R q 
 \ar[r] & q^R T_0 q
\end{tikzcd}
\]
Here, $T_0 = q T q^R$ is the descended twist. 
We have assumed the hypothesis of Lemma \ref{lem: categorical vcs criterion}, hence obtain the equivalent conditions of Lemma \ref{lem:local}, hence can invoke Theorem \ref{thm: descending spherical functor}.   Thus  $JJ^R \to q^R J_0 J_0^R q$ is an isomorphism, hence the right hand square is Cartesian.  This verifies the hypothesis of Lemma \ref{lem:magic-pullback}.  It remains only to note that the $f$-torsion of $\mathcal{P}$ is generated by the image of $JJ^R$, which is the same as the image of $J$, since $Jc$ is a retract of $JJ^RJc$; similarly for $\mathcal{Q}$. 
\end{proof}

\subsection{Proof of Theorem \ref{thm: SOV descent}}
We have assumed the hypotheses needed to invoke \cite[Thm. 1.3]{Sy}, which then asserts  $J$ is spherical. 
Per \cite{GPS2}, Viterbo restriction of Weinstein sectors is a localization, and removal of a Weinstein stop is the localization by the image of the Orlov functor.    We may check whether the desired square is a pullback after ind-completing. The result now follows from  Corollary \ref{cor:pullback from spherical}. \qed

\section{Cluster algebra}

For a matrix \(X \in \Mat_{a \times b}\), and \( R \subset \{1, \dots, a\}\), \( J \subset \{1, \dots, b\} \) we use \(X_{R, J}\) to denote the submatrix with rows \(R\) and columns \( J\), and \( X_{\bullet, J} \coloneqq X_{ \{1, \dots, a\}, J} \). Similarly, \(X_{R, \bullet} \coloneqq X_{R, \{1, \dots, b\}}\). 

\subsection{Freezing and sink covers}

Fix a seed $\widetilde B = \binom BC$, with mutable indices $M$ and frozen indices $F$. 

For $J\subset M$, the \emph{freezing} of $\widetilde B$ at $J$ is the matrix
\[
 \widetilde B_J=\widetilde B_{M\sqcup F,\,M\setminus J},
\]
with mutable index set $M\setminus J$ and frozen index set $F\sqcup J$.
Its mutable part is $B_{M\setminus J,M\setminus J}$.
Its frozen rows are those of $C_{F,M\setminus J}$ and $B_{J,M\setminus J}$.
By \cite[Proposition~3.1]{M},
\begin{equation}\label{eq:freezing-sandwich}
 \sA(\widetilde B_J)\subset\sA(\widetilde B)[z_j^{-1}:j\in J]\subset\sU(\widetilde B_J).
\end{equation}
The freezing is said to be a \emph{cluster localization} if the first inclusion is an equality \cite[Definition~3.3]{M}.
Then $\cA(\widetilde B_J)$ is the open subset $\{\prod_{j\in J}z_j\ne0\}$ of $\cA(\widetilde B)$.
Evidently, this holds whenever $\sA(\widetilde B_J)=\sU(\widetilde B_J)$.

Now let $k$ be a sink of $B$, and put
\[
 I=M\setminus(\{k\}\sqcup N(k)).
\]
 Consider the following four matrices:  
\begin{equation}
    \label{eq:matrices_B}
\widetilde B_{[k]}\coloneqq\widetilde B_{M\sqcup F,\,M\setminus\{k\}}, \qquad 
\widetilde B_{N(k)} \coloneqq \widetilde B_{M\sqcup F,\,\{k\}\sqcup I}, \qquad
 \widetilde B_{\setminus k} \coloneqq \widetilde B_{(M\setminus\{k\})\sqcup F,\,I}
, \qquad
\widetilde B_\cN \coloneqq \binom01
\end{equation}
Observe \( \cN  := \cA(\widetilde B_{\cN}) \cong \Spec\bC[x, b, t^\pm]/ (xb - 1 - t)\).

Let \(x = z_k\) be the cluster variable corresponding to the sink. The exchange relation reads
\begin{equation}
\label{eq:sink_exchange}
 xb=1+t,\qquad 
 b=z_k'/u_0,\qquad u_0=\prod_{a\in F}z_a^{[-C_{ak}]_+},\qquad
 t=u\prod_{j\in N(k)}z_j^{B_{jk}},\qquad u=\prod_{a\in F}z_a^{C_{ak}}.
\end{equation}

\begin{proposition}[See also {\cite[(4.16)]{Hlavinka}}]
\label{prop:AB_Nk}
We have 
\[
 \sA(\widetilde B_{N(k)}) \cong \sA(\widetilde B_{\setminus k})[x,b]/(xb-1-t)
\]
\end{proposition}
\begin{proof}
After freezing \(N(k)\), the vertex \(k\) is isolated from other mutable vertices, i.e. \(B_{ik} = 0 \) for all \( i \in I\).  Thus the cluster variables of \( \widetilde B_{N(k)}\) are those of \( \widetilde B_{\setminus k}\) together with \(x \) and \( u_0 b\), and \( xb = 1 +t\). Since \(u_0\) is a product of frozen variables, hence invertible, 
there is a surjection
\begin{equation}
\label{eq:pre_splitting}
 \sA(\widetilde B_{\setminus k})[x,b]/(xb-1-t)\longrightarrow\sA(\widetilde B_{N(k)}).
\end{equation}
It is the identity on $\sA(\widetilde B_{\setminus k})$ and sends $x\to z_k$ and $b\to z_k'/u_0$.

Since $xb=1+t$, the source is a free $\sA(\widetilde B_{\setminus k})$-module spanned by $x^i$ ($i\ge0$) and $b^j$ ($j\ge1$).
Their images in $ \sA(\widetilde B_{N(k)})\subset  \sC^*(\widetilde B)$ are $x^i$ and $(1+t)^jx^{-j}$, which are linearly independent over \(\sA(\widetilde B_{\setminus k})\). 
Hence \eqref{eq:pre_splitting} is injective.
\end{proof}

Consider the localizations
\[
 \cA(\widetilde B)_x\coloneqq\{x\ne0\},\qquad
 \cA(\widetilde B)_t\coloneqq\{t\ne0\}=\Bigl\{\prod_{j\in N(k)}z_j\ne0\Bigr\}.
\]
Since $(x,t) = 1$, we have a cover: 
\begin{equation}
\label{eq:sink_cover}
 \cA(\widetilde B)=\cA(\widetilde B)_x\cup \cA(\widetilde B)_t.
\end{equation}

If  $\widetilde B_{[k]}$ and $\widetilde B_{N(k)}$ are cluster localizations, then: 
\[
 \cA(\widetilde B)_x=\cA(\widetilde B_{[k]}),\qquad
 \cA(\widetilde B)_t=\cA(\widetilde B_{N(k)}).
\]

\subsection{Change of frozen variables}

\begin{lemma}
\label{lem:coefficient-change}
Let
\[
    \widetilde B_0=\binom BC,\qquad
    \widetilde B_1=\binom B{PB+QC},
\]
where
$P\in\operatorname{Mat}_{f\times m}(\bZ)$ and
$Q\in\operatorname{Mat}_{f\times f}(\bZ)$ has nonzero determinant.

Write $(x;y)$ and $(u;v)$ for the initial mutable and frozen
coordinates of the two seeds. Then
\[
\begin{aligned}
    \sA(\widetilde B_0)\otimes_{\bC[y^{\pm1}]}\bC[v^{\pm1}]
    &\cong \sA(\widetilde B_1),\\
    \sU(\widetilde B_0)\otimes_{\bC[y^{\pm1}]}\bC[v^{\pm1}]
    &\cong \sU(\widetilde B_1),
\end{aligned}
\]
where the coefficient map and the isomorphisms are given by
\[
    y^a\longmapsto v^{Qa},
    \qquad
    x_i\longmapsto u_i v^{P_i},
\]
and $P_i$ denotes the $i$th column of $P$.
In particular, if $Q\in GL_f(\bZ)$, then
$\sA(\widetilde B_0)\cong \sA(\widetilde B_1)$ and \( \sU(\widetilde B_0) \cong \sU(\widetilde B_1)\)
\end{lemma}

\begin{proof}
Choose representatives $\Gamma\subset\bZ^f$ for $\bZ^f/Q\bZ^f$.
Via $y^a\mapsto v^{Qa}$, we have
\begin{equation}
\label{eq:free_module}
 \bC[v^{\pm1}]=\bigoplus_{\gamma\in\Gamma}v^\gamma\bC[y^{\pm1}].
\end{equation}
Since $\det Q\ne0$, this is a finite free extension.
Thus both \(
\mathscr A(\widetilde B_0)
\otimes_{\mathbb C[y^{\pm1}]}\mathbb C[v^{\pm1}],
\) 
and 
\(
\mathscr U(\widetilde B_0)
\otimes_{\mathbb C[y^{\pm1}]}\mathbb C[v^{\pm1}]
\)
embed in $\bC(x,v)$.
Consider the isomorphism of fields 
\[
\phi:\bC(x,v)\longrightarrow\bC(u,v), \quad \phi(x_i)=u_iv^{P_i}, \quad \phi(v_a)=v_a. 
\]
For each mutation sequence $t$, write $x_i(t)$ and $u_i(t)$ for the corresponding cluster variables.
Since
\[
 \widetilde B_1=\begin{pmatrix}I_m&0\\P&Q\end{pmatrix}\widetilde B_0,
\]
Proposition~4.3 of~\cite{Fra} gives
\[
 \phi(x_i(t))=v^{p_i(t)}u_i(t),\qquad p_i(t)\in\bZ^f.
\]
The frozen monomials are units, so
\[
\begin{aligned}
 \phi\!\left(\sA(\widetilde B_0)\otimes_{\bC[y^{\pm1}]}\bC[v^{\pm1}]\right)
 &=
\mathbb C[v^{\pm1}]
\bigl[v^{p_i(t)}u_i(t)\mid i,t\bigr]\\
&= \bC[v^{\pm1}][u_i(t)\mid i,t]\\
 &=\sA(\widetilde B_1).
\end{aligned}
\]

Put $L_{0,t}=\bC[x_i(t)^{\pm1}\ (1\le i\le m),y^{\pm1}]$ and
$L_{1,t}=\bC[u_i(t)^{\pm1}\ (1\le i\le m),v^{\pm1}]$.
The same rescaling gives
$\phi(L_{0,t}\otimes_{\bC[y^{\pm1}]}\bC[v^{\pm1}])=L_{1,t}$.
Taking intersections over all mutation sequences gives
\[
\begin{aligned}
 \phi\!\left(\sU(\widetilde B_0)\otimes_{\bC[y^{\pm1}]}\bC[v^{\pm1}]\right)
 &=\bigcap_t\phi\!\left(L_{0,t}\otimes_{\bC[y^{\pm1}]}\bC[v^{\pm1}]\right)\\
 &=\bigcap_tL_{1,t}=\sU(\widetilde B_1).
\end{aligned}
\]
\end{proof}

\subsection{Splitting $\mathcal{A}(\widetilde{B})_t$}\label{sec:At-splitting}
Now we assume that \(\widetilde B\) is unimodular.
Fix a unimodular completion $E$ of $\widetilde B$, and a diagonal matrix $\widetilde D=\diag(d_i)_{i\in M\sqcup F}$ such that $K=E\widetilde D$ is skew-symmetric.
Let $N_{\bZ}=\bZ^{m+f}$ be the cocharacter lattice of the seed torus for $\widetilde B$, with basis $e_i$.
Let $\epsilon_i$ be the dual basis of the character lattice $N_{\bZ}^*$.
A character $a\in N_{\bZ}^*$ gives the Laurent monomial $z^a$.
The form $K$ defines a skew pairing $K(v,w)=v^{\T}Kw$ on $N_{\bQ}$.
We write $K(\cdot,v)\in N_{\bQ}^*$ for the character $w\mapsto K(w,v)$.
For $i\in M\sqcup F$ put
\begin{equation}\label{eq:dual-circles}
 c_i= E_{ji} \epsilon_j \in N_{\bZ}^*,\qquad w_i=(E^{\T})^{-1}e_i\in N_{\bZ}.
\end{equation}
A direct computation shows
\begin{equation}\label{eq:dual-circle-identities}
 K(\cdot,e_i)=d_ic_i,\qquad K(\cdot,w_i)=-d_i\epsilon_i,\qquad c_i(w_j)=\delta_{ij},\qquad \epsilon_i(w_i)=0 .
\end{equation}

Let \(k\) be a sink.
Then $c_k$ is supported on $F\sqcup N(k)$, and $t=z^{c_k}$. 
Since $c_k(w_k)=1$, every $a\in\bZ^{F\sqcup N(k)}$ decomposes as $a=a(w_k)\,c_k+\bigl(a-a(w_k)c_k\bigr)$, and
\begin{equation}\label{eq:frozen-splitting}
 \bZ^{F\sqcup N(k)}=\bZ c_k\oplus\Lambda_k,\qquad
 \Lambda_k \coloneqq \{a\in\bZ^{F\sqcup N(k)}:a(w_k)=0\}.
\end{equation}
Choosing a basis gives an identification 
$ \Lambda_k \cong \bZ^{ f + |N(k)| - 1}$. 

For $j\in I$, consider the column vector $f_j:=\widetilde B_{F\sqcup N(k),\,j} \in \bZ^{ f + |N(k)|}$. Its projection in \(\Lambda_k\) defines a vector \( [f_j] \in \bZ^{f+|N(k)|-1}\). 
We write $C_\Sigma$ for the matrix whose columns are the $[f_j]$, and define: 
\[
 \widetilde B_\Sigma=\binom{B_{I,I}}{C_\Sigma}, \qquad \Sigma\coloneqq \cA(\widetilde B_\Sigma). 
\]
\begin{lemma}
    Change of basis of $\Lambda_k$ induces isomorphisms of the corresponding $\mathscr{A}(\widetilde{B}_\Sigma)$.
\end{lemma}
\begin{proof}
A change of basis of $\Lambda_k$ replaces $C_\Sigma$ by $QC_\Sigma$ for some $Q\in GL_{f+|N(k)|-1}(\bZ)$.
So the claim follows from Lemma~\ref{lem:coefficient-change}. 
\end{proof}

\begin{proposition}
\label{prop:splitting_ABprime}
There is an isomorphism 
$\sA(\widetilde B_{\setminus k})\cong\sA(\widetilde B_\Sigma)[t^{\pm1}]$ sending $z^{c_k}\mapsto t$, and correspondingly 
$\cA(\widetilde B_{\setminus k}) \cong \bC^*_t \times \cA(\widetilde{B}_\Sigma)$. 
\end{proposition}

\begin{proof}
The seed $\widetilde B_{\setminus k}$ has mutable indices $I$, frozen indices $F\sqcup N(k)$, and
\[
 \widetilde B_{\setminus k}=\binom{B_{I,I}}{C'},\qquad C'=(f_j)_{j\in I}.
\]

Choose a basis $a_1,\dots,a_r$ of $\Lambda_k$, $r=f+|N(k)|-1$.
In this basis the columns of $C_\Sigma$ are given by
\begin{equation}\label{eq:Sigma-matrix}
 f_j-f_j(w_k)\,c_k=\sum_\alpha(C_\Sigma)_{\alpha j}\,a_\alpha\qquad(j\in I).
\end{equation}

We apply Lemma~\ref{lem:coefficient-change} to $\widetilde B_0=\widetilde B_{\setminus k}$, with new frozen coordinates $v=(t,y_1,\dots,y_r)$.
Let $R\in GL_{f+|N(k)|}(\bZ)$ be the matrix with columns $c_k,a_1,\dots,a_r$, and put
\[
 Q=R^{-1},\qquad
 P=\begin{pmatrix}\bigl((w_k)_j\bigr)_{j\in I}\\ 0\end{pmatrix}\in\Mat_{(1+r)\times I}(\bZ).
\]
The substitution $z^a\mapsto v^{Qa}$, $z_j\mapsto u_jt^{(w_k)_j}$ of Lemma~\ref{lem:coefficient-change} gives
\begin{equation}\label{eq:Sigma-variables}
 t=z^{c_k},\qquad y_\alpha=z^{a_\alpha},\qquad u_j=z_jt^{-(w_k)_j}\quad(j\in I).
\end{equation}

We compute the new frozen rows $PB_{I,I}+QC'$.
By \eqref{eq:frozen-splitting} and \eqref{eq:Sigma-matrix},
\[
 Qf_j=f_j(w_k)\,e_t+\sum_\alpha(C_\Sigma)_{\alpha j}\,e_\alpha .
\]
Since $\widetilde B_{kj}=0$ for $j\in I$, the $t$-row is
\[
 \sum_{i\in I}(w_k)_iB_{ij}+f_j(w_k)=\sum_{i\in M\sqcup F}\widetilde B_{ij}(w_k)_i=c_j(w_k)=0\qquad(j\in I).
\]
The $y$-rows are $C_\Sigma$, since $P$ has zero $y$-rows.
Hence
\[
 \binom{B_{I,I}}{PB_{I,I}+QC'}=\begin{pmatrix}B_{I,I}\\0\\C_\Sigma\end{pmatrix}.
\]
Since $Q\in GL_{f+|N(k)|}(\bZ)$, Lemma~\ref{lem:coefficient-change} gives
\[
 \sA(\widetilde B_{\setminus k})\cong\sA\begin{pmatrix}B_{I,I}\\0\\C_\Sigma\end{pmatrix}=\sA(\widetilde B_\Sigma)[t^{\pm1}],
\]
where the equality holds because $t$ appears in no exchange relation.
By \eqref{eq:Sigma-variables}, this isomorphism sends $z^{c_k}$ to $t$.
\end{proof}

\begin{corollary}
\label{cor:At-splitting} 
Assume \(B_{I, I}\) is locally acyclic. Then
 $\cA(\widetilde B)_t\cong\cN\times\Sigma$.
\end{corollary}
\begin{proof}
Since $B_{I,I}$ is locally acyclic, $\widetilde B_{N(k)}$ is a cluster localization by \eqref{eq:freezing-sandwich}.
So $\cA(\widetilde B)_t=\cA(\widetilde B_{N(k)})$, and Proposition~\ref{prop:AB_Nk} is equivalent to
$
\cA(\widetilde B)_t \cong \cA( \widetilde B_{\setminus k} ) \times_{\bC^*_t} \cN .
$
Hence the claim follows from Proposition~\ref{prop:splitting_ABprime}. 
\end{proof}

\subsection{Completion of \(\widetilde B_\Sigma\)}

Assume $\widetilde B$ is unimodular, with $E$, $\widetilde D$, $K$, $c_i$ and $w_i$ as in Section~\ref{sec:At-splitting}.  We will construct a unimodular completion of the matrix $\widetilde B_\Sigma$.

Put
\[
 N_\Sigma=\{v\in N_{\bZ}:\epsilon_k(v)=c_k(v)=0\},
\]
let $K_\Sigma$ be the restriction of $K$ to $N_\Sigma$, and let $K_k$ be the skew form on $\bZ e_k\oplus\bZ w_k$ with $K_k(w_k,e_k)=d_k$.

\begin{lemma}\label{lem:lattice-splitting}
Preserving the above notation, 
\begin{enumerate}[label=\textup{(\roman*)}]
\item There is a canonical splitting $N_{\bZ}=\bZ e_k\oplus\bZ w_k\oplus N_\Sigma$, with projection $v\mapsto v-\epsilon_k(v)e_k-c_k(v)w_k$ onto $N_\Sigma$.
\item $K(e_k,N_\Sigma)=K(w_k,N_\Sigma)=0$ and $K(w_k,e_k)=d_k$.
Hence $K=K_k\oplus K_\Sigma$.
\item For $j\in I$ the vectors $e_j$ and $w_j'=w_j-\epsilon_k(w_j)e_k$ lie in $N_\Sigma$, and
\begin{equation}\label{eq:residual-pairing}
 K_\Sigma(\cdot,e_j)=d_j\,c_j|_{N_\Sigma},\qquad K_\Sigma(\cdot,w_j')=-d_j\,\epsilon_j|_{N_\Sigma}.
\end{equation}
\item There are a basis of $N_\Sigma^*$ of the form $\{\epsilon_j|_{N_\Sigma}\}_{j\in I}\cup\{a_\alpha\}_\alpha$ with $a_\alpha(e_j)=0$, and positive rationals $q_\alpha$, such that
\[
K_\Sigma(N_\Sigma)\coloneqq\{K_\Sigma(\cdot,v):v\in N_\Sigma\}=\bigoplus_{j\in I}d_j\bZ\,\epsilon_j|_{N_\Sigma}\ \oplus\ \bigoplus_\alpha q_\alpha\bZ\,a_\alpha .
\]
Consequently, in the basis of $N_\Sigma$ dual to this one and with $D_\Sigma=\diag\bigl((d_j)_{j\in I},(q_\alpha)_\alpha\bigr)$, the matrix
\[
 E_\Sigma=K_\Sigma D_\Sigma^{-1}\in GL_{m+f-2}(\bZ)
\]
satisfies $E_\Sigma D_\Sigma=-D_\Sigma E_\Sigma^{\T}$ and $E_\Sigma e_j=c_j|_{N_\Sigma}$ for $j\in I$.
\end{enumerate}
\end{lemma}

\begin{proof}
(i) By \eqref{eq:dual-circle-identities},
\[
 \epsilon_k(e_k)=1,\qquad \epsilon_k(w_k)=0,\qquad c_k(e_k)=E_{kk}=0,\qquad c_k(w_k)=1.
\]
So the displayed map vanishes on $e_k$ and $w_k$, is the identity on $N_\Sigma$, and takes values in $N_\Sigma$.

(ii) By \eqref{eq:dual-circle-identities}, $K(\cdot,e_k)=d_kc_k$ and $K(\cdot,w_k)=-d_k\epsilon_k$.
Both vanish on $N_\Sigma$, and $K(w_k,e_k)=d_kc_k(w_k)=d_k$.

(iii) For $j\in I$, $\epsilon_k(e_j)=0$ and $c_k(e_j)=E_{jk}=0$.
Also $\epsilon_k(w_j')=0$ and $c_k(w_j')=c_k(w_j)-\epsilon_k(w_j)c_k(e_k)=0$.
Since $K(\cdot,e_k)=d_kc_k$ vanishes on $N_\Sigma$, \eqref{eq:residual-pairing} follows from \eqref{eq:dual-circle-identities}.

(iv) The vectors $e_j$, $j\in I$, are part of a basis of $N_\Sigma$, and the characters $\epsilon_j|_{N_\Sigma}$ are part of the dual basis.
Complete it by characters $a_\alpha$ with $a_\alpha(e_j)=0$.
By \eqref{eq:residual-pairing}, $d_j\epsilon_j|_{N_\Sigma}\in K_\Sigma(N_\Sigma)$.
For $v\in N_\Sigma$, the coefficient of $\epsilon_j|_{N_\Sigma}$ in $K_\Sigma(\cdot,v)$ is $K(e_j,v)=-d_jc_j(v)\in d_j\bZ$.
Hence
\[
 K_\Sigma(N_\Sigma)=\bigoplus_jd_j\bZ\epsilon_j|_{N_\Sigma}\oplus\Lambda,\qquad \Lambda=K_\Sigma(N_\Sigma)\cap\operatorname{span}_{\bQ}\{a_\alpha\}.
\]
Since $K_\Sigma$ is nondegenerate, $\Lambda$ has full rank in $\operatorname{span}_{\bQ}\{a_\alpha\}$.
Clear denominators and apply the Smith normal form.
This gives a new basis $\{a_\alpha\}$ of $\bigoplus_\alpha\bZ a_\alpha$ and $q_\alpha\in\bQ_{>0}$ with $\Lambda=\bigoplus_\alpha q_\alpha\bZ a_\alpha$.
In the dual bases, the matrix of $K_\Sigma$ is skew-symmetric, and $D_\Sigma^{-1}K_\Sigma$ is integral with image $\bZ^{m+f-2}$.
So $D_\Sigma^{-1}K_\Sigma$ is unimodular.
Therefore $E_\Sigma=K_\Sigma D_\Sigma^{-1}=-(D_\Sigma^{-1}K_\Sigma)^{\T}$ is unimodular, and
\[
 E_\Sigma D_\Sigma=K_\Sigma=-D_\Sigma E_\Sigma^{\T},\qquad E_\Sigma e_j=d_j^{-1}K_\Sigma(\cdot,e_j)=c_j|_{N_\Sigma}.\qedhere
\]
\end{proof}

For a square matrix $X$ we write $X^\vee\coloneqq-X^{\T}$. If we fix a completion \(E\) of \(\widetilde B\), we denote \( \widetilde B^\vee \coloneqq E^\vee_{\bullet, M}\). 
Then $(E^\vee,\widetilde D^{-1})$ is a completion of $\widetilde B^\vee$.
If $B^\vee$ is locally acyclic, later in Proposition~\ref{prop:finite-atlas} we will see $\cA(\widetilde B^\vee)\cong\cX(\widetilde B^\vee)$. In particular, $\cA(\widetilde B^\vee)$ does not depend on the completion, up to isomorphism.

\begin{proposition}\label{prop:Sigma-completion}
Keep the notations above. Then we have the following: 
\begin{enumerate}[label=\textup{(\roman*)}]
\item Take $\{a_\alpha\}$ in \eqref{eq:Sigma-variables} to be the basis of Lemma~\ref{lem:lattice-splitting}(iv), extended by zero on $\bZ e_k\oplus\bZ w_k$.
Then the character lattice of the seed $\widetilde B_\Sigma$ is $N_\Sigma^*$, with basis $\{\epsilon_j|_{N_\Sigma}\}_{j\in I}\cup\{a_\alpha|_{N_\Sigma}\}_\alpha$, and $(E_\Sigma,D_\Sigma)$ is a completion of $\widetilde B_\Sigma$.
So $\widetilde B_\Sigma$ is again unimodular.
\item The vertex $k$ is a sink of $B^\vee$, with the same set $I$.
There is a basis of the lattice $\Lambda_k$ of $\widetilde B^\vee$ for which $(\widetilde B^\vee)_\Sigma=(E_\Sigma^\vee)_{\bullet,I}$.
\end{enumerate}
\end{proposition}

\begin{proof}
(i) The variable $u_j$ is the monomial of the character $\epsilon_j-(w_k)_jc_k$.
This character vanishes on $e_k$ and $w_k$, because $c_k(e_k)=E_{kk}=0$ and $c_k(w_k)=1$.
It restricts to $\epsilon_j$ on $N_\Sigma$.
The variable $y_\alpha$ is the monomial of $a_\alpha$, which also vanishes on $e_k$ and $w_k$.
Together with $x=z^{\epsilon_k}$ and $t=z^{c_k}$, these characters form a basis of $N_{\bZ}^*$, by \eqref{eq:frozen-splitting}.
So the character lattice of the seed $\widetilde B_\Sigma$ is $N_\Sigma^*$, with basis $\{\epsilon_j|_{N_\Sigma}\}_{j\in I}\cup\{a_\alpha|_{N_\Sigma}\}$.
The lattice $\Lambda_k$ consists of the characters vanishing on $e_k$, $w_k$ and the $e_j$, $j\in I$.
So the basis of Lemma~\ref{lem:lattice-splitting}(iv), extended by zero, is a basis of $\Lambda_k$.
The exchange character of $\widetilde B_\Sigma$ at $j\in I$ is $c_j|_{N_\Sigma}$, and $c_j|_{N_\Sigma}=E_\Sigma e_j$ by Lemma~\ref{lem:lattice-splitting}(iv).

(ii) Since $E^\vee_{ij}=(d_j/d_i)E_{ij}$, the vertex $k$ is a sink of $B^\vee$ with the same neighbors.
The map $(E^\vee)^{-1}:N_{\bZ}^*\to N_{\bZ}$ sends $\epsilon_i\mapsto-w_i$ and $c^\vee_i=E^\vee\epsilon_i\mapsto e_i$.
So it identifies the character lattice of the dual seed with $N_{\bZ}$.
It sends the basis characters to the $-w_i$ and the exchange characters to the $e_i$.
Under this identification, the splitting of Lemma~\ref{lem:lattice-splitting} for the dual seed at $k$ removes $\bZ(-w_k)\oplus\bZ e_k$ and leaves $N_\Sigma$.
For $j\in I$, the projected basis characters and exchange characters are $-w_j'$ and $e_j$.
By \eqref{eq:residual-pairing},
\[
 K_\Sigma(\cdot,e_j)=d_jc_j|_{N_\Sigma},\qquad K_\Sigma(\cdot,-w_j')=d_j\epsilon_j|_{N_\Sigma}.
\]
Let $b_\alpha\in N_\Sigma$ satisfy $K_\Sigma(\cdot,b_\alpha)=q_\alpha a_\alpha$, with $a_\alpha$ and $q_\alpha$ as in Lemma~\ref{lem:lattice-splitting}(iv).
By Lemma~\ref{lem:lattice-splitting}(iv), $(v_i)=\{-w_j'\}_{j\in I}\cup\{b_\alpha\}$ is a basis of $N_\Sigma$.
The mutable cocharacters of the dual seed are the $-c_j$, so its lattice $\Lambda_k$ is $\{v\in N_\Sigma:c_j(v)=0,\ j\in I\}$.
This lattice has basis $\{b_\alpha\}$, because $c_j(v)=-d_j^{-1}K_\Sigma(e_j,v)$.
Since $K_\Sigma(\cdot,e_j)=d_jE_\Sigma e_j$, we get $e_j=\sum_i(D_\Sigma)_{jj}(D_\Sigma)_{ii}^{-1}(E_\Sigma)_{ij}v_i=\sum_i(E_\Sigma^\vee)_{ij}v_i$.
So $(\widetilde B^\vee)_\Sigma=(E_\Sigma^\vee)_{\bullet,I}$ for the basis $\{b_\alpha\}$ of $\Lambda_k$.
\end{proof}

\begin{corollary}\label{cor:dual-Sigma}
Assume $B^\vee_{I,I}$ is locally acyclic.
Then $\cA\bigl((\widetilde B^\vee)_\Sigma\bigr)\cong\cA\bigl((\widetilde B_\Sigma)^\vee\bigr)$, where $(\widetilde B_\Sigma)^\vee$ is taken for the completion $E_\Sigma$ of Proposition~\ref{prop:Sigma-completion}(i).
\end{corollary}
\begin{proof}
Follows from Proposition~\ref{prop:Sigma-completion}(ii), and Lemma~\ref{lem:coefficient-change}.
\end{proof}

We also have a more elementary fact:
\begin{lemma}\label{lem:dual-freezing}
The pair $(E,\widetilde D)$ is also a completion of $\widetilde B_{[k]}$ and of $\widetilde B_{N(k)}$.
For this completion, $(\widetilde B_{[k]})^\vee=(\widetilde B^\vee)_{[k]}$ and $(\widetilde B_{N(k)})^\vee=(\widetilde B^\vee)_{N(k)}$.
\end{lemma}
\begin{proof}
Since $E^\vee_{ij}=(d_j/d_i)E_{ij}$, the vertex $k$ has the same neighbors $N(k)$ in $B^\vee$.
The seed of freezing a set $J\subset M$ is $E_{\bullet,M\setminus J}$, and its dual seed is $E^\vee_{\bullet,M\setminus J}$.
Take $J=\{k\}$ and $J=N(k)$.
\end{proof}

\subsection{A finite quotient presentation for  the mirror stack}\label{sec:mirror-stack}

We give a more explicit description of \(\cX(\widetilde B^\vee)\).
Choose an integral matrix
\begin{equation}\label{alg:completion}
 T=\begin{pmatrix}-B^{\mathsf T}&-C^{\mathsf T}\\
 C'&T_{\rm ff}\end{pmatrix} \in \Mat_{m+f}(\bZ),
 \qquad
 \widetilde B_T^\vee=\binom{-B^{\mathsf T}}{C'},
\end{equation}
such that \(\det T \neq 0\).
Take the finite group \(G_T \coloneqq \Hom(\coker(T), \bC^*)\).

\begin{proposition}\label{prop:finite-atlas}
Assume that $\widetilde B$ has full rank.
There is a natural action of $G_T$ on $\cA(\widetilde B_T^\vee)$ and an isomorphism
\[
    \cX(\widetilde B^\vee)\cong[\cA(\widetilde B_T^\vee)/G_T].
\]
In particular, if $\widetilde B$ is RFR, then \(\cX(\widetilde B^\vee)\) is an affine scheme since we can choose  
$T\in GL_{m+f}(\bZ)$.
\end{proposition}
\begin{proof}
We use the notation $P^\vee$, $H$ and $\cX(\widetilde B^\vee)=[P^\vee/H]$ of~\eqref{eq:canonical-mirror}.
Let us first construct the action of $G_T$.
Let $\{e_i\}_{i\in M\sqcup F}$ be the standard basis of $\bZ^{M\sqcup F}$, and write $[c]$ for the class of $c$ in $\coker T$.
Define the action of $G_T$ on the initial coordinates $(u;v)$ of $\cA(\widetilde B_T^\vee)$ by
\[
g \cdot u_i = g([e_i])\, u_i, \qquad
g \cdot v_a = g([e_a])\, v_a.
\]
The exchange column of $\widetilde B_T^\vee$ at $j\in M$ is $Te_j$, and $[Te_j]=0$.
So both monomials of each exchange relation have the same weight under the action of $G_T$.
As with the action of $H$, the action therefore extends to $\sA(\widetilde B_T^\vee)$.

Next we build a finite cover of $P^\vee$.
Consider the homomorphism with weight $T$
\[
 \rho:(\bC^*)^{m+f}_{(h,v)}\longrightarrow (\bC^*)^{m+f}_{(Z,U)},\qquad
 Z_j=h^{-(B^{\mathsf T})_j}v^{C'_j},\qquad U_a=h^{-(C^{\mathsf T})_a}v^{(T_{\rm ff})_a}.
\]
Its kernel is $\Hom(\coker T,\bC^*)=G_T$.
An element $g\in G_T$ acts on $(\bC^*)^{m+f}_{(h,v)}$ by the translation
\[
 g\cdot h_i=g([e_i])\,h_i,\qquad g\cdot v_a=g([e_a])\,v_a .
\]

Apply Lemma~\ref{lem:coefficient-change} to the extended exchange matrix
$(-B^{\mathsf T};I_m;0)$ of $P^\vee$,
with $Q=T$ and $P=-(I_m;0)$.
The new matrix is
\[
 \begin{pmatrix}-B^{\mathsf T}\\ -\binom{I_m}{0}(-B^{\mathsf T})+T\binom{I_m}{0}\end{pmatrix}
 =\begin{pmatrix}-B^{\mathsf T}\\0\\C'\end{pmatrix}.
\]
Since the middle rows vanish,
\(\cA(-B^{\mathsf T};0;C') \cong H\times \cA(\widetilde B_T^\vee) \),
so we have a pull-back diagram
\[
\begin{tikzcd}
H\times \cA(\widetilde B_T^\vee)\arrow[r,"\pi"]\arrow[d]
\arrow[dr, phantom, "\lrcorner", very near start]
& P^\vee\arrow[d]\\
(\bC^*)^{m+f}_{(h,v)}\arrow[r,"\rho"] & (\bC^*)^{m+f}_{(Z,U)}
\end{tikzcd}
\]
where $\pi(h,u,v)=(X,Z,U)$ with $X_i=h_i^{-1}u_i$ and $(Z,U)=\rho(h,v)$, and the left vertical map is $(h,u,v)\mapsto(h,v)$.

It remains to take the quotients.
Since $\rho$ is a $G_T$-torsor, so is $\pi$, and its action is
\[
 g\cdot(h_i,u_i,v_a)=\bigl(g([e_i])\,h_i,\ g([e_i])\,u_i,\ g([e_a])\,v_a\bigr).
\]
Let $s\in H$ act on $H\times \cA(\widetilde B_T^\vee)$ by $s\cdot(h,u,v)=(s^{-1}h,u,v)$.
Then $\pi(s^{-1}h,u,v)=s\cdot\pi(h,u,v)$, because $X_i$, $Z_j$, $U_a$ are multiplied by $s_i$, $s^{(B^{\mathsf T})_j}$, $s^{(C^{\mathsf T})_a}$, which are the weights~\eqref{eq:weights}.
So $\pi$ is $H$-equivariant.
The actions of $H$ and $G_T$ on $H\times \cA(\widetilde B_T^\vee)$ commute.
Hence
\[
 \cX(\widetilde B^\vee)=[P^\vee/H]\cong\bigl[(H\times \cA(\widetilde B_T^\vee))/(H\times G_T)\bigr].
\]
Finally, $H$ acts freely and transitively on the first factor, and the projection $H\times \cA(\widetilde B_T^\vee)\to \cA(\widetilde B_T^\vee)$ is $G_T$-equivariant.
Hence
\[
 \bigl[(H\times \cA(\widetilde B_T^\vee))/(H\times G_T)\bigr]\cong[\cA(\widetilde B_T^\vee)/G_T].
\]
\end{proof}

\begin{remark}
    Proposition~\ref{prop:finite-atlas} implies that for different choices of \(T\), the stack \([ \cA(\widetilde B_T/G_T]\) are all isomorphic to each other. 
    A natural choice of \(T\) is:
    \[
    T =
    \begin{pmatrix}
        -B^\T & -C^\T \\
        CD & I_f
    \end{pmatrix}
    \]
    One can easily verify that \(T\) is non-degenerate. In particular, if \(B\) is skew-symmetric, \(\widetilde B_T^\vee \cong \widetilde B\), and \(\cX (\widetilde B^\vee) \cong [\cA(\widetilde B)/G_T]\) 
\end{remark}

Here we give a class of examples where \(\widetilde B\) is unimodular and \(B\) is skew-symmetric, but \(\cA(\widetilde B)\) is not self-mirror. 

\begin{example}\label{ex:rfr-not-self-mirror}
Fix $p\ge2$. For $c\ge1$, put
\[
\widetilde B_c=\begin{pmatrix}0&p\\-p&0\\1&0\\0&c\end{pmatrix}.
\]
This matrix is acyclic and \(\widetilde B\) is unimodular precisely when $\gcd(p,c)=1$. If $\gcd(p,c)=1$, 
take \(a\) such that \(ac\equiv 1 \pmod p\). Set \(r=(ac-1)/p\) and
\[
E=
\begin{pmatrix}
0&p&-a&0\\
-p&0&0&-ac\\
1&0&0&r\\
0&c&-r&0
\end{pmatrix}
\qquad
\widetilde D=\operatorname{diag}(1,1,a^{-1},a^{-1}),
\]
Then \(E\) is a completion of \(\widetilde B_c\), \(E \widetilde D\) is skew-symmetric, and \(\det E=(pr-ac)^2=1\). 

By~\cite[Corollary~1.21]{BFZ-III},
\[
\cA(\widetilde B_c)=\{xx'=u+y^p,\quad yy'=1+v^cx^p,\quad u,v\ne0\}.
\]
Assume $\gcd(p,c)=1$, and choose $d\ge1$ with $cd\equiv1\pmod p$.
Choose $r,s\in\bZ$ such that $pr-ds=(1-cd)/p$. Then take
\[
T=\begin{pmatrix}
0&p&-1&0\\
-p&0&0&-c\\
1&0&1&s\\
0&d&0&r
\end{pmatrix},\qquad \det T=p^2r-pds+cd=1.
\]
Proposition~\ref{prop:finite-atlas} therefore gives
\[
\cX(\widetilde B_c^\vee)\cong\cA(\widetilde B_d).
\]

For $c,e$ coprime to $p$, using the results in~\cite{BlancDolgachev} one can show
\[
\cA(\widetilde B_c)\cong\cA(\widetilde B_e)\quad\Longleftrightarrow\quad c\equiv\pm e\pmod p.
\]
Consequently, $\cA(\widetilde B_c)$ is self-mirror if and only if $c^2\equiv\pm1\pmod p$.

\end{example}

\subsection{Comparison with the Fock-Goncharov $\mathcal X$-variety}

The cluster $\mathcal X$-variety was introduced in \cite{FG09}.  We follow the presentation of  Gross--Hacking--Keel~\cite[Section~2]{GHK}.

The \emph{fixed data} $\Gamma$ consist of
\begin{itemize}
\item a lattice $N$ with a skew-symmetric form $\{\cdot,\cdot\}:N\times N\to\bQ$, and a saturated sublattice $N_{\mathrm{uf}}\subset N$;
\item an index set $I$ with $|I|=\rank N$, a subset $I_{\mathrm{uf}}\subset I$ with $|I_{\mathrm{uf}}|=\rank N_{\mathrm{uf}}$, and positive integers $d_i$, $i\in I$, with greatest common divisor $1$;
\item a sublattice $N^\circ\subset N$ of finite index, subject to integrality conditions, and $M=\Hom(N,\bZ)$, $M^\circ=\Hom(N^\circ,\bZ)$.
\end{itemize}
A \emph{seed} is a basis $\mathbf s=(e_i)_{i\in I}$ of $N$ such that $(e_i)_{i\in I_{\mathrm{uf}}}$ is a basis of $N_{\mathrm{uf}}$ and $(d_ie_i)_{i\in I}$ is a basis of $N^\circ$.
It determines
\[
 \epsilon_{ij}=\{e_i,e_j\}d_j,\qquad f_i=d_i^{-1}e_i^*\in M^\circ,
\]
and two seed tori with coordinates
\[
 \cX^{\mathrm{tori}}_{\mathbf s}=\Spec\bC[N],\quad X_i=z^{e_i},\qquad
 \cA^{\mathrm{tori}}_{\mathbf s}=\Spec\bC[M^\circ],\quad A_i=z^{f_i}.
\]
For $k\in I_{\mathrm{uf}}$, the mutation $\mu_k$ is the birational map given by \cite[(2.7), (2.8)]{GHK}
\[
 \mu_k^*X_i=\begin{cases}X_k^{-1},&i=k,\\ X_i\bigl(1+X_k^{-\operatorname{sgn}\epsilon_{ik}}\bigr)^{-\epsilon_{ik}},&i\ne k,\end{cases}
 \qquad
 A_k\,\mu_k^*A_k=\prod_{\epsilon_{kj}>0}A_j^{\epsilon_{kj}}+\prod_{\epsilon_{kj}<0}A_j^{-\epsilon_{kj}},\quad \mu_k^*A_i=A_i\ (i\ne k).
\]
The cluster varieties in~\cite{GHK} are defined by gluing the seed tori of all seeds along these maps.
We denote them by  $\cX^{\mathrm{tori}}$ and $\cA^{\mathrm{tori}}$. 
Every coordinate of a seed torus is invertible, including the frozen coordinates $X_i,A_i$ with $i\notin I_{\mathrm{uf}}$.
The ring $\Gamma(\cA^{\mathrm{tori}},\cO_{\cA^{\mathrm{tori}}})$ is the upper cluster algebra \cite[Definition~2.8, Remark~2.9]{GHK}, and the canonical map $\cA^{\mathrm{tori}}\to\Spec\Gamma(\cA^{\mathrm{tori}},\cO_{\cA^{\mathrm{tori}}})$ is an open immersion \cite[Theorem~3.14]{GHK}. 

For principal coefficients, put $\widetilde N=N\oplus M^\circ$ with
\[
 \{(n_1,m_1),(n_2,m_2)\}=\{n_1,n_2\}+\langle n_1,m_2\rangle-\langle n_2,m_1\rangle,
\]
unfrozen part $N_{\mathrm{uf}}\oplus0$, and the seed $\{(e_i,0)\}\cup\{(0,f_i)\}$.
The resulting $\cA$-variety is denoted by $\cA^{\mathrm{tori}}_{\mathrm{prin}}$.
For every seed there is an exact sequence of tori $1\to T_{N^\circ}\to\cA^{\mathrm{tori}}_{\mathrm{prin},\mathbf s}\to\cX^{\mathrm{tori}}_{\mathbf s}\to1$ compatible with mutation.
So $\cA^{\mathrm{tori}}_{\mathrm{prin}}\to\cX^{\mathrm{tori}}$ is a $T_{N^\circ}$-torsor \cite[Construction~2.11]{GHK}.

Our matrix $\widetilde B=\binom BC$ gives fixed data with $I=M\sqcup F$, $I_{\mathrm{uf}}=M$, the $d_i$ of $\widetilde D$, and
\[
 \epsilon_{ki}=\widetilde B_{ik}\qquad(k\in M,\ i\in M\sqcup F),
\]
because then \cite[(2.8)]{GHK} is the exchange relation \eqref{eq:exchange}.
The symmetrizer condition $d_i\epsilon_{ij}=-d_j\epsilon_{ji}$ is \eqref{eq:symmetrizer}.
Let $\Gamma^\vee$ be the Langlands dual fixed data of \cite[Appendix~A]{GHKK}.
It has $d_i^\vee=d_i^{-1}\delta$ with $\delta=\operatorname{lcm}(d_i)$, the form $\{\cdot,\cdot\}^\vee=\delta^{-1}\{\cdot,\cdot\}$, and the seeds $\mathbf s^\vee=(d_ie_i)_{i\in I}$.
So its matrix is
\[
 \epsilon^\vee_{ij}=\delta^{-1}\{d_ie_i,d_je_j\}\,d_j^{-1}\delta=d_i\{e_i,e_j\}=-\epsilon_{ji},
\]
where the last equality is $d_i\epsilon_{ij}=-d_j\epsilon_{ji}$.
We put
\[
 \cX^{\mathrm{tori}}(\widetilde B^\vee)=\cX^{\mathrm{tori}}_{\Gamma^\vee},
\]
which is the Fock--Goncharov dual $(\cA^{\mathrm{tori}})^\vee$ of $\cA^{\mathrm{tori}}$ \cite[Definition~A.4]{GHKK}.
If $B$ is skew-symmetric, then $\epsilon^\vee=\epsilon$, and $\cX^{\mathrm{tori}}(\widetilde B^\vee)$ is the usual $\cX$-variety of $\widetilde B$.
Since $\epsilon^\vee_{ik}=-\widetilde B_{ik}$, \cite[(2.7)]{GHK} gives the transition maps of $\cX^{\mathrm{tori}}(\widetilde B^\vee)$ at $k\in M$:
\begin{equation}\label{eq:cluster-x-mutation}
 q_k'=q_k^{-1},\qquad
 q_i'=q_i\,q_k^{[-\widetilde B_{ik}]_+}(1+q_k)^{\widetilde B_{ik}}
 \quad(i\in M\sqcup F,\ i\ne k).
\end{equation}

The mirror $\cX(\widetilde B^\vee)=[P^\vee/H]$ uses the affine scheme $P^\vee=\cA_{\mathrm{prin}}(-B^{\T})\times(\bC^*)^f$.
Its torus-glued version is
\[
 P^{\vee,\mathrm{tori}}=\cA^{\mathrm{tori}}_{\mathrm{prin}}(-B^{\T})\times(\bC^*)^f,
\]
where $\cA^{\mathrm{tori}}_{\mathrm{prin}}(-B^{\T})$ is the union of the seed tori of $\binom{-B^{\T}}{I_m}$.
The torus $H$ acts on it with the weights \eqref{eq:weights}.

\begin{proposition}\label{prop:cluster-x-mirror}
Suppose that $\widetilde B$ has full rank and that $-B^{\T}$ is locally acyclic.
There are an identification and a dense open immersion
\begin{equation}\label{eq:cluster-x-open}
 \cX^{\mathrm{tori}}(\widetilde B^\vee)\cong[P^{\vee,\mathrm{tori}}/H]\lhook\joinrel\longrightarrow\cX(\widetilde B^\vee).
\end{equation}
Restriction gives an isomorphism $\Gamma(\cX(\widetilde B^\vee),\cO)\xrightarrow{\sim}\Gamma(\cX^{\mathrm{tori}}(\widetilde B^\vee),\cO)$.
In particular, the coarse space of $\cX(\widetilde B^\vee)$ is the affinization of $\cX^{\mathrm{tori}}(\widetilde B^\vee)$:
\begin{equation}
\label{eq:cluster-x-affinization}
 \cX(\widetilde B^\vee)_{\mathrm{coarse}}\cong\Spec\Gamma(\cX^{\mathrm{tori}}(\widetilde B^\vee),\cO).
\end{equation}
\end{proposition}

\begin{proof}
We first identify $\cX^{\mathrm{tori}}(\widetilde B^\vee)$ with the quotient.
Apply the torsor $\cA^{\mathrm{tori}}_{\mathrm{prin}}\to\cX^{\mathrm{tori}}$ of \cite[Construction~2.11]{GHK} to $\Gamma^\vee$.
Its structure torus has one factor for each index in $M\sqcup F$, and we split it as $H\times(\bC^*)^f$.
The factor $(\bC^*)^f$ acts on the frozen coordinates of $\cA^{\mathrm{tori}}_{\mathrm{prin}}$ with their standard weights, and these coordinates are invertible.
So the quotient by $(\bC^*)^f$ is the slice where the frozen $\cA$-coordinates equal $1$.
On this slice, the mutable coordinates and their principal coefficients have the exchange matrix $\binom{-B^{\T}}{I_m}$.
The principal coefficients of the frozen indices remain, as the invertible coordinates $U_a$.
The residual action of $H$ has the weights \eqref{eq:weights}.
Hence $\cX^{\mathrm{tori}}(\widetilde B^\vee)\cong[P^{\vee,\mathrm{tori}}/H]$.

On the initial chart the quotient coordinates are
\begin{equation}\label{eq:cluster-x-invariants}
 q_i=Z_i\prod_{j\in M}X_j^{-B_{ij}}\quad(i\in M),
 \qquad
 q_a=U_a\prod_{j\in M}X_j^{-C_{aj}}\quad(a\in F).
\end{equation}
They are invariant under \eqref{eq:weights} and generate the invariant Laurent ring, and the $X_j$ trivialize the $H$-torsor on this chart.
The exchange relation at $k$ reads
\[
 X_kX_k'
 =Z_k\prod_jX_j^{[-B_{kj}]_+}+\prod_jX_j^{[B_{kj}]_+}
 =\Bigl(\prod_jX_j^{[B_{kj}]_+}\Bigr)(1+q_k).
\]
Substituting it into \eqref{eq:cluster-x-invariants} for the mutated seed gives \eqref{eq:cluster-x-mutation}, on the mutable and on the frozen indices.

Next we prove the open immersion.
Since $-B^{\T}$ is locally acyclic, \cite{M} gives
\[
 \sA\binom{-B^{\T}}{I_m}=\sU\binom{-B^{\T}}{I_m}.
\]
By \cite[Theorem~3.14]{GHK}, $\cA^{\mathrm{tori}}_{\mathrm{prin}}(-B^{\T})$ is a dense open subset of the spectrum of this algebra, which is $\cA_{\mathrm{prin}}(-B^{\T})$.
So $P^{\vee,\mathrm{tori}}\subset P^\vee$ is a dense open $H$-invariant subset, and \eqref{eq:cluster-x-open} follows.

It remains to compare global functions.
Put $S=\bC[U_a^{\pm1}:a\in F]$.
The global functions on $\cA^{\mathrm{tori}}_{\mathrm{prin}}(-B^{\T})$ form its upper cluster algebra \cite[Definition~2.8, Remark~2.9]{GHK}.
Taking $H$-invariants gives
\[
 \Gamma(\cX^{\mathrm{tori}}(\widetilde B^\vee),\cO)
 =\Bigl(\sU\binom{-B^{\T}}{I_m}\otimes S\Bigr)^H
 =\Bigl(\sA\binom{-B^{\T}}{I_m}\otimes S\Bigr)^H
 =\Gamma(\cX(\widetilde B^\vee),\cO).
\]
By Proposition~\ref{prop:finite-atlas}, $\cX(\widetilde B^\vee)$ is the quotient of an affine scheme by a finite group.
So its coarse space is affine with coordinate ring $\Gamma(\cX(\widetilde B^\vee),\cO)$, which proves \eqref{eq:cluster-x-affinization}.
\end{proof}

So $\cX(\widetilde B^\vee)$ contains the seed-glued $\cX^{\mathrm{tori}}(\widetilde B^\vee)$ as a dense open substack, and its coarse space is the affinization of $\cX^{\mathrm{tori}}(\widetilde B^\vee)$.
Along the added locus the stack can have finite stabilizers.
If $\widetilde B$ is really full rank, Proposition~\ref{prop:finite-atlas} identifies $\cX(\widetilde B^\vee)$ with the affine scheme \eqref{eq:cluster-x-affinization}.

\begin{example}
Let $m=f=1$ and $\widetilde B=\binom0d$ with $d\ge1$.
Then
\[
\cA(\widetilde B) = \{ xy = 1+t^d\} \subset \bC_x \times \bC_y \times \bC_t^*, \qquad
\mathcal A^{\mathrm{tori}}(\widetilde B)
=
\mathcal A(\widetilde B)
\setminus
\bigl\{(0,0,\zeta):\zeta^d=-1\bigr\}.
\]
The two torus charts of \(\cA^{\mathrm{tori}}\) are \( \{x\neq 0 \}\) and \(\{y\neq 0\}\), and their union misses \(d\) points. 

On the mirror side,
\[
 P^\vee=\{XX'=1+Z,\ Z\ne0\} \times \bC_U^*,\qquad H = \bC^* \text{ acts by } \quad
 h\cdot(X,X',Z,U)=(hX,h^{-1}X',Z,h^dU).
\]
The slice $U=1$ has the residual group $\mu_d$, so
\[
 \cX(\widetilde B^\vee)\cong\bigl[\{ab=1+z,\ z\ne0\}/\mu_d\bigr],\qquad
 \zeta\cdot(a,b,z)=(\zeta a,\zeta^{-1}b,z).
\]
With $u=a^d$ and $v=b^d$, its coarse space is
\[
 \Spec\bC[u,v,z^{\pm1}]/\bigl(uv-(1+z)^d\bigr).
\]
The two torus charts of $\cX^{\mathrm{tori}}(\widetilde B^\vee)$ are $\{u\ne0\}$ and $\{v\ne0\}$, and their union misses the stacky point $(u,v,z)=(0,0,-1)$ with stabilizer $\mu_d$.
\end{example}

\section{Liouville geometry}
\label{sec:geometry}

Let \( \widetilde B = \binom{B}{C} \) be unimodular and locally acyclic, and fix a  completion $(E,\widetilde D)$ of $\widetilde B$, such that \(E\in GL_{m+f}(\bZ) \) and $K=E\widetilde D$ is skew-symmetric. 
On a seed torus write $z_i=\exp(u_i+i\theta_i)$, and consider the Weil-Petersson form~\cite{GSV05}
\begin{equation}\label{eq:forms}
 \Omega=\frac12\sum_{i,j}K_{ij}\,\dd\log z_i\wedge \dd\log z_j,\qquad
 \ltor=u^{\T}K\,\dd\theta=\sum_{i,j}K_{ij}u_i\,\dd\theta_j.
\end{equation}
Then \(  \Im\Omega= \dd\ltor \). 
There is an identification
\begin{align}
\label{eq:identify_with_torus}
   (\bC^*)^{m+f}  & \longrightarrow  T^*(S^1)^{m+f} \\
    (z_j)_j & \longmapsto (\varphi, p) = ( -E^\T \theta, \widetilde Du) \notag
\end{align}
where \(\varphi\) and \(p\) are the base and fiber coordinates. Then \(\lambda_{\mathrm{tor}}\) is identified with \( p \dd\varphi\).  

Recall the cocharacters $w_k=(E^{\T})^{-1}e_k\in N_{\bZ}$ of \eqref{eq:dual-circles}.
 The cocharacter $w_k$ gives a circle action on the seed torus by
\begin{equation}
\label{eq:circle_action}
e^{i\tau}\cdot z_j
=
e^{i\tau(w_k)_j}z_j.
\end{equation}
If \(k \in F\) is frozen, then the action extends to the cluster variety \(\cA(\widetilde B)\). 
Indeed, \eqref{eq:dual-circle-identities} gives $c_j(w_k)=0$ for $j\in M$, so both monomials of every exchange relation \eqref{eq:exchange} have the same weight.

In Section~\ref{sec:cluster-liouville} we prove
\begin{proposition}[{See also~\cite[Theorem 4.4]{M}}]
    \label{prop:Omega}
    The \(2\)-form \(\Omega\) extends to a holomorphic symplectic form on \(X = \cA(\widetilde B)\).
\end{proposition}

We write $L_+=\{z_i>0\}$ for the Lagrangian of positive real locus. For Morse--Bott Weinstein structure we refer to Section~\ref{app:Morse-Bott}.

\begin{definition}
    We say a Liouville structure \(\lambda\) on a cluster variety \(X = \cA(\widetilde B)\) is \emph{good} if it is Morse--Bott with a Lyapunov function \(\Psi\), and satisfies the following
\begin{enumerate}[label=\textup{(P\arabic*)}, series=P]
\item
\label{P1} 
\(\dd\lambda = \Im \Omega\), and for every seed torus $T$, the form $(\lambda-\ltor)|_T$ is exact.

\item
\label{P2} 
For every frozen index $a$, $Z_{\lambda} \log|z_a|=\log|z_a|$, and $\lambda$ is invariant under the torus generated by the circles $w_a$, $a\in F$.

\item 
\label{P3}
$\lambda|_{L_+}= 0$.

\item
\label{P4} The Lyapunov  function $\Psi$ is invariant under the torus of \ref{P2}, and its critical components are isotropic tori.

\item
\label{P5}
For every seed, there are a cotangent disk bundle $D^*_\epsilon\T^{m+f}$ of its seed torus, and a compact exact deformation $\lambda_s$ on \(X\) such that
\[
 \lambda_0=\lambda,\qquad \lambda_1|_{D^*_\epsilon\T^{m+f}}=p\,\dd\varphi,\qquad \lambda_s|_{L_+}=0 .
\]
Here $p\,\dd\varphi$ is the Liouville form of \eqref{eq:identify_with_torus}.

    \end{enumerate}
\end{definition}

In this section we prove the following proposition.

\begin{proposition}\label{prop:lambda_ad} 
With the above notation,  for each sink-recurrent $X=\cA(\widetilde B)$, we can choose a good Liouville structure $\lad$ (\( \rm{ad} \) stands for ``adapted'') on \(X\), such that 
if $k$ is a sink of $B$ with $B_{M\setminus\{k\},M\setminus\{k\}}$ and $B_{I,I}$ both sink-recurrent, then $\lad$ has following additional properties:
\begin{enumerate}[resume*=P]
\item
\label{P6} 
At the sink $k$, write $V=\cA(\widetilde B)_x$ and $\cA(\widetilde B)_t\cong\cN\times\Sigma$.
Let $\lambda_V$ and $\lambda_\Sigma$ be the adapted Liouville forms on $V$ and $\Sigma$.
There is an exact symplectomorphism
\[
 (V,\lad|_V)\cong(V,\lambda_{V})
\]
equal to the identity on a Liouville domain of $V$ and preserving $L_+\cap V$.
Take the stop
\[
 \cF=T^*S^1_\theta\times\Sigma,\qquad \lambda_{\cF}=d_k\log|t|\,\dd\theta+\lambda_\Sigma,
 \qquad \theta=\arg x
\]
lying in the end $\{\log|x|\to-\infty\}$ of $\partial_\infty V$ at $\arg(-t)=0$.
Then $(X,\lad)$ is the completion of the sectorial gluing
\[
 (V,\cF)\underset{\cF}{\cup}(T^*D^2\times\Sigma).
\]

\item\label{P7}
At the sink $k$, keep the notation of \ref{P6}, and put
\[
 O=\{x\ne0,\ t\ne0\}\cong\bC^*_x\times\bC^*_t\times\Sigma,\qquad
 \lambda_O=d_k\lN+\lambda_\Sigma .
\]
There are Liouville domains $D_O\subset D_V$ with $\widehat D_O\cong(O,\lambda_O)$ and $\widehat D_V\cong(V,\lambda_V)$, whose contact boundaries coincide near $\cF$.
The circle $w_k$ acts on $V$ and preserves $\lambda_V$.
On $O$ it acts by $e^{i\alpha}\cdot(x,t,\sigma)=(x,e^{i\alpha}t,\sigma)$.
For $\alpha\in(0,2\pi)$, $e^{i\alpha}\cF\cap\cF=\emptyset$, and $\alpha\mapsto e^{i\alpha}\cF$ is a positive isotopy.
So $D_O\subset D_V$ and $\cF$ satisfy the geometric hypotheses of Lemma~\ref{lem:circle-twist} and Theorem~\ref{thm: SOV descent}.
\end{enumerate}
\end{proposition}

Proposition~\ref{prop:lambda_ad} is proved in Section~\ref{sec:cluster-liouville}.

For principal coefficients \(\binom{B}{I_m}\), take the completion
\begin{equation}\label{eq:principal-completion}
 E=\begin{pmatrix}B&-I_m\\I_m&0\end{pmatrix},\qquad \widetilde D=\diag(D,D).
\end{equation}

\begin{corollary}\label{cor:principal-geometry}
If $B$ is sink recurrent, then $\cA_{\mathrm{prin}}(B)=\cA\binom B{I_m}$ carries a good Liouville structure $\lad$, with Lyapunov function $\Psi$, for the completion \eqref{eq:principal-completion}.
\end{corollary}

In the full rank case, we do not choose a completion of $\widetilde B$.
Instead, we realize $\cA(\widetilde B)$ as a fiber of a map on a principal model.
A finite cover of the principal model is a product of $\cA(\widetilde B)$ with a cotangent torus.

Let $\widetilde B=\binom BC$ have full rank, with $B$ sink recurrent.
Rescaling $D$ by $c\in\bQ_{>0}$ replaces $\lad$ by $c\lad$, which has the same Liouville field.
So we may assume that $D^{-1}$ is integral.
Let $(\lad,\Psi)$ be the structure of Corollary~\ref{cor:principal-geometry} on $\cA_{\mathrm{prin}}(B)$ for this $D$, with positive locus $L_+$.
Put
\[
 P=\cA_{\mathrm{prin}}(B)\times(\bC^*)^f,\qquad
 \lambda_P=\lad+\sum_{a=1}^f\log|U_a|\,\dd(\arg U_a) .
\]
We again write $L_+=\{X_i,Z_i,U_a>0\}$ for the positive locus of $P$.
On the initial torus, in the coordinates $(X;Z;U)$ with $X_i=\exp({u_{X,i}+i\theta_{X,i}})$, $Z_i=\exp({u_{Z,i}+i\theta_{Z,i}})$, and $U_a=\exp({u_{U,a}+i\theta_{U,a}})$, 
\begin{equation}\label{eq:principal-form}
 \omega_P
 = \dd u^{\T} \wedge 
 \begin{pmatrix}
 BD&-D&0\\D&0&0\\0&0&I_f
 \end{pmatrix}
 \wedge \dd\theta, \qquad 
 Z_P\log|Z_i|=\log|Z_i|,\qquad
 Z_P\log|U_a|=\log|U_a| .
\end{equation}
Define
\begin{equation}\label{eq:frozen-map}
 F:P\longrightarrow(\bC^*)^m,\qquad F_i=Z_i\prod_a U_a^{-C_{ai}},\qquad
 S=D^{-1}B+C^{\T}C\in\Mat_{m\times m}(\bZ).
\end{equation}
The matrix $S$ is nondegenerate.
Indeed, $(D^{-1}B)^{\T}=-D^{-1}B$ by \eqref{eq:symmetrizer}, so $v^{\T}Sv=\norm{Cv}^2$ for $v\in\bR^m$.
Hence $Sv=0$ gives $Cv=0$, then $D^{-1}Bv=0$, and then $v=0$ because $\widetilde B$ has full rank. 

Let  \(H = (\bC^*)^m \) act on \(P\) with the weight matrix \( [ D^{-1}, (D^{-1} B)^T, -C^T]\).  Then
we have
\(
F(h\cdot z)=h^S F(z)
\).

\begin{proposition}\label{prop:fullrank-geometry}
Let $\widetilde B=\binom BC$ have full rank, with $B$ sink recurrent.
\begin{enumerate}[label=\textup{(Q\arabic*)}]
\item\label{Q1} $(P,\lambda_P)$ is a Liouville manifold with a Morse--Bott Lyapunov function.
\item\label{Q2} The fiber $F^{-1}(1)$ is isomorphic to $\cA(\widetilde B)$.
\item\label{Q3} $(\cA(\widetilde B),\alpha)$ is a Liouville manifold with a Morse--Bott Lyapunov function, where $\alpha=\lambda_P|_{\cA(\widetilde B)}$.
The primitive $\alpha$ vanishes on $L_{\widetilde B,+}=L_+\cap \cA(\widetilde B)$.

\item\label{Q4} Let $\Gamma=\bZ^m$ act on $\cA(\widetilde B)\times T^*\bR^m$ by
\[
 g\cdot(x,q,p)=\bigl(e^{-2\pi iS^{-1}g}\cdot x,\ q+2\pi S^{-1}g,\ p\bigr).
\]
There are a compact exact deformation $\lambda_P+\dd H_\tau$, $0\le\tau\le1$, with $H_0=0$ and $H_\tau|_{L_+}=0$, and a $\Gamma$-covering map $\Pi\colon\cA(\widetilde B)\times T^*\bR^m\to P$ with
\[
 \Pi^*(\lambda_P+\dd H_1)=\alpha+p\,\dd q,\qquad
 \Pi^{-1}(L_+)=\Gamma\cdot\bigl(L_{\widetilde B,+}\times T^*_0\bR^m\bigr).
\]
\end{enumerate}
\end{proposition}

Proposition~\ref{prop:fullrank-geometry} is proved in Section~\ref{sec:fullrank-geometry}.

\subsection{Review on Weinstein structures}\label{app:Morse-Bott}

For a one-form $\lambda$ on $X$ with $\omega=\dd\lambda$ symplectic, the Liouville field $Z=Z_\lambda$ is defined by $\iota_Z\omega=\lambda$, and $\phi_Z^\tau$ is its flow.

In this paper a \emph{Liouville manifold} $(X,\lambda)$ is always complete and of finite type:
$Z$ is complete, and there is a compact codimension-zero submanifold $D\subset X$ with boundary such that $Z$ points outward along $\partial D$ and $X=\bigcup_{\tau\ge0}\phi_Z^\tau(D)$.
Such a $D$ is called a \emph{Liouville domain} for $X$.
The \emph{completion} of a Liouville domain $(D,\lambda)$ is
\[
 \widehat D=D\cup_{\partial D}\bigl([0,\infty)\times\partial D\bigr),\qquad \lambda=e^r\,\lambda|_{\partial D}\ \text{ on }[0,\infty)\times\partial D,
\]
where $r$ is the coordinate on $[0,\infty)$.
If $D$ is a Liouville domain for \(X\), then $(r,y)\mapsto\phi_Z^r(y)$ identifies $\widehat D$ with $(X,\lambda)$.
For a complete Liouville field $Z$, the \emph{skeleton} $\Skel(Z)$ is the set of points whose forward orbit is relatively compact.
For a Liouville domain $D$,
$\Skel(Z)=\bigcap_{\tau\ge0}\phi_Z^{-\tau}(D)$, which is compact.
A \emph{compact exact deformation} is a family $\lambda+\dd H_\tau$, $0\le\tau\le1$, with $\bigcup_\tau\supp H_\tau$ compact; it keeps $(X,\lambda)$ a Liouville manifold.
A Lagrangian $L$ is \emph{strict exact} if $\lambda|_L=0$; then $Z|_L\in TL$.
We write $A\Subset X$ if $A$ has compact closure (i.e. relatively compact) in $X$.

We recall Morse--Bott Weinstein structures, following \cite{Eli18,AGEN}.
Let $(X^{2n},\lambda)$ be a Liouville manifold with Liouville field $Z$.
At a zero $z$ of $Z$, the linearization $D_zZ$ splits $T_zX=E^+_z\oplus E^0_z\oplus E^-_z$ into generalized eigenspaces, with eigenvalues of positive, zero and negative real part.

\begin{definition}\label{def:quantitative}
A \emph{Morse--Bott Weinstein structure} on $X^{2n}$ is a pair $(\lambda,\Psi)$, where $(X,\lambda)$ is a Liouville manifold and $\Psi$ is a proper, bounded-below Morse--Bott function such that
\begin{enumerate}[label=\textup{(\roman*)}]
\item $\Crit\Psi=\Zero(Z_\lambda)$, and $\dd\Psi(Z_\lambda)>0$ on $X\setminus\Crit\Psi$;
\item $\Crit\Psi$ is a finite union of compact isotropic submanifolds $C$;
\item for some metric and some $c>0$, $\dd\Psi(Z_\lambda)\ge c\bigl(|Z_\lambda|^2+|\dd\Psi|^2\bigr)$ near $\Crit\Psi$;
\item $\ind_C\Psi+\dim C\le n$ for every component $C$, where $\ind_C\Psi$ is the normal Morse--Bott index.
\end{enumerate}
Then $\Psi$ is a \emph{Morse--Bott Lyapunov function} for $\lambda$.
If moreover $\Psi$ is a Morse function, then $(\lambda,\Psi)$ is a \emph{Weinstein structure} \cite[Chapter~11]{CE}.
\end{definition}

At a point of a component $C$ condition (iii) gives
\[
 \Hess\Psi(v,DZv)\ge c\bigl(|DZv|^2+|\Hess\Psi(v,\cdot)|^2\bigr)\qquad(v\in TX).
\]
Hence $\ker DZ=TC$, and $DZ$ has no nonzero purely imaginary eigenvalues, so $E^0|_C=TC$.

A Morse--Bott Weinstein structure becomes a Weinstein structure after a homotopy supported near $\Crit\Psi$ \cite[Lemma~2.18]{AGEN}, \cite[Proposition~2.7]{St}.
We apply results stated for Weinstein manifolds, such as those of \cite{GPS2}, directly to Morse--Bott Weinstein structures.

\subsection{A local model}\label{sec:local-models}

We first describe a local model.\label{sec:node-model}
Let
\[
 \cN=\{(x,y,t)\in\bC^2\times\bC^*:\ xy=1+t\}.
\]
This is  the cluster variety with quiver 
\begin{tikzpicture}[scale = .6]
    \node[draw, circle, minimum size=5mm, inner sep=0pt]    (x) at (0,0)   {$x$};
    \node[draw, rectangle, minimum size=4mm, inner sep=0pt] (y) at (1.3,0) {$t$};
    \draw (y) -- (x);
\end{tikzpicture}. A seed torus is given by \(\bC^*_x \times \bC^*_t\). Put 
\[
x = \exp( u + i\theta), \qquad t = \exp( v + i\eta).  
\]
The holomorphic symplectic form $\Omega_{\cN}=\dd\log t\wedge \dd\log x=\dd y\wedge \dd x/(xy-1)$ extends over $x=0$, and
\(
 \Im\Omega_{\cN}= 
 \dd(v \dd\theta - u \dd \eta).
\)

Choose constants and smooth cutoff functions with
\[
 \begin{gathered}
 L>a>0,\qquad h\le0,\qquad 0\le h'\le1,\qquad u-h(u)<0\ \ (u<0),\\
 h(u)=\begin{cases}u+a,&u\le-L,\\ 0,&u\text{ near }[0,\infty);\end{cases}\\
 \rho(v)\in C_c^\infty(\bR,[0,1]),\qquad
 \rho(v)=1\text{ near }0;\\
 \beta\in C^\infty(\bR/2\pi\bZ),\qquad\supp\beta\subset(\pi-\varepsilon,\pi+\varepsilon),\qquad
 \beta(\eta)=\eta-\pi\text{ near }\pi, \\ (\eta-\pi)\beta(\eta)\ge0,\qquad\beta'\le1 .
 \end{gathered}
\]
On $(\pi-\varepsilon,\pi+\varepsilon)$, the symbol $\eta-\pi$ denotes the lift to $(-\varepsilon,\varepsilon)$.
Below, $h$, $\rho$ and $\beta$ are evaluated at $u$, $v$ and $\eta$, respectively.
Define
\begin{equation}\label{eq:node-primitive}
 \lN=v\,\dd\theta-u\,\dd\eta+\dd(h\rho\beta).
\end{equation}
The constant $a$ can be fixed before $L$ and $h$. Let \( L_{\cN, +} = \{ x, t > 0 \}\). Then \( \lambda_\cN|_{L_{\cN,+}} = 0\). 

\begin{lemma}\label{lem:node}
The form \(\lambda_\cN\) extends to a Liouville form of $\Im\Omega_{\cN}$, with skeleton
\[
 \Skel(Z_{\lN})=\{|x|=|t|=1\}\ \cup\ \{y=0,\ t=-1,\ |x|\le1\} .
\]
\end{lemma}

\begin{proof}
Near $x=0$ we have $u\le-L$ and $t$ near $-1$, so \(u \approx -\infty\), \(v \approx 0\), and \( \eta\approx\pi\). Hence
$\lN=v\,\dd\theta+(\eta-\pi)\,\dd u+a\,\dd\eta$.
Let $\Log$ be the principal branch of the logarithm on $\bC\setminus(-\infty,0]$.
Since $\Log(1-xy)=\Log(-t)=v+i(\eta-\pi)$ near $t=-1$,
\[
 \lN=\Im\Bigl(\frac{\Log(1-xy)}{x}\,\dd x\Bigr)+a\,\dd(\arg(1-xy)).
\]
The function $\Log(1-xy)/x$ extends holomorphically over $x=0$ with value $-y$.
Hence \(\lN\) can be extended to \(\cN\). 

Let \(Z = Z_{\lN}\). 
For a function $f$, write $\dot f=Z f$.
On the torus, the equation $\iota_{Z} \Im\Omega_{\cN}=\lN$ gives
\begin{equation}\label{eq:node-flow}
 \begin{aligned}
 \dot u&=u-h\rho\beta',&\qquad \dot v&=v,\\
 \dot\theta&=-h\rho'\beta,&\qquad \dot\eta&=h'\rho\beta.
 \end{aligned}
\end{equation}
So $|\dot u|+|\dot v|+|\dot\theta|+|\dot \eta|\le C(1+|u|+|v|)$.
On $\{x=0\}$ the field is $Z=(1+a)y\,\partial_y$, so $y(\tau)=e^{(1+a)\tau}y(0)$.
Thus $Z$ is complete.

Next we determine the skeleton.
By \eqref{eq:node-flow} and the choice of $h$ and $\beta$,
\[
 \dot u\le u-h<0\ \ (u<0),\qquad\dot u=u\ \ (u\ge0),\qquad (\eta-\pi)\dot \eta\ge0 .
\]
Let a forward orbit be relatively compact.
Since $v(\tau)=e^\tau v(0)$, we get $v=0$.
Suppose $u(0)<0$.
Then $u(\tau)\to-\infty$ while $|y|$ stays bounded, so $t=xy-1\to-1$ and $\eta\to\pi$.
Since $(\eta-\pi)\dot \eta\ge0$, this forces $\eta=\pi$, hence $t=-1$ and $y=0$.
If $x=0$, boundedness of $e^{(1+a)\tau}y(0)$ gives $y=0$.
If $u\ge0$, boundedness gives $u=v=0$.
So the skeleton is the stated set.

Next we construct a Liouville domain.
Fix $v_0>0$ and $R>1$, and consider
\[
 D=\{|x|\le R,\ |y|\le Y,\ e^{-v_0}\le|t|\le e^{v_0}\} .
\]
On the faces $|x|=R$ and $|t|=e^{\pm v_0}$,
\[
 Z\log|x|=\log R>0,\qquad Zv=v=\pm v_0,
\]
so $Z$ points outward there.
Choose $Y>e^{L}(1+e^{v_0})$.
Then $u\le-L$ on the face $|y|=Y$.
From $\log|y|=\log|1+t|-u$ and \eqref{eq:node-flow},
\[
 Z\log|y|=\Re\frac{t\,(v+i\rho\beta)}{1+t}+a+(-u-a)\bigl(1-\rho\beta'\bigr) .
\]
The first term extends continuously to $t=-1$ with value $1$.
So it is bounded on the closed annulus $e^{-v_0}\le|t|\le e^{v_0}$, by a bound independent of $a$.
Choose $a$ larger than the negative of its minimum.
The last term is nonnegative, since $-u-a\ge0$ and $\rho\beta'\le1$.
So $Z\log|y|>0$ on the face $|y|=Y$.
Hence $Z$ points outward on every face of $D$, and rounding the corners gives a Liouville domain.
Every backward orbit enters $D$.
On the torus this follows from \eqref{eq:node-flow}, and on $\{x=0\}$ from $y(-\tau)\to0$.
Hence $\cN=\widehat D$. \end{proof}

Finally, we explicitly describe  \(\cN\) as a handle attachment to  \( (\bC^*)^2\).
Choose $0<\delta\ll1$ with $\log\delta<-L$.
Near $\{|x|\le\delta,\ y=0\}$ use the holomorphic coordinates
\[ 
 \left(x, \xi=\frac{\Log(1-xy)}{x} \right), \qquad
\Omega_{\mathcal N}=\dd\xi\wedge \dd x.
\]
When $|x\xi|$ is small, $\rho=1$, $\beta=\eta-\pi$, and $x\xi=\Log(1-xy)=v+i(\eta-\pi)$.
Hence, in this region,
\begin{equation}\label{eq:cap-coordinates}
 \lN=\Im(\xi\,\dd x)+a\,\dd\Im(x\xi),\qquad Z(x,\xi)=\bigl(-ax,\ (1+a)\xi\bigr).
\end{equation}
Write $x=q_1+iq_2$ and $\xi=p_2+ip_1$.
Then \eqref{eq:cap-coordinates} implies
\[
 \lN=p\cdot \dd q+a\,\dd(p\cdot q),\qquad Zq=-aq,\qquad Zp=(1+a)p,
\]
and
\[
 Z|q|^2=-2a|q|^2,\qquad Z|p|^2=2(1+a)|p|^2,\qquad
 Z\bigl(-|q|^2+|p|^2\bigr)=2a|q|^2+2(1+a)|p|^2 .
\]
So the box $\{|x|\le\delta,\ |\xi|\le P_0\}$ with $\delta P_0$ small is a Weinstein handle of index two \cite[Chapter~11]{CE}.
The Liouville field points inward on the face $\{|x|=\delta\}$ and outward on the face $\{|\xi|=P_0\}$, and the core of the handle is the disk $\{\xi=0\}$.
Hence if we take the domain of \( (\bC^*)^2_{x, t} \) to be \( D_{(\bC^*)^2} = \{ \delta \le |x| \ll \infty , \  | \log |t| | \ll \infty   \}\), and \(\cF_\cN \coloneqq \{ |x| = \delta, \ \arg (t) = \pi \} \), then 
\begin{equation}
\label{eq:handle_N}
\cN \cong \widehat {D_{(\bC^*)^2} \cup_{\cF_\cN} T^*D^2}. 
\end{equation}

\subsection{Gluing Liouville forms}
We give a general construction. 
The input is a manifold $X$ glued from Liouville manifolds $V$ and $\cN \times \Sigma$ along $\{x\ne0\}$, where \(\Sigma\) is another Liouville manifold. 
The output is a Liouville form on $X$.

\begin{setup}\label{setup:attachment}
The data consist of the following.
\begin{enumerate}[label=\textup{(A\arabic*)}]
\item\label{A1} Liouville manifolds $(V,\lambda_V)$ and $(\Sigma,\lambda_\Sigma)$ with Morse--Bott Lyapunov functions $\Psi_V$ and $\Psi_\Sigma$, and a constant $\kappa>0$.

\item\label{A2} A manifold $X=V\cup U$ with $U=\cN\times\Sigma$ and $V\cap U=\cN^\circ\times\Sigma$, where $\cN^\circ=\cN\cap\{x\ne0\}\cong(\bC^*)^2_{(x,t)}$.

\item\label{A3} A symplectic form $\omega$ on $X$ with $\omega|_V=\dd\lambda_V$ and $\omega|_U=\dd\lambda_U$, where $\lambda_U=\kappa\lN+\lambda_\Sigma$.

\item\label{A4} A smooth map $x:V\to\bC^*$ that agrees with the coordinate $x$ of $\cN$ on $V\cap U$, with $Z_V\log|x|=\log|x|$.
A circle action on $V$ with generator $\partial_\eta$ satisfies $\iota_{\partial_\eta}\omega=\kappa\,\dd\log|x|$, and on $V\cap U$ it is $e^{i\alpha}\cdot(x,t,\sigma)=(x,e^{i\alpha}t,\sigma)$.

\item\label{A5} A compact torus $G$ acting on $X$ by symplectomorphisms that preserve $V$, $U$, $\lambda_V$ and $\lambda_U$.
On $U = \cN \times \Sigma $ it acts by $g\cdot(x,t,\sigma)=(e^{i\alpha(g)}x,\,t,\,g\cdot\sigma)$, for a character $\alpha$ of $G$ and an action on $\Sigma$ preserving $\lambda_\Sigma$ and $\Psi_\Sigma$.
For every $\xi$ in the Lie algebra of $G$, $\lambda_V(\xi)=\lambda_U(\xi)$ on $V\cap U$.

\item\label{A6} A Lagrangian $L_+\subset X$ such that $L_+\cap V$ is properly embedded in $V$, $\lambda_V|_{L_+\cap V}=0$, and $L_+\cap U=L_{\cN, +} \times L_{\Sigma, +}$.
Here $L_{\Sigma, +}\subset\Sigma$ is a proper strict exact Lagrangian, and $L_+\cap V\cap U$ is connected.
\item\label{A7} The closed form $\lambda_U-\lambda_V$ on $V\cap U$ is exact.
\end{enumerate}

\end{setup}

By \ref{A7}, $\lambda_U-\lambda_V=\dd F$ for a function $F$ on $V\cap U$.
Both $\lambda_V$ and $\lambda_U$ vanish on $L_+\cap V\cap U$, which is connected by \ref{A6}.
So $F$ is constant on $L_+\cap V\cap U$, and we normalize
\begin{equation}\label{eq:normalization_F}
 \lambda_U-\lambda_V=\dd F,\qquad F|_{L_+\cap V\cap U}=0 .
\end{equation}
The function $F$ is $G$-invariant.
Indeed, for $\xi$ in the Lie algebra of $G$, \ref{A5} gives $\dd F(\xi)=\lambda_U(\xi)-\lambda_V(\xi)=0$.

By \ref{A4}, $\Zero(Z_V)\subset\{\log|x|=0\}$, and
\[
 \lambda_V(\partial_\eta)=\omega(Z_V,\partial_\eta)=-\kappa\,Z_V\log|x|=-\kappa\log|x| .
\]

Since $G$ is compact, there is a proper $G$-invariant function $R_\Sigma\ge1$ on $\Sigma$.
Let $L$ be the constant of \eqref{eq:node-primitive}.
Choose $A_0<A_1<-L$, choose $0<\epsilon_0<\epsilon_1$ and $C_0>0$, and choose cutoffs with values in $[0,1]$ such that
\[
\begin{array}{c|cc}
 &\text{equals }1\text{ on}&\text{equals }0\text{ on}\\ \hline
 \chi_x(\log|x|)&(-\infty,A_0]&[A_1,\infty)\\
 \chi_t(\log|t|)&[-\epsilon_0,\epsilon_0]&\bR\setminus(-\epsilon_1,\epsilon_1)\\
 \chi_\Sigma(R_\Sigma|x|/C_0)&[0,1]&[2,\infty)
\end{array}
\]
Put $\chi=\chi_x\chi_t\chi_\Sigma$ on $V\cap U$, and define
\[
 \lambda_X= \lambda_V+\dd(\chi F)
\]
Near \( \{x = 0\} \) we have $\chi=1$ and $\lambda_V+\dd F=\lambda_U$. Hence \(\lambda_X\) is defined globally on \(X\) and \(\dd \lambda_X = \omega\).  

\begin{proposition}
\label{prop:gluing_Liouville_form}
Under Setup~\ref{setup:attachment}, the one-form $\lambda_X$ on $X$  has the following properties.
\begin{enumerate}[label=\textup{(\roman*)}]
\item 
\label{item:i}
$(X,\lambda_X)$ is a Liouville manifold.
For small $\delta>0$,
\[
 \Skel(Z_X)=\Skel(Z_V)\ \cup\ \overline{\textstyle\bigcup_{\tau\ge0}\phi_{Z_X}^{-\tau}(K_{\capp})},
 \qquad K_{\capp}=\{|x|\le\delta,\ y=0,\ t=-1\}\times\Skel(Z_\Sigma).
\]
\item\label{item:ii} $\lambda_X=\lambda_V$ on a neighborhood of a Liouville domain $D_V$ of $V$, and $(V,\lambda_X|_V)$ is the completion of $(D_V,\lambda_V)$.

\item\label{item:iii} $L_+$ is a proper strict exact Lagrangian for $\lambda_X$.

\item\label{item:iv} $\lambda_X$ is $G$-invariant, and $\lambda_X(\xi)=\lambda_V(\xi)$ on $V$ for $\xi$ in the Lie algebra of $G$.

\item\label{item:v} $X$ has a $G$-invariant Morse--Bott Lyapunov function $\Psi_X$ for $\lambda_X$.
Its critical components are those of $\Psi_V$, with the same germs, and the components $\{x=y=0,\ t=-1\}\times C_\Sigma$ for $C_\Sigma\subset\Crit\Psi_\Sigma$.
The new components have normal index $2+\ind_{C_\Sigma}\Psi_\Sigma$.
\item\label{item:vii} Let \(\cF \coloneqq \cF_\cN \times \Sigma\) be a stop for \(O \coloneqq V \cap U\). Then we have inclusion of sectors
\begin{equation}\label{eq:geometric-square}
 \begin{tikzcd}[column sep=large]
 (D_O,\cF)\arrow[r,hook]\arrow[d] & (D_V,\cF)\arrow[d]\\
 D_O\arrow[r,hook] & D_V,
 \end{tikzcd}
 \qquad D_O\subset D_V\subset D_X.
\end{equation}
Here $O=V\cap U$ and
\[
 \widehat D_O\cong(O,\lambda_O),\quad
 \widehat D_V\cong(V,\lambda_V),\quad
 \widehat D_X\cong(X,\lambda_X),\quad
 \lambda_O=\kappa \lambda_\cN + \lambda_\Sigma. 
\]
The group action \ref{A4} defines a positive wrapping of \(\cF\) in both \(\partial_\infty O\) and \( \partial_\infty V\). 

Moreover,
there is a sectorial gluing
\begin{equation}\label{eq:abstract-sectorial-gluing}
 (X,\lambda_X)\cong
 \widehat{\,(V,\cF)\underset{\cF}{\cup}(T^*D^2\times\Sigma)\,}.
\end{equation}
\end{enumerate}
\end{proposition}

We use the following criterion for finite type.

\begin{lemma}\label{lem:finite-type}
Let $\lambda$ be a one-form on $X$ with $\dd\lambda$ symplectic and complete Liouville field $Z$.
Suppose that $\Skel(Z)$ is compact and that
\[
 \bigcup_{\tau\ge0}\phi_Z^{-\tau}(K)\Subset X\qquad\text{for every compact }K\subset X .
\]
Then $(X,\lambda)$ is a Liouville manifold (of finite type).
\end{lemma}

\begin{proof}
We first show that the backward flow is uniformly attracted to the skeleton.
Let $K$ be compact and $U\supset\Skel(Z)$ open.
We claim that there is $T>0$ with $\phi_Z^{-\tau}(K)\subset U$ for all $\tau\ge T$.
Otherwise there are $p_j\in K$ and $\tau_j\to\infty$ with $\phi_Z^{-\tau_j}(p_j)\notin U$.
By hypothesis, a subsequence of $\phi_Z^{-\tau_j}(p_j)$ converges to some $q\notin U$.
For every $\sigma\ge0$,
\[
 \phi_Z^{\sigma}(q)=\lim_j\phi_Z^{\sigma-\tau_j}(p_j)\in\overline{\textstyle\bigcup_{\tau\ge0}\phi_Z^{-\tau}(K)} .
\]
The right side is compact, so $q\in\Skel(Z)\subset U$, a contradiction.

Next we construct a proper function $R$ with $ZR>0$ outside a compact set.
Choose a proper function $R_0\ge0$, a compact set $K$, and a smooth function $\psi\ge0$ with
\[
 \psi=0\ \text{near }\Skel(Z),\qquad \psi>1+|ZR_0|\ \text{on }X\setminus K.
\]
Put
\[
 R(p)=\int_{-\infty}^0\psi\bigl(\phi_Z^{\tau}(p)\bigr)\,\dd\tau .
\]
By uniform attraction, the integrand vanishes for $\tau\le-T(p)$, with $T(p)$ locally bounded in $p$.
So $R$ is smooth and $ZR=\psi$.
Let $p\notin K$, and let $\tau_0\le0$ be the last time with $\phi_Z^{\tau_0}(p)\in K$, or $\tau_0=-\infty$ if there is none.
Then
\[
 R_0(p)-\max_KR_0\le\int_{\tau_0}^0|ZR_0|\bigl(\phi_Z^\tau(p)\bigr)\,\dd\tau\le\int_{\tau_0}^0\psi\bigl(\phi_Z^\tau(p)\bigr)\,\dd\tau\le R(p).
\]
So $R$ is proper.
For large $c$, the set $D=\{R\le c\}$ is compact and $ZR=\psi>0$ on $\partial D$.
By uniform attraction, $X=\bigcup_{\tau\ge0}\phi_Z^\tau(D)$, so $D$ is a Liouville domain and $X=\widehat D$.
\end{proof}

\begin{proof}[Proof of Proposition~\ref{prop:gluing_Liouville_form}]
Let $u=\log|x|$, and let $\Phi_X^\tau$ be the flow of $Z_X$.

\ref{item:i} We first prove that $Z_X$ is complete.
Let $X_\chi$ be defined by $\iota_{X_\chi}\omega=\dd\chi$.
On $V\cap U$,
\[
 Z_X=Z_V+\chi(Z_U-Z_V)+F\,X_\chi .
\]
The functions $|x|$, $|t|$ and $R_\Sigma$ Poisson-commute with $u$.
For $|x|$ and $|t|$ this holds because $\omega=\kappa(\dd v\wedge \dd\theta+\dd\eta\wedge \dd u)+\omega_\Sigma$ pairs moduli with arguments.
For $R_\Sigma$ it holds because $R_\Sigma$ is a function on the factor $\Sigma$.
So $X_\chi u=0$, and \eqref{eq:node-flow} together with \ref{A4} and $Z_{\kappa\lN}=Z_{\lN}$ gives
\begin{equation}\label{eq:attachment-monotone}
 Z_Xu=u-\chi\,h\rho\beta' .
\end{equation}
Since $h\le0$ and $\chi\rho\beta'\le1$,
\[
 Z_Xu\le u-h<0\ \ (u<0),\qquad Z_Xu=u\ \ (u\ge0),\qquad |Z_Xu|\le C_1(1+|u|) \text{ for some constant } C_1.
\]
Moreover there is $c_1>0$ with
\begin{equation}\label{eq:attachment-speed}
 -Z_Xu\ge c_1\qquad(u\le A_1),
\end{equation}
because $h(u)-u>0$ on $(-\infty,A_1]$ and $h(u)-u=a$ for $u\le-L$.
Fix $u_-<u_+$.
We have $\chi_t(\log|t|)=0$ for $|\log|t||\ge\epsilon_1$, and $\chi_\Sigma(R_\Sigma|x|/C_0)=0$ for $R_\Sigma|x|\ge2C_0$.
On $\{u\ge u_-\}$ the second condition gives $R_\Sigma\le2C_0e^{-u_-}$ on $\supp\chi$.
Hence
\[
 \supp \dd(\chi F)\cap\{u_-\le u\le u_+\}\subset
 K_{u_-,u_+}=\bigl\{e^{u_-}\le|x|\le e^{u_+},\ e^{-\epsilon_1}\le|t|\le e^{\epsilon_1},\ R_\Sigma\le2C_0e^{-u_-}\bigr\}.
\]
The set $K_{u_-,u_+}$ is compact, since $R_\Sigma$ is proper on $\Sigma$.
Since $\iota_{Z_X-Z_V}\omega=\dd(\chi F)$, we have $Z_X=Z_V$ on $\{u_-\le u\le u_+\}\setminus K_{u_-,u_+}$.
By $|Z_Xu|\le C_1(1+|u|)$, every trajectory satisfies $1+|u(\tau)|\le(1+|u(0)|)e^{C_1|\tau|}$.
So for $\tau\in[-T,T]$, a trajectory starting in $V$ stays in $\{u_-\le u\le u_+\}$ with $u_\pm=\pm\bigl((1+|u(0)|)e^{C_1T}-1\bigr)$.
There $Z_X=Z_V$ outside the compact set $K_{u_-,u_+}$, and $Z_V$ is complete.
On $\{x=0\}$ the field is $Z_X=(-ax,(1+a)\xi)+Z_\Sigma$ in the coordinates \eqref{eq:cap-coordinates}, which is complete.
So $Z_X$ is complete.

Next we show that for every compact $K\subset X$ (recall \( \Subset\) means relatively compact),
\begin{equation}\label{eq:attachment-saturation}
 \bigcup_{\tau\ge0}\Phi_X^{-\tau}(K)\Subset X.
\end{equation}

Suppose first that $K\subset V$.
While a backward trajectory lies in $\{u\le A_1\}$, its coordinate $u$ increases at speed at least $c_1$ by \eqref{eq:attachment-speed}.
So it reaches $\{u=A_1\}$ in bounded time.
On $\{u\ge A_1\}$ we have $\chi=0$, so $Z_X=Z_V$ there, and $Z_Vu=u$ keeps the backward trajectory in $\{u\ge A_1\}$.
The rest of the trajectory is a backward $Z_V$-trajectory starting in a compact set.
Since $V$ is of finite type, the union of these trajectories is relatively compact.
The finite segments depend continuously on the initial point, so their union is relatively compact too.

Now let $K$ meet $\{x=0\}$.
It remains to treat a neighborhood of $K\cap\{x=0\}$.
In the coordinates \eqref{eq:cap-coordinates}, cover it by finitely many product boxes $\{|x|\le\delta,\ |\xi|\le b\}\times K_\Sigma$ with $K_\Sigma\subset\Sigma$ compact.
Let $K_\Sigma^-$ be the closure of $\bigcup_{\tau\ge0}\phi_{Z_\Sigma}^{-\tau}(K_\Sigma)$.
It is compact, because $\Sigma$ is of finite type by \ref{A1}.
Choose $\delta$ so small that
\[
 \log\delta<A_0,\qquad \delta\sup_{K_\Sigma^-}R_\Sigma<C_0,
\]
and that $|x\xi|\le\delta b$ lies in the plateaux of $\rho$, $\beta$ and $\chi_t$.
Then $\lambda_X=\lambda_U$ on the box $\{|x|\le\delta,|\xi|\le b\}\times K_\Sigma^-$.
By \eqref{eq:cap-coordinates}, a backward trajectory starting in the smaller box satisfies
\[
 x(-\tau)=e^{a\tau}x(0),\qquad\xi(-\tau)=e^{-(1+a)\tau}\xi(0),\qquad\sigma(-\tau)=\phi_{Z_\Sigma}^{-\tau}\sigma(0)
\]
until it leaves through $\{|x|=\delta\}$.
Points with $x=0$ never leave.
The exit points with $x\ne0$ form a compact subset of $V$, and the first case applies to them.
This proves \eqref{eq:attachment-saturation}.

Next we determine the skeleton.
Let $z$ have a relatively compact forward trajectory.
If $z\in V$ and $u(z)>0$, then $u$ grows exponentially, which is impossible.
If $u(z)=0$, then the trajectory is a $Z_V$-trajectory in $\{u=0\}$, and it is relatively compact if and only if $z\in\Skel(Z_V)$.
If $u(z)<0$, then $u\to-\infty$ by \eqref{eq:attachment-monotone}.
Relative compactness and $xy=1+t$ then give $t\to-1$ and a bounded coordinate in $\Sigma$.
So the tail of the trajectory lies in a product box.
There the computation in the proof of Lemma~\ref{lem:node} gives $y=0$ and $t=-1$, and the coordinate in $\Sigma$ lies in $\Skel(Z_\Sigma)$.
The same holds on $\{x=0\}$, where $\xi$ expands unless $\xi=0$.
The set $K_{\capp}$ is compact and forward invariant.
By \eqref{eq:attachment-saturation}, $\bigcup_{\tau\ge0}\Phi_X^{-\tau}(K_{\capp})$ has compact closure.
This closure is forward invariant, so its points have relatively compact forward trajectories.
Hence
\[
 \Skel(Z_X)=\Skel(Z_V)\cup
 \overline{\textstyle\bigcup_{\tau\ge0}\Phi_X^{-\tau}(K_{\capp})}.
\]
This set is compact.
By Lemma~\ref{lem:finite-type}, this and \eqref{eq:attachment-saturation} prove that \(\lambda_X\) is a Liouville structure of finite type.

\ref{item:ii} Since $\log|x|$ is bounded below on compact subsets of $V$, we may choose $A_1$ below $\min_{D_V}\log|x|$ for a Liouville domain $D_V$ of $V$.
Then $\chi=0$ near $D_V$, so $\lambda_X=\lambda_V$ there.
By the first case in the proof of \ref{item:i}, every backward $Z_X$-trajectory in $V$ enters $\{u\ge A_1\}$, and then $D_V$.
So $(V,\lambda_X|_V)$ and $(V,\lambda_V)$ are both the completion of $(D_V,\lambda_V)$.

\ref{item:iii} Since $x>0$ on $L_{\cN, +}$, we have $L_+\subset V$.
By \eqref{eq:normalization_F}, $(\chi F)|_{L_+}=0$, so $\lambda_X|_{L_+}=\lambda_V|_{L_+}=0$ by \ref{A6}.
The set $L_+$ is closed in $V$ by \ref{A6}.
A limit point of $L_+$ on $\{x=0\}$ would have $t=-1$, whereas $t>0$ on $L_+\cap U$.
So $L_+$ is properly embedded in $X$.

\ref{item:iv} The cutoffs are functions of $|x|$, $|t|$ and $R_\Sigma$, so they are $G$-invariant by \ref{A5}.
Since $\lambda_V$ and $F$ are $G$-invariant, so are $\chi F$ and $\lambda_X$, and $\dd(\chi F)(\xi)=\xi(\chi F)=0$.

\ref{item:v} We construct $\Psi_X$.
All zeros of $Z_V$ lie in $\{u=0\}$, where $Z_X=Z_V$.
By \eqref{eq:attachment-monotone}, $Z_X$ has no zeros in $V\cap\{u\ne0\}$.
On $\{x=0\}$ the field is the product field, so
\[
 \Zero(Z_X)=\Zero(Z_V)\ \sqcup\ \bigl(\{x=y=0,\ t=-1\}\times\Zero(Z_\Sigma)\bigr).
\]
In the real coordinates $(q,p)$ of Section~\ref{sec:local-models}, put $\Psi_{\capp}=-|q|^2+|p|^2+\Psi_\Sigma$.
Since $\kappa\lN$ and $\lN$ have the same Liouville field,
\[
 Z_X\Psi_{\capp}=2a|q|^2+2(1+a)|p|^2+Z_\Sigma\Psi_\Sigma
\]
on the product box.
Let $D$ be a Liouville domain of $X$ containing $\Zero(Z_X)$.
Inside $D$, a neighborhood of $\{x=0\}$ lies in the product box, since $y$ is bounded on $D$.
Choose a cutoff $\psi_0$ on $D$ with $\psi_0=1$ near $\{x=0\}$ and $\supp\psi_0$ in the product box, and put $\Psi_0=\psi_0\Psi_{\capp}+(1-\psi_0)\Psi_V$.
Then $\{Z_X\Psi_0\le0\}\setminus\Zero(Z_X)$ lies in a compact band $K_0\subset D\cap\{u_-\le u\le u_+\}$ with $u_+<0$.
Choose a smooth nonincreasing function $b_0$ with
\[
 b_0(u)=b_0(u_--1)\ \ (u\le u_--1),\qquad b_0=0\ \text{near }[u_+/2,\infty),\qquad b_0'\le-1\ \text{on }[u_-,u_+].
\]
It extends over $\{x=0\}$ by its constant value.
Put
\[
 c_0=\min_{K_0}(-Z_Xu)>0,\qquad M>c_0^{-1}\max_{K_0}(-Z_X\Psi_0)+1,\qquad \Psi_D=\Psi_0+Mb_0(u).
\]
Then $Z_X\Psi_D>0$ on $D\setminus\Zero(Z_X)$.
Indeed, $Z_X(b_0(u))=b_0'(u)Z_Xu\ge0$ everywhere, and $Z_X(Mb_0(u))\ge Mc_0$ on $K_0$.
Near $\Zero(Z_X)$ the function $b_0(u)$ is constant.
Extend $\Psi_D$ across the collar $\bigcup_{\tau\ge0}\phi_{Z_X}^\tau(\partial D)$ by a function increasing in $\tau$.
The result $\Psi_X$ is proper and bounded below, with $\Crit\Psi_X=\Zero(Z_X)$ and $\dd\Psi_X(Z_X)>0$ elsewhere.
Near the old components $\Psi_X=\Psi_V$, so the germs are those of $\Psi_V$.
Near a new component $C=\{x=y=0,\ t=-1\}\times C_\Sigma$ we have $\Psi_X=\Psi_{\capp}+Mb_0(u_--1)$, so
\[
 \begin{gathered}
 \ind_C\Psi_X=2+\ind_{C_\Sigma}\Psi_\Sigma,\qquad\dim C=\dim C_\Sigma,\\
 \ind_C\Psi_X+\dim C\le2+\tfrac12\dim_\bR\Sigma=\tfrac12\dim_\bR X .
 \end{gathered}
\]
The gradient-like estimate of Definition~\ref{def:quantitative}(iii) holds near the old components by \ref{A1}.
Near a new component it follows from the estimate for $\Psi_\Sigma$ and the formula for $Z_X\Psi_{\capp}$.
All choices can be made $G$-invariant.
Indeed, $\Psi_V$ and $\Psi_\Sigma$ are $G$-invariant, and $g\cdot(x,\xi)=(e^{i\alpha(g)}x,e^{-i\alpha(g)}\xi)$, so $\Psi_{\capp}$ is $G$-invariant.
The cutoff $\psi_0$ and the collar extension can be averaged over $G$.

\ref{item:vii} Fix a Liouville domain $W_\Sigma$ of $\Sigma$ and a constant $M_O>L$, and let $D_O\subset O$ be the rounding of
\[
 \{|u|\le M_O,\ |v|\le1\}\times W_\Sigma .
\]
For $0\le s\le1$ put $\lambda_O^s=\lambda_O-s\kappa\,\dd(h\rho\beta)$.
Then $\lambda_O^0=\lambda_O$, and $\lambda_O^1=\kappa(v\,\dd\theta-u\,\dd\eta)+\lambda_\Sigma$ is the product primitive.
By \eqref{eq:node-flow}, the Liouville field $Z^s$ of $\lambda_O^s$ has $\Sigma$-component $Z_\Sigma$ and satisfies
\[
 Z^su=u-(1-s)h\rho\beta',\qquad Z^sv=v .
\]
Since $h\le0$ and $(1-s)\rho\beta'\le1$, we get $Z^su=u$ for $u\ge0$ and $Z^su\le u-h<0$ for $u<0$, with $u-h=-a$ for $u\le-L$.
So $Z^s$ points outward along $\partial D_O$, and $(D_O,\lambda_O^s)$ is a Liouville domain.
As in the proof of Lemma~\ref{lem:node}, $Z^s$ is complete.
Every backward orbit of $Z^s$ in $O$ enters $D_O$: $v\to0$, the $\Sigma$-component enters $W_\Sigma$, and $u$ enters $[-M_O,M_O]$.
Hence $\widehat D_O\cong(O,\lambda_O^s)$ for every $s$.

Next we nest $D_O$ in $D_V$ and $D_X$.
Let $B\subset O$ be a compact neighborhood of $D_O$, and let $\chi_O$ be a cutoff supported in $B$ and equal to one near $D_O$.
Every compact subset of $V$ lies in $\phi_{Z_V}^T(D)$ for a Liouville domain $D$ of $V$ and large $T$, and $\phi_{Z_V}^T(D)$ is again a Liouville domain.
So we may take the domain $D_V$ in the proof of \ref{item:ii} with $B\subset\operatorname{int}D_V$, and then choose $A_1$.
By \eqref{eq:normalization_F}, $\lambda_O^s-\lambda_V=\dd(F-s\kappa h\rho\beta)$ on $O$.
This potential vanishes on $L_+\cap O$, since $\eta=0$ there and $\beta(0)=0$.
Fix $s$ and put $K=\chi_O\,(F-s\kappa h\rho\beta)$, extended by zero.
The compact exact deformations
\[
 \lambda_V+r\,\dd K,\qquad \lambda_X+r\,\dd K,\qquad 0\le r\le1,
\]
are supported in $\operatorname{int}D_V$, agree near $D_V$, vanish on $L_+$, and end at $\lambda_O^s$ near $D_O$.
By \ref{item:i}, $D_V$ lies in the interior of a Liouville domain $D_X$ of $X$.
The deformations do not change $\lambda_V$ and $\lambda_X$ near $\partial D_V$ and $\partial D_X$, so $D_V$ and $D_X$ remain Liouville domains.
This gives Liouville inclusions $D_O\subset\operatorname{int}D_V\subset\operatorname{int}D_X$ for the endpoint primitives.
Let $Y$ be the field with $\iota_Y\omega=-\dd K$, and $\varphi_r$ its flow.
Then
\[
 \varphi_1^*(\lambda_X+\dd K)-\lambda_X=\dd\int_0^1\varphi_r^*\bigl((\lambda_X+r\,\dd K)(Y)\bigr)\,\dd r,
\]
and the integrand is compactly supported.
The field $Y$ is tangent to $L_+$, since $K$ vanishes on $L_+$.
So $\widehat D_X\cong(X,\lambda_X)$ by an exact symplectomorphism preserving $L_+$, and in the same way $\widehat D_V\cong(V,\lambda_V)$.

Take $s=1$.
Since $Z_Vu=u$ by \ref{A4}, the flow tube of the stop in $V$ has $u(\tau)=-M_Oe^\tau$ and does not return.
Pushing $\partial D_O$ and $\partial D_V$ along this tube until they agree near the stop, and removing the same neighborhood of the stop from both, gives the stopped pairs with identical side collars and the square \eqref{eq:geometric-square}.
The circle of \ref{A4} moves $\cF$ positively, since $\lambda_V(\partial_\eta)=-\kappa\log|x|>0$ near $\cF$.
The gluing \eqref{eq:abstract-sectorial-gluing} is \eqref{eq:handle_N} times $\Sigma$, attached to $D_V$ by the backward flow tube of the face $\{|x|=\delta\}$, along which $u$ is monotone by \eqref{eq:attachment-speed}.
\end{proof}

\subsection{Liouville structures on unimodular cluster varieties}\label{sec:cluster-liouville}

We first record how the completion changes under mutation.

\begin{lemma}\label{lem:completion-mutation}
Let $(E,\widetilde D)$ be a completion of $\widetilde B$ and $k\in M$.
Define $L,F\in\Mat_{m+f}(\bZ)$ by
\[
 L=J+\sum_{i\ne k}[-E_{ik}]_+\,e_ie_k^{\T},\qquad
 F=J+\sum_{j\ne k}[E_{kj}]_+\,e_ke_j^{\T},\qquad
 J=I_{m+f}-2e_ke_k^{\T}.
\]
Then
\[
 \mu_kE=LEF,\qquad L^2=I_{m+f},\qquad F\widetilde D=\widetilde DL^{\T},\qquad (\mu_kE)\widetilde D=L(E\widetilde D)L^{\T}.
\]
In particular $(\mu_kE,\widetilde D)$ is a completion of $\mu_k\widetilde B$.
Moreover $\mu_k(E^\vee)=(\mu_kE)^\vee$.
\end{lemma}

\begin{proof}
The identity $\mu_kE=LEF$ is \cite[(3.2) in the proof of Lemma~3.2]{BFZ-III} with $\varepsilon=1$.
It uses only $E_{kk}=0$.
The matrix $L$ differs from $I_{m+f}$ only in column $k$, and $Le_k=-e_k+\sum_{i\ne k}[-E_{ik}]_+e_i$.
So $L^2e_k=e_k$, and $L^2=I_{m+f}$.
The identity $E_{kj}d_j=-d_kE_{jk}$ gives $[E_{kj}]_+d_j=d_k[-E_{jk}]_+$, which is the $(k,j)$ entry of $F\widetilde D=\widetilde DL^{\T}$.
The other entries agree trivially.
Hence $(\mu_kE)\widetilde D=LE\widetilde DL^{\T}$, which is skew-symmetric.
Since $\det L=\det F=-1$, the matrix $\mu_kE$ lies in $GL_{m+f}(\bZ)$.
Its columns in $M$ are $\mu_k\widetilde B$, because the mutation rule for a column $j\in M$ only involves the entries of $\widetilde B$.
The matrices $L$ and $F$ built from $E^\vee$ are $F^{\T}$ and $L^{\T}$.
So $\mu_k(E^\vee)=F^{\T}E^\vee L^{\T}=(LEF)^\vee$.
\end{proof}

Now we prove Proposition~\ref{prop:Omega}. 
\begin{proof}[Proof of Proposition~\ref{prop:Omega}]
    We first show that \(\Omega\) is invariant under mutation, so it is defined on the union of seed tori \(\cA^{\mathrm{tori}}(\widetilde B)\). Since \(B\) is locally acyclic, by~\cite[Proposition 7.6]{Lam-Speyer}, \(\codim_\bC \left( \cA(\widetilde B) \setminus \cA^{\mathrm{tori}}(\widetilde B) \right) \ge 2\). Hence \(\Omega\) is globally defined.  
    
    Let $L$ be as in Lemma~\ref{lem:completion-mutation}, and put $K=E\widetilde D$.
    Put $\eta=\sum_iE_{ik}\,\dd\log z_i$ and $f=z^{c_k}/(1+z^{c_k})$, where \(c_k\) is given in ~\eqref{eq:dual-circles}. 
    Let us consider a mutation at \(k\). 
    Taking $\dd\log$ of \eqref{eq:exchange}  gives 
    $\dd\log z'=L^{\T}\dd\log z+f\,e_k\,\eta$.
    Using $L^2=I_{m+f}$ and $\eta\wedge\eta=0$, we get
    \[
    \Omega'=\tfrac12(\dd\log z')^{\T}LKL^{\T}(\dd\log z')
    =\Omega+f\big((\dd\log z)^{\T}K(L^{\T}e_k)\wedge \eta+\eta\wedge(L^{\T}e_k)^{\T}K\, \dd \log z\big)
    \]
    Now $L^{\T}e_k=-e_k$, and $K e_k=d_kc_k$ by \eqref{eq:dual-circle-identities}.
    So $(\dd\log z)^{\T}Ke_k=d_k\eta$, and both terms in the bracket are multiples of $\eta\wedge\eta=0$.
    Hence $\Omega'=\Omega$.
    \end{proof}
    
\begin{lemma}\label{lem:form-invariance}
Let \(k\) be a sink of \(B\). On $\cA(\widetilde B)_t\cong\cN\times\Sigma$,
\begin{equation}\label{eq:form-splitting}
 \Omega=d_k\,\dd\log t\wedge \dd\log x+\Omega_\Sigma,\qquad
 \ltor=d_k\bigl(\log|t|\,\dd(\arg x)-\log|x|\,\dd(\arg t)\bigr)+\ltor^\Sigma .
\end{equation}
Here $\Omega_\Sigma$ and $\ltor^\Sigma$ are the forms \eqref{eq:forms} of the completion $(E_\Sigma,D_\Sigma)$ of Proposition~\ref{prop:Sigma-completion}.
\end{lemma}

\begin{proof}
Let the $(x,t,u_j,y_\alpha)$ be as in \eqref{eq:Sigma-variables}, which are monomials in the seed coordinates, given by the basis $\{\epsilon_k,c_k\}\cup\{\epsilon_j|_{N_\Sigma}\}\cup\{a_\alpha|_{N_\Sigma}\}$ of $N_{\bZ}^*$.
Both expressions in \eqref{eq:forms} are covariant under such a change of basis.
So \eqref{eq:form-splitting} follows from $K=K_k\oplus K_\Sigma$ and $K_k(w_k,e_k)=d_k$ in Lemma~\ref{lem:lattice-splitting}(ii).
\end{proof}

Next we compute the pairing of a primitive with the frozen circles.

For $a\in M\sqcup F$, let $\xi_a$ be the vector field given by the circle $e^{i\alpha}\cdot z_j=e^{i\alpha(w_a)_j}z_j$ of the seed torus.
For $a\in F$ it acts on $\cA(\widetilde B)$.
On the seed torus, $\xi_a=\sum_j(w_a)_j\partial_{\theta_j}$, and \eqref{eq:forms} and \eqref{eq:dual-circle-identities} give
\begin{equation}\label{eq:frozen-hamiltonian}
 \iota_{\xi_a}\Im\Omega=-\sum_{i,j}K_{ij}(w_a)_j\,\dd u_i=d_a\,\dd u_a=d_a\,\dd\log|z_a| .
\end{equation}
Let $\lambda$ be a $\xi_a$-invariant primitive of $\Im\Omega$ with $Z_\lambda\log|z_a|=\log|z_a|$.
Then
\begin{equation}\label{eq:frozen-pairing}
 \lambda(\xi_a)=\Im\Omega(Z_\lambda,\xi_a)=-d_a\,Z_\lambda\log|z_a|=-d_a\log|z_a| .
\end{equation}

Next we compare the torus primitives of adjacent seeds.
Let $\mathcal D$ be the Bloch--Wigner dilogarithm.
It is a single-valued real-analytic function on $\bC\setminus\{0,1\}$ with
\[
 \dd\mathcal D(z)=\log|z|\,\dd(\arg(1-z))-\log|1-z|\,\dd(\arg z) .
\]

\begin{lemma}
\label{lem:mutation-potential}
Let $j\in M$, and let $\ltor$ and $\ltor'$ be the torus primitives \eqref{eq:forms} of the seed and of its mutation at $j$.
On their common torus,
\[
 \ltor'-\ltor=-d_j\,\dd\mathcal D(-z^{c_j}).
\]
A closed $1$-form on either seed torus is exact as soon as its restriction to the common torus is exact.
\end{lemma}

\begin{proof}
We first treat $\cN = \{xb=1+t\}$, where $\Omega=d_k\,\dd\log t\wedge \dd\log x=d_k\,\dd\log b\wedge \dd\log t$.
The torus primitives of the charts $(x,t)$ and $(b,t)$ are
\[
 \ltor^{(x,t)}=d_k(\log|t|\,\dd(\arg x)-\log|x|\,\dd(\arg t)),\qquad
 \ltor^{(b,t)}=d_k(\log|b|\, \dd (\arg t)-\log|t|\,\dd(\arg b)).
\]
Since $\log|b|=\log|1+t|-\log|x|$ and $\arg b=\arg(1+t)-\arg x$,
\begin{equation}\label{eq:mutation-potential}
 \ltor^{(b,t)}-\ltor^{(x,t)}
 =d_k\bigl(\log|1+t|\,\dd(\arg t)-\log|t|\,\dd(\arg(1+t))\bigr)
 =-d_k\,\dd\mathcal D(-t).
\end{equation}

Now let $j\in M$ be arbitrary.
Lemma~\ref{lem:lattice-splitting}(i) and (ii) at $j$ give coordinates $x=z_j$, $t=z^{c_j}$ and coordinates on $N_\Sigma$, defined with $j$ in place of $k$, with $\Omega=d_j\,\dd\log t\wedge \dd\log x+\Omega_\Sigma$.
The mutation at $j$ replaces $x$ by $b=(1+t)/x$ times a monomial in $t$ and the coordinates on $N_\Sigma$, and fixes the other coordinates.
Both expressions in \eqref{eq:forms} are covariant under monomial changes of coordinates.
So \eqref{eq:mutation-potential}, with $d_j$ in place of $d_k$, gives $\ltor'-\ltor=-d_j\,\dd\mathcal D(-z^{c_j})$.
Finally, the complement of the common torus in either seed torus is a hypersurface, of real codimension two.
So restriction from the seed torus to the common torus is injective on $H^1(\,\cdot\,;\bR)$.
\end{proof}

\begin{proof}[Proof of Proposition~\ref{prop:lambda_ad}]

We prove by induction.
If $M=\emptyset$, then $\cA(\widetilde B)$ is the seed torus, and we choose \(\lambda_{ad} = p \dd \varphi\). 

Let $k$ be a sink of $B$. 
Put
\[
 V=\cA(\widetilde B_{[k]}) = \cA(\widetilde B)_x,\qquad U=\cA(\widetilde B)_t\cong\cN\times\Sigma,\qquad \kappa=d_k,
\]
with $G$ the torus generated by the circles $\xi_a$, $a\in F$. By induction we have good Liouville forms for \(V\) and \(\Sigma\). 
Let us verify Setup~\ref{setup:attachment}. 

\ref{A1} and \ref{A2} are evident. \ref{A3} follows from~\eqref{eq:form-splitting}.

\ref{A4}: The variable $x=z_k$ is frozen in $V$, so $Z_V\log|x|=\log|x|$ by \ref{P2} for $V$.
Since $k$ is frozen in $\widetilde B_{[k]}$, 
we can consider the corresponding circle action in~\eqref{eq:circle_action} on \(V = \cA(\widetilde B_{[k]})\). 
Recall \(\xi_k\) is the corresponding vector field.
Then $\iota_{\xi_k}\omega=d_k\,\dd\log|x|$ by \eqref{eq:frozen-hamiltonian}.
By \eqref{eq:dual-circle-identities}, $\epsilon_k(w_k)=0$ and $c_k(w_k)=1$, so $\xi_k$ fixes $x$ and rotates $t$ with weight one.
It fixes $u_j$ and $y_\alpha$, because their characters vanish on $w_k$.
So $\xi_k=\partial_\eta$, and $e^{i\alpha}\cdot(x,t,\sigma)=(x,e^{i\alpha}t,\sigma)$.

\ref{A5}: Let $a\in F$.
Since $c_k(w_a)=0$, the circle $\xi_a$ fixes $t$ and rotates $x$ with weight $\epsilon_k(w_a)$.
On $\Sigma$ it acts through the cocharacter $w_a-\epsilon_k(w_a)e_k\in N_\Sigma$.
This cocharacter satisfies $c_j(w_a-\epsilon_k(w_a)e_k)=-\epsilon_k(w_a)E_{kj}=0$ for $j\in I$, so it lies in the torus generated by the frozen circles of $\Sigma$.
By \ref{P2} and \ref{P4} for $V$ and $\Sigma$, these circles preserve $\lambda_V$, $\lambda_\Sigma$ and $\Psi_\Sigma$.
The form $\lN$ does not depend on $\theta=\arg x$, so $\lambda_U$ is $G$-invariant.
By \eqref{eq:frozen-pairing}, $\lambda_V(\xi_a)=-d_a\log|z_a|=\lambda_U(\xi_a)$.
For $\lambda_U$ we use $Z_U\log|z_a|=\log|z_a|$, which holds because $z_a$ is a monomial in $t$ and the $y_\alpha$.

\ref{A6} follows from \ref{P3} for $V$ and $\Sigma$.

\ref{A7}: By \eqref{eq:node-primitive}, $d_k\lN=\ltor^{(x,t)}+d_k\,\dd(h\rho\beta)$ with $\ltor^{(x,t)}=d_k(\log|t|\,\dd(\arg x)-\log|x|\,\dd(\arg t))$.
By \ref{P1} for $V$ and $\Sigma$ and by \eqref{eq:form-splitting}, $\lambda_U-\lambda_V$ is exact on the seed torus $T\subset V\cap U$.
The complement of $T$ in $V\cap U$ has real codimension two.
So $\lambda_U-\lambda_V$ is exact on $V\cap U$.

Now we apply Proposition~\ref{prop:gluing_Liouville_form}. 
Put $\lad=\lambda_X$ and $\Psi=\Psi_X$. 
By parts (i) and (v), $\lad$ is Morse--Bott with Lyapunov function $\Psi$.
We check \ref{P1}--\ref{P7}.

\ref{P1}: $\dd\lad=\omega=\Im\Omega$.
On the seed torus, $\lad-\ltor=(\lambda_V-\ltor)+\dd(\chi F)$, which is exact by \ref{P1} for $V$.
For the other seed tori, apply Lemma~\ref{lem:mutation-potential} along a sequence of mutations.

\ref{P2}: By Proposition~\ref{prop:gluing_Liouville_form}\ref{item:iv}, $\lad$ is $G$-invariant and $\lad(\xi_a)=\lambda_V(\xi_a)=-d_a\log|z_a|$ on $V$.
By \eqref{eq:frozen-hamiltonian}, $\lad(\xi_a)=-d_a\,Z_{\lad}\log|z_a|$.
So $Z_{\lad}\log|z_a|=\log|z_a|$ on the dense set $V$, hence on $\cA(\widetilde B)$.

\ref{P3}: This is Proposition~\ref{prop:gluing_Liouville_form}\ref{item:iii}.

\ref{P4}: By Proposition~\ref{prop:gluing_Liouville_form}\ref{item:v}, $\Psi$ is $G$-invariant.
Its critical components are those of $\Psi_V$ and the sets $\{x=y=0,\ t=-1\}\times C_\Sigma\cong C_\Sigma$, which are isotropic tori by induction.

\ref{P5}: Fix a seed torus, identified with $T^*\T^{m+f}$ by \eqref{eq:identify_with_torus}, and fix $\epsilon>0$.
Since $\lad|_{L_+}=\ltor|_{L_+}=0$, \ref{P1} gives $\lad-\ltor=\dd G_0$ on the seed torus, with $G_0|_{L_+}=0$.
Let $\chi_0$ be a cutoff with $\supp\chi_0\subset D^*_{2\epsilon}\T^{m+f}$ and $\chi_0=1$ on $D^*_\epsilon\T^{m+f}$, and put
\[
 \lambda_s=\lad-s\,\dd(\chi_0G_0),\qquad0\le s\le1 .
\]
Then $\lambda_1|_{D^*_\epsilon\T^{m+f}}=p\,\dd\varphi$ and $\lambda_s|_{L_+}=0$.
The support of $\chi_0G_0$ is compact, so $\lambda_s$ is a compact exact deformation.

\ref{P6} and
\ref{P7}: Follows from Proposition~\ref{prop:gluing_Liouville_form}\ref{item:vii}. 
\end{proof}

\subsection{The full rank case}\label{sec:fullrank-geometry}

We prove Proposition~\ref{prop:fullrank-geometry}.
We keep $P$, $\lambda_P$, $L_+$, $F$, $S$ and the action of $H=(\bC^*)^m$ on $P$ from Section~\ref{sec:geometry}.
We identify $F^{-1}(1)$ and split a finite cover of $P$, then deform its Liouville form to the product form.

Since $(P,\lambda_P)=(\cA_{\mathrm{prin}}(B),\lad)\times(T^*\T^f,p\,\dd q)$, \ref{Q1} follows from Corollary~\ref{cor:principal-geometry}, with
\begin{equation}\label{eq:PsiP}
 \Psi_P=\Psi+\sum_{a\in F}(\log|U_a|)^2,\qquad \Crit\Psi_P=\Crit\Psi\times\T^f .
\end{equation}

Next we write out the action of $H$.
For $h\in H$ and $D=\diag(d_i)_{i\in M}$,
\begin{equation}\label{eq:h-action}
 h\cdot(X_i,Z_i,U_a)=\Bigl(h_i^{d_i^{-1}}X_i,\ \prod_jh_j^{(D^{-1}B)_{ij}}Z_i,\ \prod_jh_j^{-C_{aj}}U_a\Bigr).
\end{equation}
The exchange character $\prod_iX_i^{B_{ik}}Z_k$ of $\cA_{\mathrm{prin}}(B)$ at $k$ has weight $(D^{-1}B+B^{\T}D^{-1})e_k=0$, by \eqref{eq:symmetrizer}.
So \eqref{eq:h-action} preserves the exchange relations.
By \eqref{eq:frozen-map},
\[
 F(h\cdot z)=h^SF(z),\qquad (h^S)_i=\prod_jh_j^{S_{ij}} .
\]
For the completion \eqref{eq:principal-completion}, the frozen cocharacters are $w_{m+j}=(E^{\T})^{-1}e_{m+j}=(-e_j;B^{\T}e_j)$.
For $h=e^{i\phi}$ with $\phi\in\bR^m$, \eqref{eq:h-action} acts on $\cA_{\mathrm{prin}}(B)$ by $z\mapsto e^{iw}z$ with $w=-\sum_j(D^{-1}\phi)_jw_{m+j}$, and by $U_a\mapsto e^{-i(C\phi)_a}U_a$.
Hence $\{|h|=1\}$ preserves both $\lambda_P$ and $\Psi_P$, by \ref{P2} and \ref{P4}.

Let $\widehat P$ be the fiber product of $F$ and $h\mapsto h^S$:
\[
 \widehat P=\{(z,h)\in P\times(\bC^*)^m:F(z)=h^S\},\qquad \pi(z,h)=z,\qquad \lambda_{\widehat P}=\pi^*\lambda_P .
\]
The map $\pi$ is a finite cover with deck group $\Delta=\bZ^m/S\bZ^m$.

\begin{lemma}\label{lem:phase-cover}
Let
\[
 \iota:\cA(\widetilde B)\longrightarrow P,\qquad \iota(z_M;z_F)=(z_M;z_F^C;z_F),\qquad (z_F^C)_i=\prod_{a\in F}z_a^{C_{ai}}.
\]
\begin{enumerate}[label=\textup{(\roman*)}]
\item The map $\iota$ is an isomorphism onto $F^{-1}(1)$.
\end{enumerate}
Identify $\cA(\widetilde B)$ with $F^{-1}(1)$ by $\iota$, and put $\alpha=\lambda_P|_{\cA(\widetilde B)}$ and $L_{\widetilde B,+}=L_+\cap\cA(\widetilde B)$.
\begin{enumerate}[label=\textup{(\roman*)},resume]
\item There is a diffeomorphism $\Phi:\cA(\widetilde B)\times T^*\T^m\to\widehat P$ with
\[
 \Phi^*\lambda_{\widehat P}=\alpha+p\,\dd q+\dd G(x,p),\qquad G(x,0)=0 .
\]
The deck group of $\pi$ is $\Delta=\{h\mapsto e^{iq_g}h\}$, with $q_g=2\pi S^{-1}g$ for $g\in\bZ^m/S\bZ^m$.
In the coordinates $(x,q,p)$, it acts by $(x,q,p)\mapsto(e^{-iq_g}\cdot x,\,q+q_g,\,p)$.
For $v\in\bZ^m$, the loop $\gamma_v(t)=\pi\Phi(x,q+2\pi tv,p)$, $t\in[0,1]$, satisfies $[F\circ\gamma_v]=Sv\in\pi_1\bigl((\bC^*)^m\bigr)$.
\item $(\cA(\widetilde B),\alpha)$ is a Liouville manifold with Morse--Bott Lyapunov function $\psi=\Psi_P|_{\cA(\widetilde B)}$, and
\[
 \pi^{-1}\Skel(Z_{\lambda_P})=\Phi\bigl(\Skel(Z_\alpha)\times\T^m\times\{0\}\bigr).
\]
\item We have $\Phi^{-1}\pi^{-1}L_+=\bigsqcup_{g\in\Delta}L_g$ with
\[
 L_g=L_{\widetilde B,g}\times T^*_{q_g}\T^m,\qquad L_{\widetilde B,g}=e^{-iq_g}\cdot L_{\widetilde B,+} .
\]
Each $L_{\widetilde B,g}$ is properly embedded with $\alpha|_{L_{\widetilde B,g}}=0$, and $G|_{L_g}=0$.
The deck group preserves $\alpha+p\,\dd q$ and $G$.
\end{enumerate}
\end{lemma}

\begin{proof}
(i) Consider the map $H\times\cA(\widetilde B)\to P$, $(h,x)\mapsto h\cdot\iota(x)$.
In initial coordinates it is $X_i=h_i^{d_i^{-1}}z_i$, together with the homomorphism of frozen tori
\[
 \binom{\log Z}{\log U}=R\binom{\log h}{\log z_F},\qquad
 R=\begin{pmatrix}D^{-1}B&C^{\T}\\-C&I_f\end{pmatrix},\qquad \det R=\det S\ne0 .
\]
We apply Lemma~\ref{lem:coefficient-change} to the seed $(B;I_m;0)$ of $P$, with new frozen coordinates $(h,z_F)$.
The matrices of the lemma are $\binom{D^{-1}}{0}$ and $R^{\T}$, in place of its $P$ and $Q$.
The new frozen rows are
\[
 \binom{D^{-1}}{0}B+R^{\T}\binom{I_m}{0}=\binom{D^{-1}B+B^{\T}D^{-1}}{C}=\binom{0}{C}.
\]
So the new seed is $(B;0;C)$, and $\cA(B;0;C)=H\times\cA(\widetilde B)$, because $h$ appears in no exchange relation.
Hence $H\times\cA(\widetilde B)$ is the pullback of $P$ along the isogeny $(h,z_F)\mapsto(Z,U)$, with projection $(h,x)\mapsto h\cdot\iota(x)$.
The kernel is $\{(\eta,\eta^C):\eta^S=1\}$ and acts by $(h,x)\mapsto(\eta h,\eta^{-1}\cdot x)$.
Since $F(h\cdot\iota(x))=h^S$,
\[
 F^{-1}(1)\cong\bigl(\{h^S=1\}\times\cA(\widetilde B)\bigr)/\{\eta^S=1\}\cong\cA(\widetilde B),
\]
where the inverse of the last isomorphism is induced by $x\mapsto(1,x)$.
The fiber is smooth: at $z\in F^{-1}(1)$, $\left.\frac{\dd}{\dd t}\right|_{t=0}F(e^{tv}\cdot z)=Sv$ for $v\in\bC^m$, and $S$ is invertible.

(ii) We construct $\Phi$.
For $s,q\in\bR^m$, let $\varrho_{s,q}$ act on $P$ by
\begin{equation}\label{eq:phase-coordinates}
 \varrho_{s,q}(X_i,Z_i,U_a)=\Bigl(e^{d_i^{-1}(s_i+iq_i)}X_i,\ e^{(D^{-1}B(s+iq))_i}Z_i,\ e^{(C(s-iq))_a}U_a\Bigr).
\end{equation}
On $(X;Z)$ this is the action \eqref{eq:h-action} of $h=e^{s+iq}$, and $U$ is a product factor of $P$.
So $\varrho_{s,q}$ is an automorphism of $P$, and $\varrho_{0,\theta}$ is the action of $e^{i\theta}$.
Since $D^{-1}$, $D^{-1}B$ and $C$ are integral, $\varrho_{s,q+2\pi v}=\varrho_{s,q}$ for $v\in\bZ^m$.
By \eqref{eq:frozen-map} and $-S^{\T}=D^{-1}B-C^{\T}C$,
\[
 F\circ\varrho_{s,q}=e^{-S^{\T}s+iSq}\,F .
\]
Put $s=-(S^{\T})^{-1}p$ and
\[
 \Phi(x,q,p)=\bigl(\varrho_{s,q}(x),\ e^{S^{-1}p+iq}\bigr)\qquad\bigl((x,q,p)\in\cA(\widetilde B)\times T^*\T^m\bigr).
\]
Then $F(\varrho_{s,q}(x))=e^{p+iSq}=(e^{S^{-1}p+iq})^S$, so $\Phi$ maps into $\widehat P$.
In particular $F\circ\gamma_v(t)=e^{p+iS(q+2\pi tv)}$, so $[F\circ\gamma_v]=Sv$.
Its inverse is $(z,h)\mapsto\bigl(\varrho_{s,q}^{-1}(z),\arg h,S\log|h|\bigr)$.
So $\Phi$ is a diffeomorphism, and $\pi\circ\Phi(x,q,p)=\varrho_{s,q}(x)$.
The deck transformations of $\pi$ are $h\mapsto\zeta h$ with $\zeta^S=1$, that is, $\zeta=e^{iq_g}$.
Since $\varrho_{s,q+q_g}=\varrho_{s,q}\circ\varrho_{0,q_g}$, they act by $(x,q,p)\mapsto(e^{-iq_g}\cdot x,q+q_g,p)$.

Next we compute $\Phi^*\dd\lambda_{\widehat P}$, where $\lambda_{\widehat P}=\pi^*\lambda_P$.
On the initial torus write $z_j=\exp(u_j+i\theta_j)$, as in \eqref{eq:forms}.
By \eqref{eq:phase-coordinates},
\[
 u_X=D^{-1}s+u_M,\quad \theta_X=D^{-1}q+\theta_M,\quad
 u_Z=D^{-1}Bs+C^{\T}u_F,\quad \theta_Z=D^{-1}Bq+C^{\T}\theta_F,\]\[
 u_U=Cs+u_F,\quad \theta_U=-Cq+\theta_F .
\]
Substitute this into \eqref{eq:principal-form}.
The coefficients of $\dd s\wedge \dd\theta_M$, $\dd s\wedge \dd\theta_F$, $\dd u_M\wedge \dd q$ and $\dd u_F\wedge \dd q$ are $D^{-1}(BD+DB^{\T})=0$, $C^{\T}-C^{\T}=0$, $B-B=0$ and $C-C=0$.
The coefficient of $\dd s\wedge \dd q$ is $B^{\T}D^{-1}-C^{\T}C=-S$.
With $p=-S^{\T}s$, this gives
\[
 \Phi^*\dd\lambda_{\widehat P}=\dd u^{\T}\begin{pmatrix}BD&-DC^{\T}\\CD&I_f\end{pmatrix}\wedge \dd\theta-\dd s^{\T}S\wedge \dd q=\dd\alpha+\dd p\wedge \dd q .
\]
The identity holds everywhere by density, and shows that $\dd\alpha$ is nondegenerate.

Put $\lambda=\Phi^*\lambda_{\widehat P}$, and let $Y=Z_\lambda$ be its Liouville field.
We have $\pi\Phi(x,q+\theta,p)=e^{i\theta}\cdot\pi\Phi(x,q,p)$, and $e^{i\theta}$ preserves $\lambda_P$.
So $\lambda$ is invariant under the translations of $q$.
Moreover $p_i=\log|F_i\circ\pi\Phi|$, and $F_i$ is a frozen monomial, so $Yp_i=p_i$ by \eqref{eq:principal-form}.
Hence
\[
 \lambda(\partial_{q_i})=\dd\lambda(Y,\partial_{q_i})=-(\iota_{\partial_{q_i}}\dd\lambda)(Y)=\dd p_i(Y)=p_i .
\]
So $\beta=\lambda-\alpha-p\,\dd q$ is closed, invariant under the translations of $q$, and satisfies $\beta(\partial_{q_i})=0$.
Thus $\beta$ is pulled back from $\cA(\widetilde B)\times\bR^m_p$.
Its restriction to $\{q=0,\ p=0\}$ vanishes, because $\pi\Phi(x,0,0)=x$.
The radial homotopy in $p$ gives $\beta=\dd G(x,p)$ with $G(x,0)=0$.

(iii) Write $Y=(Y_x,Y_q,Y_p)$.
By (ii), $\iota_Y(\dd\alpha+\dd p\wedge \dd q)=\alpha+p\,\dd q+\dd G$, that is,
\[
 \iota_{Y_x}\dd\alpha=\alpha+\dd_xG,\qquad Y_q=-\partial_pG,\qquad Y_p=p .
\]
At $p=0$ we have $G=0$ and $\dd_xG=0$, so $Y=(Z_\alpha,-\partial_pG(x,0),0)$ there.
The field $Y$ is complete, as the lift of $Z_{\lambda_P}$.
So $Z_\alpha$ is complete.
An orbit with $p\ne0$ leaves every compact set, since $p(t)=e^tp(0)$.
An orbit in $\{p=0\}$ is relatively compact if and only if its projection to $\cA(\widetilde B)$ is, since $\T^m$ is compact.
So $\Skel(Y)=\Skel(Z_\alpha)\times\T^m\times\{0\}$, and $\pi^{-1}\Skel(Z_{\lambda_P})=\Phi(\Skel(Y))$ because $\pi$ is a finite cover.
Put $\widetilde\Psi_P=\Psi_P\circ\pi\circ\Phi$ and $\psi=\Psi_P|_{\cA(\widetilde B)}$.
The function $\widetilde\Psi_P$ is Morse--Bott and independent of $q$, with $\widetilde\Psi_P(x,q,0)=\psi(x)$.
Thus
\[
 \dd\psi(Z_\alpha)=\dd\widetilde\Psi_P(Y)|_{p=0},\qquad
 \Crit\widetilde\Psi_P=\Crit\psi\times\T^m\times\{0\},\qquad \Crit\psi=\Zero(Z_\alpha).
\]
Indeed, if $\dd\psi=0$ or $Z_\alpha=0$, the first identity forces $Y=0$; conversely, $Y=0$ implies $p=0$ and $\dd\widetilde\Psi_P=0$.
Hence $\Crit\psi=\Crit\widetilde\Psi_P/\T^m$ is compact and smooth.
The function $\psi$ is proper, and restricting the gradient-like estimate for $\widetilde\Psi_P$ gives that for $\psi$.

At a critical point, put $\mathcal H=\Hess\widetilde\Psi_P$.
The quadratic part of the gradient-like estimate gives
\[
 \mathcal H(v,DY\,v)\ge c\bigl(\norm{DY\,v}^2+\norm{\mathcal Hv}^2\bigr).
\]
If $\Hess\psi(w,-)=0$ and $v=(w,0,0)$, then $\mathcal H(v,-)$ has only a $p$-component, whereas $DY\,v$ has no $p$-component.
Thus $\mathcal H(v,DY\,v)=0$, hence $\mathcal Hv=0$ and $w\in T\Crit\psi$.
This proves that $\psi$ is Morse--Bott.
For a component $C$ of $\Crit\psi$, the component $C\times\T^m\times\{0\}$ of $\Crit\widetilde\Psi_P$ is isotropic.
Since $\Hess\psi$ is a restriction of $\mathcal H$,
\[
 \ind_C\psi+\dim C
 \le \ind_{C\times\T^m\times\{0\}}\widetilde\Psi_P+\dim C
 \le\tfrac12\dim_\bR\cA(\widetilde B).
\]
Also $\dd\alpha|_{TC}=0$, since $\dd\lambda=\dd\alpha+\dd p\wedge\dd q$.
Every backward $Z_\alpha$-orbit enters a sufficiently large regular sublevel of $\psi$, by properness and $\dd\psi(Z_\alpha)>0$ outside $\Crit\psi$.
So $\alpha$ is of finite type.

(iv) A point $\Phi(x,q,p)$ lies over $L_+$ if and only if $\varrho_{s,q}(x)\in L_+$.
Then $e^{p+iSq}=F(\varrho_{s,q}(x))>0$, so $Sq\in2\pi\bZ^m$ and $q=q_g$ for some $g$.
The map $\varrho_{s,0}$ preserves $L_+$, and $\varrho_{0,q_g}$ is the action of $e^{iq_g}$.
So the condition reads $x\in e^{-iq_g}\cdot L_{\widetilde B,+}=L_{\widetilde B,g}$, and $\Phi^{-1}\pi^{-1}L_+=\bigsqcup_gL_g$.
The element $e^{-iq_g}$ preserves $\cA(\widetilde B)$, since $(e^{-iq_g})^S=1$, and it preserves $\lambda_P$.
So $L_{\widetilde B,g}$ is properly embedded and $\alpha|_{L_{\widetilde B,g}}=0$.
By (ii), the deck group preserves $\alpha+p\,\dd q$.
It preserves $\lambda$, as a group of deck transformations.
Hence it preserves $\dd G$, and also $G$, since $G=0$ on $\{p=0\}$.
Finally $\dd G|_{L_g}=\lambda|_{L_g}-(\alpha+p\,\dd q)|_{L_g}=0$, because $\lambda_P|_{L_+}=0$, $\alpha|_{L_{\widetilde B,g}}=0$ and $q=q_g$ on $L_g$.
Since $L_g$ is connected and meets $\{p=0\}$, we get $G|_{L_g}=0$.
\end{proof}

The correction $\dd G$ need not have compact support.
Put $\lambda_0=\alpha+p\,\dd q$ on $\cA(\widetilde B)\times T^*\T^m$.

\begin{lemma}\label{lem:phase-product}
There are a compact exact deformation $\lambda_P+\dd H_\tau$, $0\le\tau\le1$, with $H_0=0$ and $H_\tau|_{L_+}=0$, and a $\Delta$-equivariant diffeomorphism $\Theta:\cA(\widetilde B)\times T^*\T^m\to\widehat P$ such that
\[
 \Theta^*\pi^*(\lambda_P+\dd H_1)=\lambda_0,\qquad \Theta(L_g)=\Phi(L_g).
\]
Moreover $\Theta=\Phi$ on a Liouville domain $D_0$ of $\lambda_0$, and $\Theta_*=\Phi_*$ on $\pi_1(\T^m)$.
\end{lemma}

\begin{proof}
We first construct the deformation.
Choose a $\Delta$-invariant $\chi\in C^\infty_c(\cA(\widetilde B)\times\bR^m_p,[0,1])$ equal to $1$ near $\Skel(Z_\alpha)\times\{0\}$.
Since $G$ is $\Delta$-invariant, the formula
\[
 H_\tau\circ\pi\circ\Phi=-\tau\chi G
\]
defines $H_\tau$ on $P$, supported in the compact image of $\supp\chi\times\T^m$.
Also $H_\tau|_{L_+}=0$, since $G|_{L_g}=0$.
Put $\widehat\lambda_\tau=\pi^*(\lambda_P+\dd H_\tau)$.
Then
\[
 \Phi^*\widehat\lambda_\tau=\lambda_0+\dd\bigl((1-\tau\chi)G\bigr).
\]
These forms agree with $\Phi^*\lambda_{\widehat P}$ outside a fixed compact set, so their Liouville fields are complete and have a common outer Liouville domain.
Since $(1-\tau\chi)G$ is independent of $q$ and vanishes at $p=0$, the calculation in Lemma~\ref{lem:phase-cover}(iii) gives
\[
 Z_{\Phi^*\widehat\lambda_\tau}p=p,\qquad
 \Skel(Z_{\Phi^*\widehat\lambda_\tau})=\Skel(Z_\alpha)\times\T^m\times\{0\}.
\]

Next we construct $\Theta$.
The function $\psi=\Psi_P|_{\cA(\widetilde B)}$ is $\Delta$-invariant.
Take a sufficiently large regular sublevel $D$ of $\psi+|p|^2$, a $\Delta$-invariant Liouville domain for $\lambda_0$.
Let $\phi_0$ and $\phi_1$ be the flows of $Z_{\lambda_0}$ and $Z_{\Phi^*\widehat\lambda_1}$.
The domains $\phi_0^{-T}(D)$, $T\ge0$, decrease to $\Skel(Z_\alpha)\times\T^m\times\{0\}$.
So $D_0=\phi_0^{-T}(D)$ lies in $\{\chi=1\}$ for large $T$, and $\Phi^*\widehat\lambda_1=\lambda_0$ near $D_0$.
Every backward $\phi_1$-orbit enters $D_0$.
Indeed, a backward orbit is relatively compact once it enters the common outer Liouville domain.
If it avoided $\operatorname{int}D_0$, its backward limit set would be a compact invariant set disjoint from $\Skel(Z_\alpha)\times\T^m\times\{0\}\subset\operatorname{int}D_0$, which is impossible.
So $(\cA(\widetilde B)\times T^*\T^m,\lambda_0)$ and $(\cA(\widetilde B)\times T^*\T^m,\Phi^*\widehat\lambda_1)$ are both completions of $(D_0,\lambda_0)$.
Put
\[
 \Theta=\Phi\circ\Theta_1,\qquad \Theta_1|_{D_0}=\id,\qquad \Theta_1(\phi_0^\sigma z)=\phi_1^\sigma z\quad(z\in\partial D_0,\ \sigma\ge0).
\]
Then $\Theta_1^*\Phi^*\widehat\lambda_1=\lambda_0$, that is, $\Theta^*\widehat\lambda_1=\lambda_0$.
Both forms vanish on $L_g$, so both flows preserve $L_g$ and $\Theta(L_g)=\Phi(L_g)$.
The forms and $D_0$ are $\Delta$-invariant, so $\Theta_1$ commutes with $\Delta$.
Finally, $D_0\hookrightarrow\cA(\widetilde B)\times T^*\T^m$ is a homotopy equivalence and $\Theta|_{D_0}=\Phi|_{D_0}$, so $\Theta_*=\Phi_*$.
\end{proof}

\begin{proof}[Proof of Proposition~\ref{prop:fullrank-geometry}]
\ref{Q1} follows from Corollary~\ref{cor:principal-geometry} and the product decomposition of $P$.
\ref{Q2} is Lemma~\ref{lem:phase-cover}(i).
\ref{Q3} is Lemma~\ref{lem:phase-cover}(iii), (iv).

For \ref{Q4}, take $H_\tau$ and $\Theta$ from Lemma~\ref{lem:phase-product}, and put
\[
 \Pi(x,q,p)=\pi\Theta(x,q\bmod2\pi\bZ^m,p).
\]
The $\Delta$-equivariance of $\Theta$ identifies the deck group of $\Pi$ with the action of $\Gamma$ in~\ref{Q4}.
The two pullback identities follow from Lemma~\ref{lem:phase-product} and Lemma~\ref{lem:phase-cover}(iv).
\end{proof}

\section{Mirror symmetry}

In this section we could work throughout with Fukaya categories with $\bZ$ coefficients, and correspondingly on the $B$-side with algebraic varieties over $\bZ$.  We write $\C$ instead for essentially typographical reasons. 

\subsection{B-side}

Let $S$  be  a smooth finitely generated commutative $\bC$-algebra, and 
suppose $A \subset S[u^{\pm 1}, y^{\pm 1}]$ is a smooth finitely generated $\C$-subalgebra such that $y \in A^\times$ and $S[u^{\pm 1}, y^{\pm 1}] = A[g^{-1}]$ for some $g \in A \cap S[u^{\pm 1}, y]$.

We write: 
\[ Y=\Spec A\cup_{\Spec S[u^{\pm 1}, y^{\pm 1}]} \Spec S[u^{\pm 1}, y].
\]
\begin{lemma}
\label{lem:Y}
$Y$ is a smooth separated scheme of finite type, and it is quasi-affine.
\end{lemma}

\begin{proof}
$Y$ is smooth because both charts are.
It is separated because $A\otimes_{\bC}S[u^{\pm1}, y]\to S[u^{\pm1}, y^{\pm 1}]$ is surjective, as $y^{-1}\in A$.
Choose finite generating sets $(a_i)$ of $S[u^{\pm1}, y]$ and $(b_j)$ of $A$.
Since $S[u^{\pm1}, y^{\pm 1}]=A[g^{-1}]=S[u^{\pm1}, y][y^{-1}]$, there are exponents with $g^{N_i}a_i\in A$ and $y^{M_j}b_j\in S[u^{\pm1}, y]$.
The ring $R_Y=\bC[y,g,g^{N_i}a_i,y^{M_j}b_j]\subset S[u^{\pm1}, y] \cap A$ satisfies $(R_Y)_g=S[u^{\pm1}, y]$ and $(R_Y)_y=A$.
So $Y\cong D(g)\cup D(y)\subset\Spec R_Y$ is quasi-affine.
\end{proof}

We set
\[
 i:D=\{y=0\}=\Spec S[u^{\pm1}]\hookrightarrow Y,\qquad
 j:H=\{u=-1\}=\Spec S\hookrightarrow D .
\]

\begin{lemma}
\label{lem:blowup}
Let $D'$ be the strict transform of $D$ in the blowup of $Y$ along $H$.
Then
\[
 \Bl_HY\setminus D'
 \cong\Spec A\cup_{\Spec S[u^{\pm 1}, y^{\pm 1}]}\Spec\frac{S[u^{\pm1}][y,y']}{(yy'-1-u)} .
\]
The blowup projection restricts to $\pi:\Bl_HY\setminus D'\to Y$, with $\pi|_{\Spec A}=\id$.
\end{lemma}
\begin{proof}
The center $H$ lies in $D\subset \Spec S[u^{\pm1}, y]$ and is disjoint from $\Spec A\setminus D$.
So over $\Spec A$ the modification is the identity.
Over $\Spec S[u^{\pm1}, y]$, the blowup of the ideal $(y,1+u)$ has two charts.
On the chart where $1+u$ generates the exceptional ideal, the strict transform of $D$ is $\{y/(1+u)=0\}$.
Off this locus, $y$ also generates the exceptional ideal.
So removing $D'$ leaves exactly the chart where $y$ generates, with coordinate $y'=(1+u)/y$.
That chart is $\Spec S[u^{\pm1}][y,y']/(yy'-1-u)$.
Its intersection with $\Spec A$ is 
$\Spec S[u^{\pm1},y^{\pm1}]$.
\end{proof}

\begin{proposition} 
\label{prop: general B-model diagram}
We have the following commutative diagram.
    \[
\begin{tikzcd}
& \Perf_A \arrow[r] & \Perf_{S[u^{\pm1},y^{\pm1}]}  \\
\Coh(\Bl_H Y \setminus D')  \arrow[dr, phantom, "\lrcorner", very near start] & \Coh(Y) \arrow[ur, phantom, "\urcorner", very near start] \arrow[l] \arrow[r] \arrow[u] & \Perf_{S[u^{\pm1}, y]} \arrow[u] \\
\Perf_S \arrow[u] & \Perf_{S[u^\pm]} \arrow[r, equal]  \arrow[l, "j^*"] \arrow[u, "i_*"] &  \Perf_{S[u^\pm]} \arrow[u]
\end{tikzcd}
\]
As indicated, the lower left square is a pushout and the upper right square is a pullback. 
\end{proposition}
\begin{proof}
    The lower left square being a pushout holds by \cite[Theorem~4.13]{GLe}, since $H\subset D \subset Y$ are inclusions of smooth divisors in smooth spaces.  The upper right square is a pullback by Zariski descent. 
\end{proof}

\begin{corollary}
    Fix a unimodular seed \( \widetilde B = \binom{B}{C}\). Let \(k\) be a sink and \(x\) be the corresponding cluster variable as in \eqref{eq:sink_exchange}. 
    Assume $B_{I,I}$ is locally acyclic (e.g. $B$ sink recurrent).
    Fix a completion \(E\) and the submatrices of $\widetilde{B}$ 
as in \eqref{eq:matrices_B}. Recall by Proposition \ref{prop:splitting_ABprime}  there are splittings \( \cA(\widetilde{B})_{xt} \cong  \bC^*_x \times \cA(\widetilde B_{\setminus k}) \) and \(\cA(\widetilde B_{\setminus k}) \cong \bC^*_t \times \cA(\widetilde B_\Sigma ) \). 

\begin{equation}\label{eq:s4-B-diagram}
\begin{tikzcd}[column sep=large,row sep=large]
&\Coh(\cA( \widetilde B_{[k]})) \arrow[r]
&\Coh( \bC_x^* \times  \bC^*_t \times \cA(\widetilde B_\Sigma ) )\\
\Coh(\mathcal{A}(\widetilde B)) \arrow[dr, phantom, "\lrcorner", very near start]
&\Coh(Y)\arrow[l]\arrow[u]\arrow[r] \arrow[ur, phantom, "\urcorner", very near start]
&\Coh(\bC_x \times  \bC^*_t \times \cA(\widetilde B_\Sigma )) \arrow[u]\\
\Coh( \cA(\widetilde B_{\Sigma}))  \arrow[u]
&\Coh(  \bC^*_t \times \cA(\widetilde B_\Sigma ))      \arrow[l,"t=-1"]
\arrow[u]\arrow[r,equal]
&\Coh(  \bC^*_t \times \cA(\widetilde B_\Sigma ) ) \arrow[u,"x=0"]
\end{tikzcd}
\end{equation}
\end{corollary}
\begin{proof}
There is a sink cover given by \eqref{eq:sink_cover}. 
Take \(A = \sA( \widetilde B_{[k] })\), \(S = \sA (\widetilde B_\Sigma )\), \(y = x\) and \(u = t\). 
By Proposition~\ref{prop:AB_Nk} and Proposition~\ref{prop:splitting_ABprime}, all the assumptions in Lemma~\ref{lem:Y} and Lemma~\ref{lem:blowup} are satisfied. In particular, Lemma~\ref{lem:blowup} and the sink cover~\eqref{eq:sink_cover} imply that \(\Bl_H Y \setminus D' \cong \cA(\widetilde B)\). The diagram follows from Proposition~\ref{prop: general B-model diagram}. 
\end{proof}

\subsection{A-side}

We take $\cA(\widetilde B)$ with the symplectic form $\Im\Omega$, and let $I$ be its complex structure.
The determinant of the holomorphic symplectic form gives $\det_\bC(T\cA(\widetilde B))$ a canonical trivialization.
We use its square to give Maslov data on $\cA(\widetilde B)$.
Another way to describe this is as follows. 
The complex symplectic bundle $(T\cA(\widetilde B),\Omega)$ has structure group $Sp_{m+f}(\bC)$.
Choose an almost complex structure $J$ on $\cA(\widetilde B)$ compatible with $\Im\Omega$ and with $IJ=-JI$.
Then $(T\cA(\widetilde B),J)$ has structure group $Sp_{m+f}(\bC)\cap U(m+f)$, which is compact and simply connected.
So the determinant map to $U(1)$ is trivial on it, giving $\det_\bC(T\cA(\widetilde B),J)$ a trivialization.
We take the trivial background class in $H^2(\cA(\widetilde B);\bZ/2)$, so to work over $\bZ$, we should equip
Lagrangian branes with spin structures.

On a seed torus $(\bC^*)^{m+f}$, we can take $J$ constant in the frame $(\partial_u,\partial_\theta)$.
The fiber polarization $\operatorname{span}\{\partial_{p_j}\}=\operatorname{span}\{\partial_{u_j}\}$ of $T^*\T^{m+f}$ under \eqref{eq:identify_with_torus} is also constant in this frame.
So on the seed torus our Maslov data is the one given by the fiber polarization.
The central torus $\T^{m+f}$ is transverse to the fibers, and $L_+=\{\varphi=0\}$ is a fiber.
We grade both by their constant phase relative to the fiber polarization.
We give $\T^{m+f}$ the Lie group spin structure, i.e.\ the spin structure of the frame $(\partial_{\varphi_j})$.
Its restriction to each circle $\{|z_i|=1\}$, with the other $z_j$ fixed, does not extend over a disk.
The Lagrangian $L_+\cong\bR^{m+f}$ has a unique spin structure.

In the full rank case we make the same choices on $P$, and pull them back along the covering $\Pi$ of~\ref{Q4}.

\begin{theorem}
       Fix a unimodular, sink recurrent seed \( \widetilde B = \binom{B}{C}\). Let \(k\) be a sink and \(x\) be the corresponding cluster variable as in \eqref{eq:sink_exchange}. 
There is a sink cover given by \eqref{eq:sink_cover}.
Fix a completion \(E\) and the submatrices of $\widetilde{B}$ 
as in \eqref{eq:matrices_B}. Recall by Corollary~\ref{cor:At-splitting} there are splittings \( \cA(\widetilde{B})_{xt} \cong  \bC^*_x \times \cA(\widetilde B_{\setminus k}) \) and \(\cA(\widetilde B_{\setminus k}) \cong \bC^*_t \times \cA(\widetilde B_\Sigma ) \). Use the adapted Liouville form \(\lambda_{ad}\) in Proposition~\ref{prop:lambda_ad}, then we have the following diagram. 

\footnotesize
\begin{equation}\label{eq:s4-A-diagram}
\begin{tikzcd}[%
row sep=large]
&\Fuk\left(\cA( \widetilde B_{[k]})\right)\arrow[r]
&\Fuk\left(  \bC_x^* \times  \bC^*_t \times \cA(\widetilde B_\Sigma )  \right)\\
\Fuk\left( \cA(\widetilde B) \right) \arrow[dr, phantom, "\lrcorner", very near start]
&
\Fuk\left( \cA( \widetilde B_{[k]}), T^*S^1 \times \cA(\widetilde B_\Sigma) \right) \arrow[l]\arrow[u]  
\arrow[ur, phantom, "\urcorner", very near start]
\arrow[r,"q"]
&\Fuk \left(\bC_x^* \times  \bC^*_t \times \cA(\widetilde B_\Sigma ) , T^*S^1 \times \cA(\widetilde B_\Sigma) \right) \arrow[u]\\
\Fuk\left(T^*D^2 \times \cA(\widetilde B_\Sigma)\right) \arrow[u]
&\Fuk\left(T^*S^1 \times \cA(\widetilde B_\Sigma)\right)\arrow[l]
\arrow[u,"J"]\arrow[r,equal]
&\Fuk\left(T^*S^1 \times \cA(\widetilde B_\Sigma)\right)\arrow[u,"J_0"']
\end{tikzcd}
\end{equation}
\normalsize

\end{theorem}

\begin{proof}
    Each of the maps above is either a Viterbo restriction, stop removal, or Orlov functor.  The existence of Liouville forms adapted to the definitions of such maps is established in Proposition~\ref{prop:lambda_ad}.
    By \ref{P6}, the lower left square is a sectorial gluing, hence it is a pushout diagram by \cite{GPS2}. 
    By \ref{P7}
    the map 
    \( \Fuk\left( \cA(\widetilde B_{[k]}) \right) \to \Fuk(\bC_x^* \times  \bC^*_t \times \cA(\widetilde B_\Sigma)) \) and the map \(q\) can be defined as Viterbo restrictions, and
    \( V, O, \cF\) satisfy the assumptions in Lemma~\ref{lem:circle-twist} and Theorem~\ref{thm: SOV descent}. Thus
    the upper right square is a pullback. 
\end{proof}

Also in Proposition \ref{prop:lambda_ad}, we showed that the central torus in each cluster chart is exact for $\lambda_{ad}$, and we moreover chose fixed parameterizations by $(S^1)^{m+f}$ by 
\eqref{eq:identify_with_torus}. 

We recall: 

\begin{theorem} \label{thm: compact cluster}
    \cite{PT, GLe} For a full rank sink recurrent seed $\widetilde{B}$, the collection of parameterized exact Lagrangian tori 
    is a Lagrangian cluster ensemble. 
\end{theorem}

\begin{remark}
    Strictly speaking, one should check that the Liouville structures / skeleta used in the cited references are compatible with those in our setting, at least in a neighborhood of a pair of such tori; this can be done similarly to \cite[Sec. 4.7]{PT}.  Alternatively, one could just rewrite the proofs in our setting.  We leave this to the reader.
\end{remark}

Finally, we will need: 

\begin{lemma}\label{lem: viterbo respects positive}
    The Viterbo restriction $\Fuk(\cA( \widetilde B_{[k]})) 
     \to \Fuk( \bC_x^* \times  \bC^*_t \times \cA(\widetilde B_\Sigma ) )$ carries the positive real locus $L_+$ to the positive real locus $L_+$.  The same holds for the map $q: \Fuk\left( \cA( \widetilde B_{[k]}), T^*S^1 \times \cA(\widetilde B_\Sigma) \right) \to \Fuk \left(\bC_x^* \times  \bC^*_t \times \cA(\widetilde B_\Sigma ) , T^*S^1 \times \cA(\widetilde B_\Sigma) \right)$.
\end{lemma}
\begin{proof}
    By \ref{P7},  we have an embedding of subdomains \( \overline{\bC_x^* \times  \bC^*_t \times \cA(\widetilde B_\Sigma )} \subset \overline{\cA( \widetilde B_{[k]})} \), and by \ref{P3} \(L_+\) is strictly exact.  The intersection of \(L_+\) for the larger domain with the smaller domain is still the real loci Lagrangian. 
    The claim follows from~\cite[Section 2.6]{Sy}, which proves that the \cite{GPS2} Viterbo functor coincides with that of \cite[Section 4.11]{AbouzaidSeidel} for such Lagrangians. 
\end{proof}

\subsection{Model case}

Here we consider the example: $$X = \cA \binom{0}{1} = \{ x y = 1 + t\} \subset \bC_x \times \bC_y \times \bC^*_t.$$ Then \(x\) is a cluster variable corresponding to a sink. 

On the B-side, \eqref{eq:s4-B-diagram} becomes
\begin{equation}
\label{eq:node-B-diagram}
\begin{tikzcd}[column sep=large,row sep=large]
&\Coh(\bC_{\c x}^*\times\bC_{\c t}^*) \arrow[r, equal]
&\Coh( \bC_{\c x}^*\times\bC_{\c t}^* )\\
\Coh(X) \arrow[dr, phantom, "\lrcorner", very near start]
&\Coh(\bC_{\c x}\times\bC_{\c t}^*)\arrow[l]\arrow[u]\arrow[r, equal] 
&\Coh(\bC_{\c x}\times\bC_{\c t}^*) \arrow[u]\\
\Coh(\mathrm{pt})  \arrow[u]
&\Coh(  \bC_{\c t}^*)      \arrow[l,"\c t=-1"]
\arrow[u]\arrow[r,equal]
&\Coh(  \bC_{\c t}^* ) \arrow[u,"\c x=0"]
\end{tikzcd}
\end{equation}

On the A-side, diagram~\eqref{eq:s4-A-diagram} reads
\begin{equation}\label{eq:node-A-diagram}
\begin{tikzcd}[%
row sep=large]
&\Fuk\left(\bC_x^*\times\bC_t^*\right)\arrow[r, equal]
&\Fuk\left(  \bC_x^*\times\bC_t^*  \right)\\
\Fuk\left( X \right) \arrow[dr, phantom, "\lrcorner", very near start]
&
\Fuk\left( \bC_x^*\times\bC_t^*, T^*S^1 \right) \arrow[l]\arrow[u]
\arrow[r, equal]
&\Fuk \left(\bC_x^*\times\bC_t^* , T^*S^1 \right) \arrow[u]\\
\Fuk\left(T^*D^2\right) \arrow[u]
&\Fuk\left(T^*S^1\right)\arrow[l]
\arrow[u,"J"]\arrow[r,equal]
&\Fuk\left(T^*S^1\right)\arrow[u,"J_0"']
\end{tikzcd}.
\end{equation}
Here \(T^*S^1\) is embedded in the ideal boundary of \(\bC^*_x \times \bC^*_t\) by \( \arg(-t) = 0, \log|x| = -\infty \). 

We describe explicitly how the two sides match. 
The completion  of \( \binom{0}{1} \) is \( \begin{pmatrix}
    0 & -1 \\
    1 & 0
\end{pmatrix}
\). Hence the map~\eqref{eq:identify_with_torus} is given by
\begin{align*}
    \bC^*_x \times \bC^*_t & \longrightarrow T^*(S^1_{\c x} \times S^1_{\c t}) \\
    (x, t) & \longmapsto ( -\arg t, \arg x, \log |x|,  \log |t| ).  
\end{align*}
In particular, the circle direction of \(x\) is sent to the circle direction of \(\c t\), and vice versa.  
This phenomenon happens because we use the symplectic structure \(\Im \Omega\) on the A-side. Given this identification, we have \( (\bC_x^* \times \bC_t^* , T^*S^1) \cong\big( T^* S^1_{\c x}, { (\pi, -\infty)} \big)  \times T^*S^1_{\c t}  \). 

Now we describe the maps matching two sides (which is a special case of \cite{GLe}). The isomorphism
\[
\Coh(\bC_{\c x}^* \times \bC_{\c t}^*)\xlongrightarrow{\sim}\Fuk( T^*S^1_{\c x} \times T^*S^1_{\c t} ) \cong \Fuk(\bC^*_x \times \bC^*_t)
\]
is the product of the standard identifications \( \Coh(\bC^*_{\c x}) \simeq \Fuk(T^*S^1_{\c x})\) and \(\Coh(\bC^*_{\c t} ) \simeq \Fuk(T^*S^1_{\c t}) \).  Similarly, the isomorphism 
\[
    \Coh( \bC_{\c x} \times \bC_{\c t}^*) \xlongrightarrow{\sim} \Fuk\big( (T^*S^1_{\c x}, (\pi, -\infty)) \times T^*S^1_{\c t}  \big) \cong \Fuk(\bC^*_x \times \bC^*_t , T^*S^1)  
\]
is given by \( \Coh ( \bC_{\c x}) \simeq \Fuk\big(T^*S^1_{\c x}, (\pi, -\infty) \big) \) and  \(\Coh(\bC^*_{\c t} ) \simeq \Fuk(T^*S^1_{\c t}) \). Finally, the map \(\Fuk(T^*S^1) \simeq \Loc(S^1) \longrightarrow \Fuk(T^*D^2)\) is the left adjoint of the functor \( \Loc(D^2) \longrightarrow \Loc(S^1)\), sending the trivial local system to local system with monodromy \(-1\). Here the minus sign comes because the spin structure naturally associated to the fiber polarization of the cotangent bundle $T^*(S^1 \times S^1)$ is the Lie group spin structure; its restriction to one of the $S^1$ gives the spin structure which does not extend over the disk. 

It matches with the B-side functor
\( \Coh(\bC_{\c t}^* ) \xrightarrow{\c t = -1} \Coh(\mathrm{pt}) \).

\subsection{Unimodular case}

Suppose given affine algebraic varieties $X, Y$, a map $\eta: \Coh(X) \to \Coh(Y)$, and an isomorphism $\eta (\mathcal{O}_X) \cong \mathcal{O}_Y$.  Then $\eta$ induces 
$\eta_{\cO}: \cO(X) = \Hom(\cO_X, \cO_X) \to \Hom(\cO_Y, \cO_Y) = \cO(Y)$, and evidently $\eta$ is the pullback along $\Spec(\eta_{\cO})$.  Note that $\eta_{\cO}$ is determined by its composition with the restriction of $\cO(Y)$ to any Zariski dense subset.

It follows that the  isomorphism as requested by the theorem is unique if it exists.  Indeed, in the unimodular case $\mathcal{X}(\widetilde B^\vee)$ is an affine variety, and we have fixed the image of $\mathcal{O}$ and required that the desired isomorphism is compatible with a prescribed isomorphism on each cluster chart (one would be enough).

Fix some unimodular sink recurrent cluster variety.  Choose a seed $\widetilde{B}$ with a sink for which the sink cover is by (unimodular) sink recurrent cluster varieties.  Fix a unimodular completion  of $\widetilde{B}$ to $(E,\widetilde D)$; we write  $E^\vee =-E^{\T}$ and take $\widetilde{B}^\vee$ to be the dual seed it completes.  
We will prove the theorem by comparing \eqref{eq:s4-A-diagram} for $\widetilde{B}$ with \eqref{eq:s4-B-diagram} for $\widetilde{B}^\vee$. 

Note by Corollary \ref{cor:dual-Sigma} we have $(\widetilde B^\vee)_\Sigma = (\widetilde B_\Sigma)^\vee$; also 
$(\widetilde B^\vee)_{[k]} = (\widetilde B_{[k]})^\vee$ by Lemma~\ref{lem:dual-freezing}.  We henceforth suppress the order of these operations. 

The seeds $\widetilde{B}_{[k]}$ and $\widetilde{B}_\Sigma$ have fewer mutable variables, so we may assume by induction that we have established the theorem for these. 
That is, we have isomorphisms $\Coh(\cA(\widetilde B^\vee_{[k]})) \simeq \Fuk(\cA(\widetilde B_{[k]}))$ and $\Coh(\cA(\widetilde B^\vee_{\Sigma}))  \simeq  \Fuk(\cA(\widetilde B_{\Sigma}))$, in each case matching $\mathcal{O}$ with the positive real locus, and inducing the expected cluster transformations on the moduli of objects arising from the compact tori.

The right column and bottom row of the diagrams \eqref{eq:s4-B-diagram} and \eqref{eq:s4-A-diagram}   are products of the corresponding parts of \eqref{eq:node-B-diagram} and \eqref{eq:node-A-diagram} with 
$\Coh(\cA(\widetilde B^\vee_{\Sigma}))$ and $\Fuk(\cA(\widetilde B_{\Sigma}))$
respectively; we induce compatible isomorphisms by induction and  K\"unneth. 

To match the diagrams it remains to show that the following diagram commutes: 
\begin{equation}\label{eq:compatibility check}
\begin{tikzcd}[column sep=large,row sep=large]
\Coh(\cA( \widetilde B_{[k]}^\vee)) \arrow[r] \arrow[d]
&\Coh( \bC_{\c x}^* \times  \bC^*_{\c t} \times \cA(\widetilde B_\Sigma^\vee ) ) \arrow[d]\\
\Fuk(\cA( \widetilde B_{[k]})) \arrow[r] \arrow[u]
&\Fuk( \bC_x^* \times  \bC^*_t \times \cA(\widetilde B_\Sigma ) ) \arrow[u]
\end{tikzcd}
\end{equation}
Indeed, having done so, the remaining vertices of the diagrams are obtained by pullback and pushout. 

Since the vertical arrows of \eqref{eq:compatibility check} match $L_+$ and $\mathcal{O}$ and the top horizontal arrow is the restriction of sheaves, it suffices to show (as explained at the outset of this proof) that the bottom horizontal arrow sends $L_+ \mapsto L_+$ (this is Lemma \ref{lem: viterbo respects positive} above) and that the induced map $\Hom_{\cA( \widetilde B_{[k]})}(L_+, L_+) \to \Hom_{\bC_x^* \times  \bC^*_t \times \cA(\widetilde B_\Sigma )}(L_+, L_+)$ is as expected.  But again, we may check this after further Viterbo restriction to a cluster torus, where the images match by the inductive hypothesis.  (Here we implicitly use functoriality of Viterbo restriction; see \cite[Sec. 2]{HS26} for a proof.) 

We have now matched the diagrams and hence produced an isomorphism 
$\Fuk(\cA( \widetilde B)) \cong \Coh(\cA( \widetilde B^\vee))$. It remains to establish the assertions regarding seed tori and the positive real locus.  For the first, some of the seed Lagrangian tori were already in $\cA( \widetilde B_{[k]})$, functoriality of Viterbo implies that restriction to these tori remains compatible as desired.  Compatibility with the remaining tori follows from the known mutation formula for compact Lagrangians, recalled in Theorem \ref{thm: compact cluster} above. 

Regarding the positive locus, it suffices to show that the identification $\Fuk( \cA( \widetilde B_{[k]}), T^*S^1 \times \cA(\widetilde B_\Sigma)) \cong \Coh(Y)$ sends $L_+ \mapsto \mathcal{O}$, since the maps from these categories to the ultimate ones of interest are a stop removal on one side and a pullback on the other. Consider therefore: 
\begin{equation} \label{eq: follow L plus}
    \begin{tikzcd}[column sep=large,row sep=large]
\Fuk\left(\cA( \widetilde B_{[k]})\right)\arrow[r]
&\Fuk\left(  \bC_x^* \times  \bC^*_t \times \cA(\widetilde B_\Sigma )  \right)\\
\Fuk\left( \cA( \widetilde B_{[k]}), T^*S^1 \times \cA(\widetilde B_\Sigma) \right) \arrow[u]  
\arrow[ur, phantom, "\urcorner", very near start]
\arrow[r,"q"]
&\Fuk \left(\bC_x^* \times  \bC^*_t \times \cA(\widetilde B_\Sigma ) , T^*S^1 \times \cA(\widetilde B_\Sigma) \right) \arrow[u]
\end{tikzcd}
\end{equation}

In general, objects of $\mathcal{L} \times_{\mathcal{M}} \mathcal{N}$ are given by tuples
$(\ell, n, \phi: \ell|_{\mathcal{M}} \cong n|_{\mathcal{M}})$, so it would suffice for us to show that the $L_+$ in the bottom corner of \eqref{eq: follow L plus} is given as $(L_+, L_+, 1)$ where the meaning of ``$1$'' is determined using the fixed isomorphisms of the $L_+$ with $\mathcal{O}$.  The fact that $L_+$ restricts to $L_+$ is again Lemma \ref{lem: viterbo respects positive}.  

To identify the isomorphism as ``$1$'', we could again compare after further Viterbo restriction to tori.  Alternatively, note that on the $B$-side, because they are related by an affine blowup, $\operatorname{Pic}(Y) = \operatorname{Pic}(\mathcal{A}(\widetilde B^\vee))$; the latter is trivial by \cite[Thm. 4.13]{Cao-Keller-Qin}.  Thus $\mathcal{O}_Y$ is already determined by the triviality of its restrictions to an open cover, and we do not have to identify the isomorphism. \qed

\subsection{Full rank case}

Let $\widetilde B=\binom BC$ have full rank, with $B$ sink recurrent.
We use $P$, $F$, $S$ and the primitive $\alpha$ of
Proposition~\ref{prop:fullrank-geometry}. Thus
\[
 P=\cA_{\mathrm{prin}}(B)\times(\bC^*)^f,\qquad
 F_i=Z_i\prod_aU_a^{-C_{ai}},\qquad
 F^{-1}(1)=\cA(\widetilde B).
\]
The principal completion~\eqref{eq:principal-completion} is unimodular.
The preceding subsection and K\"unneth therefore give
\begin{equation}\label{eq:fullrank-principal-hms}
 \Ind\Fuk(P)\simeq\QCoh(P^\vee),\qquad
 P^\vee=\cA_{\mathrm{prin}}(-B^{\T})\times(\bC^*)^f,
 \qquad L_+\longmapsto\cO_{P^\vee}.
\end{equation}

We recall the categorical Pontrjagin duality described in \cite[Lemma 8]{Sh21}: for a finitely generated abelian group $\Gamma$ acting on a presentable dg category $\mathcal{C}$, there is a $H = \Hom(\Gamma, \mathbb{G}_m)$ action on $\mathcal{C}^\Gamma$, and a canonical identification
$(\mathcal{C}^\Gamma)^{H} \simeq \mathcal{C}$. 

By Lemmas~\ref{lem:phase-cover} and~\ref{lem:phase-product}, after appropriate compact exact deformation there is a $\Gamma$-cover, respecting Liouville structures, 
$\cA(\widetilde B) \times T^*\bR^m \to P$.  Here,  the domain is an ind Liouville sector.

We claim that
$ \Ind \Fuk ( \cA(\widetilde B) \times T^*\bR^m)^\Gamma \cong \Ind \Fuk(P)$.  One can see this either by passing to microsheaves via \cite{GPS3} (since then the sheaf of categories of microsheaves will be equivariant for the $\Gamma$ action, and limits commute with limits) or alternatively by taking a cover by small domains with transverse boundaries of the  torus $\bR^m/(2\pi S^{-1}\bZ^m)$  to a cover of  $ \cA(\widetilde B)  \times T^*\bR^m$, and then applying sectorial descent.

Taking $H$ invariants and composing with K\"unneth, we obtain: 
$$\Ind \Fuk ( \cA(\widetilde B)) \simeq \Ind \Fuk ( \cA(\widetilde B) \times T^* \bR^m)
\simeq \Ind \Fuk(P)^H \simeq \QCoh(P^\vee)^H \simeq \QCoh([P^\vee / H])
$$

The last isomorphism still requires explanation: we must show that the action of $H$ transported by mirror symmetry to $\QCoh(P^\vee)$ in fact arises from the geometric action of $H$ on $P^\vee$.  
Since $P^\vee$ is smooth and affine, we see by Lemma \ref{lem:gm-action-geometric}, that the transported action necessarily arises from {\em some} geometric action.  
We may identify which action after Viterbo restriction to some seed torus $T^* \T \subset P$ and restriction to a corresponding open in $P^\vee$.  After passing to these tori, we may determine the action instead of the larger group $\Hom(\pi_1(\T), \mathbb{G}_m)$ corresponding to the universal cover.  By factoring, we reduce to the case of the circle, where we are asking to know that the canonical action of $\mathbb{G}_m = \Hom(\mathbb{Z}, \mathbb{G}_m)$ on $\Loc(\R)^{\mathbb{Z}} = B\mathbb{Z}-\operatorname{mod} = \QCoh(\mathbb{G}_m)$ is what you think it is.  We leave this meditation for the interested reader. 

Finally, the matching of $\mathcal{O}$ with $L_+$ and compatibility with the Lagrangian cluster ensemble can be descended similarly. 

The final assertion, for the case where $\cX(\widetilde B^\vee)$ is a scheme, follows similarly by a descent argument. 
Complete the rows of $-\widetilde B^{\T}$ to an integral matrix
\begin{equation}\label{eq:reverse-skew-completion}
 T=\begin{pmatrix}-B^{\T}&-C^{\T}\\C'&T_{\mathrm{ff}}\end{pmatrix},
 \qquad \det T\ne0,
 \qquad \widetilde B'=\binom{-B^{\T}}{C'}.
\end{equation}
The matrix $\widetilde B'$ has full rank, since its columns are the first $m$ columns of $T$.
The groups
\[
 \Gamma:=\Hom(\coker T,\bC^*),\qquad
 K:=\Hom(\coker(-T^{\T}),\bC^*).
\]
are canonically Pontrjagin dual, by the pairing:
\[
 \coker T\times\coker(-T^{\T})\longrightarrow\bC^*,\qquad
 ([a],[b])\longmapsto\exp\bigl(2\pi i\,b^{\T}T^{-1}a\bigr).
\]

The negative transpose is a completion for the dual of $\widetilde B'$:
\[
 -T^{\T}=\begin{pmatrix}B&-C'^{\T}\\C&-T_{\mathrm{ff}}^{\T}\end{pmatrix}.
\]
Applying Proposition~\ref{prop:finite-atlas} twice gives
\begin{equation}\label{eq:reverse-skew-quotients}
 \cX(\widetilde{B}^\vee) \cong[\cA(\widetilde B')/\Gamma],\qquad \cX((\widetilde B')^\vee) \cong[\cA(\widetilde{B})/K].
\end{equation}

So we have
$$\Ind\Fuk(\cX(\widetilde B^\vee)) \simeq \Ind\Fuk(\cA(\widetilde{B'}))^{\Gamma} \simeq \QCoh(\cX((\widetilde B')^\vee))^{\Gamma} \simeq (\QCoh(
\cA(\widetilde{B}))^{K})^{\Gamma} \simeq \QCoh(\cA(\widetilde{B}))$$
where the second identification is the mirror symmetry we have already proven, and the last equality is the categorical Pontrjagin duality.  \qed

\appendix

\section{Skeleton in the acyclic case}
\label{app:acyclic-skeleton}

The  matrix $B$ is \emph{acyclic} if there are no $i_1,\dots,i_r\in M$ with $B_{i_1i_2},\dots,B_{i_{r-1}i_r},B_{i_ri_1}>0$.
In the case where $\widetilde B$ is unimodular and $B$ is acyclic, we can compute the skeleton of \(\cA(\widetilde B)\) directly, and prove Theorem~\ref{thm:main} using microlocal sheaf theory.

The skeleton calculation in Theorem~\ref{thm:acyclic-skeleton} has been formalized in Lean~4 (\url{https://github.com/Hu-Mingyuan/cluster-skeleton}) using Mathlib and~\cite{FS-FLTZ-Lean}. 

We fix a completion $(E,\widetilde D)$, and denote $K=E\widetilde D$. Let \(c_i\) be the characters of \eqref{eq:dual-circles},  and $d_i=\widetilde D_{ii}$. Let \(\Omega\) be the holomorphic symplectic form constructed in Section~\ref{sec:geometry}. 
A subset $I\subset M$ is \emph{independent} if $B_{ij}=0$ for all $i,j\in I$.

In the acyclic case,  $\sA(\widetilde B)$ can be presented by generators and relations.
For $i\in M$ put
\[
 P_i^\pm=\prod_{j\in M\sqcup F}z_j^{[\pm\widetilde B_{ji}]_+},\qquad t_i=P_i^+/P_i^-=z^{c_i}.
\]
By \cite[Corollary~1.21]{BFZ-III},
\begin{equation}\label{eq:acyclic-exchange-model}
 \sA(\widetilde B)=
 \frac{\bC[z_i,z_i'\ (i\in M),\ z_j^{\pm1}\ (j\in F)]}
 {(z_i z_i'-P_i^+-P_i^-)_{i\in M}} .
\end{equation}

\begin{lemma}
[{See also~\cite[Proposition~2.1]{LiuLiuVectorBundles}}]
\label{lem:acyclic-node-charts}
We have the following.
\begin{enumerate}[label=\textup{(\roman*)}]
\item At every point \(p\in\cA(\widetilde B)\), the set $\{i\in M \mid z_i(p)=0\}$ is independent.
\item Let $I$ be independent, and let $U_I=\{z_j\ne0\ (j\in M\setminus I)\}\subset\cA(\widetilde B)$.
Then $U_I\cong\cN^{I}\times(\bC^*)^{m+f-2|I|}$, where the $i$th factor $\cN$ has coordinates
\[
 (x,y,t)=\bigl(z_i,\ z_i'/P_i^-,\ t_i\bigr)\qquad(i\in I).
\]
\item The open sets $U_I$ cover $\cA(\widetilde B)$.
\end{enumerate}
\end{lemma}

\begin{proof}
(i) Fix a point \(p\), let $Z=\{i\in M:z_i(p)=0\}$, and suppose that $B_{Z,Z}\ne0$.
Since \(B\) is acyclic, we can take a maximal path $j_1,\dots,j_r$ in $Z$ with $r\ge2$ and $B_{j_sj_{s+1}}>0$.
Put \(i = j_r\). 
By maximality,  $B_{ji}\ge0$ for all $j\in Z$, and $B_{j_{r-1}i}>0$.
So $P_i^+(p)=0$ since \(z_{j_{r-1}}(p) = 0\) .
Every mutable $j$ in $P_i^-$ has $B_{ji}<0$, so $j\notin Z$; the frozen $z_j$ are units.
So $P_i^-(p)\ne0$. This contradicts the relation $z_iz_i'=P_i^+ + P_i^-$.

(ii) On $U_I$, the $z_j$ with $j\in(M\sqcup F)\setminus I$ are units. 
For $i\in I$, the monomials $P_i^\pm$ involve only these $z_j$.
For $j\in M\setminus I$, we have $z_j'=(P_j^++P_j^-)/z_j$.
So \eqref{eq:acyclic-exchange-model} gives
\[
 \cO(U_I)=\bC\bigl[z_j^{\pm1}\ (j\in M \sqcup F \setminus I)\bigr]\bigl[z_i,y_i\ (i\in I)\bigr]\big/\bigl(z_iy_i-1-t_i\bigr)_{i\in I},
 \qquad y_i=z_i'/P_i^-.
\]
The vectors $c_i$, $i\in I$, are columns of $E\in GL_{m+f}(\bZ)$ and lie in $\bZ^{(M\sqcup F)\setminus I}$.
So they extend to a basis of $\bZ^{(M\sqcup F)\setminus I}$, which implies $\bC[z_j^{\pm1}\ (j\notin I)]=\bC[t_i^{\pm1}\ (i\in I)]\otimes\cO\bigl((\bC^*)^{m+f-2|I|}\bigr)$. This proves (ii).

(iii) For any point \(p\in \cA(\widetilde B)\), take $I = \{ i \in M \mid z_i(p) = 0\}$. Then \( p \in U_I\). 
\end{proof}

\begin{theorem}\label{thm:acyclic-skeleton}
There is a Liouville form $\lambda_{\mathrm{ac}}$ of $\Im\Omega$ on $\cA(\widetilde B)$ such that
\begin{equation}\label{eq:acyclic-skeleton}
 \begin{gathered}
 \Lambda=\Skel(Z_{\lambda_{\mathrm{ac}}})=\bigcup_{I\subset M\text{ independent}}\Lambda_I,\\
 \Lambda_I=\bigl\{z_i'=0,\ |z_i|\le1\ (i\in I),\quad |z_j|=1\ (j\in(M\sqcup F)\setminus I)\bigr\}.
 \end{gathered}
\end{equation}

\end{theorem}

\begin{proof}

Use the constants $L>a>0$ and the cutoffs $h,\rho,\beta$ of Section~\ref{sec:node-model}. Let \(\lambda_{\mathrm{tor}} = u^\T K \dd \theta \) as in~\eqref{eq:forms}. 
Fix $\delta_0>0$ with $h=0$ on $(-\delta_0,\infty)$.
Choose a smooth nonincreasing $\vartheta\colon\bR\to[0,1]$ with $\vartheta=1$ on $(-\infty,-1]$ and $\vartheta=0$ on $[0,\infty)$.
On the initial torus $(\bC^*)^{m+f}$ write
\[
 u_j=\log|z_j|,\qquad \theta_j=\arg z_j,\qquad
 v_i=\log|t_i|=c_i(u),\qquad \eta_i=\arg t_i=c_i(\theta),
\]
and set
\[
 \chi_i \coloneqq \prod_{j\in M:\,B_{ji}\ne0}\vartheta(u_i-u_j),\qquad
 H \coloneqq \sum_{i\in M}d_i\,h(u_i)\,\chi_i\,\rho(v_i)\,\beta(\eta_i),\qquad
 \lambda_{\mathrm{ac}} \coloneqq\ltor+\dd H .
\]
An empty product is $1$. Note that 
if $B_{ij}\ne0$, then $\chi_i\chi_j=0$.

Let $I\subset M$ be independent.
For $i,i'\in I$, $K_{ii'}=d_{i'}B_{ii'}=0$.
For $i\in M$,
\[
 \sum_{j}K_{ji}\,u_j=d_i\sum_jE_{ji}u_j=d_iv_i,\qquad
 \sum_lK_{il}\,\dd\theta_l=-\sum_lK_{li}\,\dd\theta_l=-d_i\,\dd\eta_i .
\]
The terms of $\ltor$ with $l\in I$ give $\sum_{i\in I}d_iv_i\,\dd\theta_i$.
The terms with $j\in I$ and $l\notin I$ give $-\sum_{i\in I}d_iu_i\,\dd\eta_i$, since $K_{il}=0$ for $i,l\in I$.
Hence
\begin{equation}\label{eq:acyclic-ltor-split}
 \ltor=\sum_{i\in I}d_i\bigl(v_i\,\dd\theta_i-u_i\,\dd\eta_i\bigr)+\sum_{j,l\in(M\sqcup F)\setminus I}K_{jl}\,u_j\,\dd\theta_l .
\end{equation}

We first show that $\lambda_{\mathrm{ac}}$ is globally defined.
Let $q\in\cA(\widetilde B)\setminus(\bC^*)^{m+f}$, and let $I=\{i\in M\mid z_i(q)=0\}$.
Then $I\ne\emptyset$, $I$ is independent by Lemma~\ref{lem:acyclic-node-charts}(i), and $q\in U_I$.
Let $i\in I$.
The monomials $P_i^\pm$ only involve $z_j$ with $j\notin I$, so they are units on $U_I$.
Hence $t_i=P_i^+/P_i^-$ is a unit on $U_I$, and $0=z_i(q)z_i'(q)=P_i^-(q)\bigl(1+t_i(q)\bigr)$ gives $t_i(q)=-1$.
So, as $p\to q$ with $p\in(\bC^*)^{m+f}$,
\[
 u_i(p)\to-\infty,\qquad v_i(p)\to0,\qquad \eta_i(p)\to\pi,\qquad u_j(p)\to\log|z_j(q)|\in\bR\quad\bigl(j\in(M\sqcup F)\setminus I\bigr).
\]
Choose a neighborhood $W\subset U_I$ of $q$ such that, for $i\in I$ and $j\in M$ with $B_{ji}\ne0$,
\[
 \rho(v_i)=1,\quad \beta(\eta_i)=\eta_i-\pi\quad\text{on }W,\qquad
 u_i\le-L,\quad u_i\le u_j-1\quad\text{on }W\cap(\bC^*)^{m+f}.
\]
Then, on $W\cap(\bC^*)^{m+f}$, $h(u_i)=u_i+a$ and $\chi_i=1$.
Also $\chi_j=0$ there, because $\chi_j$ has the factor $\vartheta(u_j-u_i)$ and $u_j-u_i\ge1$.
So, on $W\cap(\bC^*)^{m+f}$,
\begin{equation}\label{eq:acyclic-H-near}
 H=\sum_{i\in I}d_i(u_i+a)(\eta_i-\pi)+H',
\end{equation}
where $H'$ is the sum of the summands of $H$ with index $k\in M\setminus I$ such that $B_{ik}=0$ for all $i\in I$.
For such $k$, the functions $u_k$, $\chi_k$, $v_k=c_k(u)$ and $\eta_k=c_k(\theta)$ only involve $z_j$ with $j\notin I$, by the support of $c_k$.
These $z_j$ are units on $U_I$.
So $H'$ does not depend on $u_i$ and $\theta_i$, $i\in I$.

By \eqref{eq:acyclic-ltor-split} and \eqref{eq:acyclic-H-near}, on $W\cap(\bC^*)^{m+f}$,
\[
 \lambda_{\mathrm{ac}}=\sum_{j,l\notin I}K_{jl}\,u_j\,\dd\theta_l+\dd H'+\sum_{i\in I}d_i\bigl(v_i\,\dd\theta_i+(\eta_i-\pi)\,\dd u_i+a\,\dd\eta_i\bigr).
\]
Here we used $v_i\,\dd\theta_i-u_i\,\dd\eta_i+\dd\bigl((u_i+a)(\eta_i-\pi)\bigr)=v_i\,\dd\theta_i+(\eta_i-\pi)\,\dd u_i+a\,\dd\eta_i$.
The first two terms are smooth on $W$.
On $W$, $1-z_iy_i=-t_i$ and $\Log(-t_i)=v_i+\mathrm i(\eta_i-\pi)$.
So, as in the proof of Lemma~\ref{lem:node},
\[
 v_i\,\dd\theta_i+(\eta_i-\pi)\,\dd u_i+a\,\dd\eta_i
 =\Im\Bigl(\frac{\Log(1-z_iy_i)}{z_i}\,\dd z_i\Bigr)+a\,\dd\arg(1-z_iy_i).
\]
The right side is smooth on $W$, because $\Log(1-z_iy_i)/z_i$ is holomorphic on $W$, with value $-y_i$ on $\{z_i=0\}$.
Hence $\lambda_{\mathrm{ac}}$ extends smoothly over $W$,  hence also extends to $\cA(\widetilde B)$.

Now we consider the Liouville field $Z=Z_{\lambda_{\mathrm{ac}}}$, given by $\iota_Z\Im\Omega=\lambda_{\mathrm{ac}}$.
A dot denotes the derivative along $Z$.
On the initial torus, $\Im\Omega=\sum_{j,l}K_{jl}\,\dd u_j\wedge\dd\theta_l$.
Comparing the coefficients of $\dd\theta$ and of $\dd u$ in $\iota_Z\Im\Omega=\ltor+\dd H$ gives
\[
 K^{\T}\dot u=K^{\T}u+\partial_\theta H,\qquad -K\dot\theta=\partial_uH .
\]
The function $H$ depends on $\theta$ only through $\eta_i=c_i(\theta)$, and $d_ic_i=Ke_i$.
So
\[
 \partial_\theta H=\sum_{i\in M}d_i\,h(u_i)\chi_i\rho(v_i)\beta'(\eta_i)\,c_i
 =K\sum_{i\in M}h(u_i)\chi_i\rho(v_i)\beta'(\eta_i)\,e_i .
\]
Since $K^{\T}=-K$, this gives
\begin{equation}\label{eq:acyclic-flow}
 \begin{aligned}
 \dot u_j&=u_j &&(j\in F),\\
 \dot u_i&=u_i-h(u_i)\,\chi_i\,\rho(v_i)\,\beta'(\eta_i) &&(i\in M),\\
 \dot\theta&=-K^{-1}\partial_uH .
 \end{aligned}
\end{equation}
For $i\in M$, $c_i(K^{-1}x)=-x_i/d_i$, because $d_ic_i^{\T}=e_i^{\T}K^{\T}=-e_i^{\T}K$.
Also $c_i(e_k)=E_{ki}=B_{ki}$ for $k\in M$.
So applying $c_i$ to \eqref{eq:acyclic-flow} gives
\begin{equation}\label{eq:acyclic-angle-flow}
 \dot\eta_i=\frac1{d_i}\,\partial_{u_i}H,\qquad
 \dot v_i=v_i-\sum_{k\in M}B_{ki}\,h(u_k)\,\chi_k\,\rho(v_k)\,\beta'(\eta_k)\qquad(i\in M).
\end{equation}

We now bound $\dot u$.
Since $h\le0$, $0\le\chi_i\rho\le1$ and $\beta'\le1$, we have $h\chi_i\rho\beta'\ge h\chi_i\rho\ge h$.
With $h=0$ on $(-\delta_0,\infty)$, \eqref{eq:acyclic-flow} gives, for $i\in M$,
\begin{equation}\label{eq:acyclic-radial}
 \dot u_i\le u_i-h(u_i)\quad\text{on }\{u_i<0\},\qquad \dot u_i=u_i\quad\text{on }\{u_i>-\delta_0\}.
\end{equation}
Recall that $u-h(u)<0$ for $u<0$, and $u-h(u)=-a$ for $u\le-L$.
Also $|h(u)|\le|u|+L$, since $h(u)=u+a$ for $u\le-L$, and $a-L\le h\le0$ on $[-L,\infty)$.
So \eqref{eq:acyclic-flow} gives $|\dot u_j|\le C(1+|u_j|)$ for all $j\in M\sqcup F$, with $C=(1+\max|\beta'|)(1+L)$.
By density, \eqref{eq:acyclic-radial} holds on $\{z_i\ne0\}$, and the bound holds on $\{z_j\ne0\}$.

Next we compute $Z$ near $\{z_i=0\}$.
Let $q$, $I$ and $W$ be as above, and let $i\in I$.
On $W\cap(\bC^*)^{m+f}$, $h(u_i)=u_i+a$, $\chi_i=\rho(v_i)=\beta'(\eta_i)=1$, and $\chi_k=0$ for $B_{ki}\ne0$.
Also $\partial_{u_i}H=d_i(\eta_i-\pi)$ there, by \eqref{eq:acyclic-H-near}.
So \eqref{eq:acyclic-flow} and \eqref{eq:acyclic-angle-flow} give
\[
 \dot u_i=-a,\qquad \dot v_i=v_i,\qquad \dot\eta_i=\eta_i-\pi\qquad\text{on }W\cap(\bC^*)^{m+f}.
\]
Since $\Log(-t_i)=v_i+\mathrm i(\eta_i-\pi)$, by density
\begin{equation}\label{eq:acyclic-divisor-flow}
 Z|z_i|^2=-2a|z_i|^2,\qquad \dot t_i=t_i\Log(-t_i)\qquad\text{on }W.
\end{equation}
Every point of $\{z_i=0\}$ is such a $q$, so these sets $W$ cover $\{z_i=0\}$.
Along an orbit near $\{z_i=0\}$, $|z_i|^2$ solves $\frac{\dd}{\dd\tau}|z_i|^2=-2a|z_i|^2$.
By uniqueness of solutions, $\{z_i=0\}$ and its complement are invariant under the flow of $Z$.

On $\{z_i=0\}$, the $z_j$ with $B_{ji}\ne0$ are nonzero, by Lemma~\ref{lem:acyclic-node-charts}(i), so $P_i^-\ne0$.
So $y_i=z_i'/P_i^-$ is defined near $\{z_i=0\}$, and $y_i=(1+t_i)/z_i$ where $z_i\ne0$.
Put $\phi(t)=t\Log(-t)/(1+t)$.
It is holomorphic near $t=-1$, with $\phi(-1)=1$.
On $W\cap\{z_iy_i\ne0\}$, \eqref{eq:acyclic-divisor-flow} gives
\[
 Z\log|y_i|=\Re\frac{\dot t_i}{1+t_i}-Z\log|z_i|=\Re\phi(t_i)+a .
\]
By density, $Z|y_i|^2=2\bigl(\Re\phi(t_i)+a\bigr)|y_i|^2$ on $W$.
On $W\cap\{z_i=0\}$, $t_i=-1$, so
\begin{equation}\label{eq:acyclic-partner-flow}
 Z|y_i|^2=2(1+a)|y_i|^2\qquad\text{on }\{z_i=0\}.
\end{equation}
Hence $|y_i(\phi_Z^\tau(p))|=e^{(1+a)\tau}|y_i(p)|$ for $p\in\{z_i=0\}$ and all $\tau$.

Next we show that $Z$ is complete.
Let $\gamma\colon(\tau_-,\tau_+)\to\cA(\widetilde B)$ be a maximal integral curve of $Z$, and let $I_0=\{i\in M\mid z_i(\gamma(0))=0\}$.
By invariance, $z_i=0$ on $\gamma$ for $i\in I_0$, and $z_j\ne0$ on $\gamma$ for $j\in(M\sqcup F)\setminus I_0$.
Let $T>0$.
For $\tau\in(\tau_-,\tau_+)$ with $|\tau|\le T$, the bound $|\dot u_j|\le C(1+|u_j|)$ and \eqref{eq:acyclic-partner-flow} give
\[
 |u_j(\gamma(\tau))|\le\bigl(1+|u_j(\gamma(0))|\bigr)e^{CT}\quad(j\notin I_0),\qquad
 |y_i(\gamma(\tau))|\le e^{(1+a)T}|y_i(\gamma(0))|\quad(i\in I_0).
\]
So the coordinates $z_j^{\pm1}$ $(j\notin I_0)$, $z_j'=(P_j^++P_j^-)/z_j$ $(j\in M\setminus I_0)$ and $z_i'=P_i^-y_i$ $(i\in I_0)$ are bounded on these $\gamma(\tau)$.
By \eqref{eq:acyclic-exchange-model}, $\cA(\widetilde B)$ is closed in $\bC^{2m}\times(\bC^*)^f$, so these $\gamma(\tau)$ lie in a compact set.
Hence $\tau_\pm=\pm\infty$, since a maximal integral curve with $\tau_+<\infty$ or $\tau_->-\infty$ leaves every compact set.

We now prove \eqref{eq:acyclic-skeleton}.

(i) \(\Skel(Z) \subset \cup_I \Lambda_I\).
Let $p\in\Skel(Z)$, so that the closure of $\gamma(\tau)=\phi_Z^\tau(p)$, $\tau\ge0$, is compact.
Put
\[
 I_0=\{i\in M\mid z_i(p)=0\},\qquad I=I_0\cup\{i\in M\mid z_i(p)\ne0,\ u_i(p)<0\}.
\]
We show that $p\in\Lambda_I$.

We first determine $u$ along $\gamma$.
For $j\in F$, $u_j(\gamma(\tau))=e^\tau u_j(p)$ is bounded, so $u_j(p)=0$.
Let $i\in M\setminus I_0$.
If $u_i(p)>0$, then $u_i(\gamma(\tau))=e^\tau u_i(p)$ by \eqref{eq:acyclic-radial}, which is unbounded.
If $u_i(p)=0$, then $u_i=0$ on $\gamma$, by \eqref{eq:acyclic-radial} and uniqueness for $\dot u_i=u_i$.
If $u_i(p)<0$, then $u_i$ decreases along $\gamma$, and \eqref{eq:acyclic-radial} gives
\[
 \dot u_i\le\sup_{u\le u_i(p)}\bigl(u-h(u)\bigr)<0\qquad\text{on }\gamma .
\]
The supremum is negative because $u-h(u)<0$ for $u<0$ and $u-h(u)=-a$ for $u\le-L$.
So $u_i(\gamma(\tau))\to-\infty$.
In summary,
\[
 u_j=0\ \text{on }\gamma\ \ \bigl(j\in(M\sqcup F)\setminus I\bigr),\qquad
 u_i(\gamma(\tau))\to-\infty\ \ (i\in I\setminus I_0),\qquad
 z_i=0\ \text{on }\gamma\ \ (i\in I_0).
\]
Let $q$ be a limit point of $\gamma$.
Then $\{i\in M\mid z_i(q)=0\}=I$, so $I$ is independent by Lemma~\ref{lem:acyclic-node-charts}(i).
Also $|z_j|=1$ on $\gamma$ and at $q$ for $j\in M\setminus I$.
So the closure of $\gamma$ is a compact subset of $U_I$.
On $U_I$, $z_i'=P_i^-y_i$ with $P_i^-$ a unit, for $i\in I$.
Since $|z_i(p)|\le1$ for $i\in I$ and $|z_j(p)|=1$ for $j\notin I$, it remains to show that $y_i(p)=0$ for $i\in I$.

Let $i\in I_0$.
By \eqref{eq:acyclic-partner-flow}, $|y_i(\gamma(\tau))|=e^{(1+a)\tau}|y_i(p)|$, which is bounded only if $y_i(p)=0$.

Let $i\in I\setminus I_0$.
We show that $t_i=-1$ on $\gamma$.
The character $c_i$ is supported on $F\sqcup\{j\in M\mid B_{ji}\ne0\}$, which is disjoint from $I$.
So $v_i=c_i(u)=0$ on $\gamma$.
Consider the open set
\[
 O=\bigl\{z_i\ne0,\ |z_j|>e^{-\delta_0}\ (j\in M,\ B_{ji}\ne0)\bigr\}\supset\gamma .
\]
Besides the $i$th summand of $H$, only the summands with index $j$, $B_{ji}\ne0$, depend on $u_i$, through $\chi_j$ and $v_j=c_j(u)$.
These vanish on $O\cap(\bC^*)^{m+f}$, because $h(u_j)=0$ there.
Also $\partial_{u_i}v_i=E_{ii}=0$.
So \eqref{eq:acyclic-angle-flow} gives, on $O\cap(\bC^*)^{m+f}$, and then on $O$ by density,
\[
 \dot\eta_i=\bigl(h'(u_i)\,\chi_i+h(u_i)\,\partial_{u_i}\chi_i\bigr)\,\rho(v_i)\,\beta(\eta_i).
\]
The factor in front of $\beta(\eta_i)$ is $\ge0$, because $h'\ge0$, $h\le0$, $\chi_i\ge0$, $\rho\ge0$, and $\partial_{u_i}\chi_i\le0$ as $\vartheta'\le0$.
Since $(\eta-\pi)\beta(\eta)\ge0$ and $\supp\beta\subset(\pi-\varepsilon,\pi+\varepsilon)$, the distance $|\eta_i-\pi|$ in $\bR/2\pi\bZ$ is nondecreasing along $\gamma$.
On the other hand, $z_i(\gamma(\tau))\to0$, and $y_i$ is bounded on $\gamma$.
So $1+t_i=z_iy_i\to0$, and $|\eta_i-\pi|\to0$.
Hence $\eta_i=\pi$ on $\gamma$.
With $v_i=0$, this gives $t_i=-1$ on $\gamma$, so $y_i=(1+t_i)/z_i=0$ on $\gamma$.
Thus $p\in\Lambda_I$, and $\Skel(Z)\subset\bigcup_I\Lambda_I$.

(ii) \( \Lambda_I \subset \Skel(Z)\). 

Let $I$ be independent, and let $p\in\Lambda_I$.
The set $\Lambda_I$ is closed, and on it the coordinates $z_j^{\pm1}$ $(j\notin I)$, $z_i$, $z_i'=0$ $(i\in I)$ and $z_j'=(P_j^++P_j^-)/z_j$ $(j\in M\setminus I)$ are bounded.
So $\Lambda_I$ is compact, and it suffices to show $\phi_Z^\tau(p)\in\Lambda_I$ for $\tau\ge0$.
Let $\gamma(\tau)=\phi_Z^\tau(p)$.
As above, $u_j=0$ on $\gamma$ for $j\notin I$.
For $i\in I$, $z_i=0$ on $\gamma$ if $z_i(p)=0$, and $u_i\le0$ on $\gamma$ if $z_i(p)\ne0$.
It remains to show $y_i=0$ on $\gamma$ for $i\in I$.
If $z_i(p)=0$, this follows from \eqref{eq:acyclic-partner-flow} and $y_i(p)=0$.
If $z_i(p)\ne0$, then $v_i=0$ on $\gamma$, and $\eta_i$ solves the equation for $\dot\eta_i$ above, with $\eta_i(p)=\pi$.
Since $\beta(\pi)=0$, uniqueness gives $\eta_i=\pi$ on $\gamma$.
So $t_i=-1$ and $y_i=0$ on $\gamma$.
This proves \eqref{eq:acyclic-skeleton}.
\end{proof}

With the skeleton in hand, we prove homological mirror symmetry by microlocal sheaf theory.

\begin{theorem}\label{thm:acyclic-skeletal-hms}
\(
 \mu\Sh(\Lambda)^c \simeq\Coh\bigl(\cX(\widetilde B^\vee)\bigr).
\)
\end{theorem}

Here $\mu\Sh(\Lambda)$ is the  (large) category of microsheaves on $\Lambda$ \cite{NS}, and $\mu\Sh(\Lambda)^c$ is its subcategory of compact objects. 
By \cite{NZ, GPS3}, Theorem~\ref{thm:acyclic-skeletal-hms} gives Theorem~\ref{thm:main} for acyclic $B$.

We first set up the A-side.
Fix $0<\varepsilon<\delta\ll 1$.
For independent $I\subset J \subset M$ with \(J\) independent, put
\[
 V_{I, J} \coloneqq \bigl\{p\in\Lambda \mid |z_j(p)|<\delta+\varepsilon\ (j\in J),\ \ |z_i(p)|>\delta-\varepsilon\ (i\in M\setminus I)\bigr\}.
\]
Denote \(V_I \coloneqq V_{I, I}\). 
We have \( V_{\emptyset} \simeq \Lambda \cap (\bC^*)^{m+f}\).

For each $i\in J\setminus I$, there are open embeddings
\[
 V_{I\cup\{i\},J}\xleftarrow{\ h_i\ }V_{I,J}\xrightarrow{\ f_i\ }V_{I,J\setminus\{i\}},
\]
which induce corestriction functors (i.e., left adjoint to the restrictions)
\[
\mu\Sh(V_{I\cup\{i\},J})
\xleftarrow{(h_i)_!}
\mu\Sh(V_{I,J})
\xrightarrow{(f_i)_!}
\mu\Sh(V_{I,J\setminus\{i\}}) 
\]
Consider the poset
\[
\mathcal P=\{(I,J) \mid I\subset J\subset M,\ J\text{ independent}\}. 
\]
There is an arrow \( (I, J) \to (I', J')\) if and only if \(I \subset I' \subset J'\subset J\). 
Note that 
morphisms of \(\cP\) are generated by \( (I, J) \to (I \cup \{i\}, J)\) and \( (I, J) \to (I, J\setminus \{i\} )\), for all \(i \in J\setminus I\).  

\begin{lemma}\label{lem:acyclic-microlocal}
    The sets $V_{I}$ cover $\Lambda$, and
\[
 \mu\Sh(\Lambda)\simeq\operatorname*{colim}^{\mathrm{Pr}^\mathrm{L}}_{\cP}\bigl(\mu\Sh(V_{I,J}),\ (f_i)_!,\ (h_i)_!\bigr).
\]
\end{lemma}

\begin{proof}
Let $p\in\Lambda_{I}$, and put 
$I_p=\{i\in M\mid |z_i(p)|<\delta+\varepsilon\}$. If \(i\in I_p\), then \(|z_i(p)| \le \delta + \epsilon < 1\). Hence \(i\in I\), and \(I_p \subset I\). In particular, \(I_p\) is independent. For \(j \notin I_p\), we have \(|z_j(p)| \ge \delta + \epsilon\), so \(p\in V_{I_p}\). This shows that \(\{V_I\}\) cover \(\Lambda\).

For independent $I_1,\dots,I_k$, put $I=\bigcap_\alpha I_\alpha$ and $J=\bigcup_\alpha I_\alpha$.
Then $V_{I_1}\cap\dots\cap V_{I_k}$ is empty unless $J$ is independent, and in that case it is $V_{I,J}$.
Since $\mu\Sh$ is a cosheaf~\cite{NS}, $\mu\Sh(\Lambda)$ is the colimit of the $\mu\Sh(V_{I,J})$ along the corestrictions.
\end{proof}

Now we consider the B-side.
Let $z_j$, $j\in M\sqcup F$, be the coordinates of $(\bC^*)^{M\sqcup F}$, and for $i\in M$ put
\[
  t_i=z^{Ee_i}=\prod_{j\in M\sqcup F}z_j^{E_{ji}} .
\]
Let
\[
 Y_\Sigma=\Spec\bC\bigl[z_i\ (i\in M),\ z_j^{\pm1}\ (j\in F)\bigr]\setminus\bigcup_{B_{ij}\ne0} \{ z_i = 0, z_j = 0\}
\]
It is the toric variety of the fan $\Sigma=\{\sigma_S\mid S\subset M\text{ independent}\}$, $\sigma_S=\sum_{i\in S}\bR_{\ge0}e_i$, with affine charts $\{z_l\ne0\ (l\in M\setminus S)\}$.

For each pair $I\subset J$ with \(J\) independent, take
\[
 Y_{I,J} \coloneqq  \{z_i=0\ (i\in J),\ \ t_i=-1\ (i\in I)\}\subset Y_\Sigma .
\]
On $Y_{I,J}$, the $z_l$ with $B_{il}\ne0$ for some $i\in I$ are units.
The exponents of the $t_i$, $i\in I$, are columns of $E\in GL_{m+f}(\bZ)$, so the $t_i$ are part of a coordinate system on these units.
Hence $Y_{I,J}$ is smooth.
For each $i\in J\setminus I$, there are closed embeddings
\[
 Y_{I\cup\{i\},J}\xrightarrow{\ \kappa_i\ }Y_{I,J}\xrightarrow{\ \iota_i\ }Y_{I,J\setminus\{i\}},
\]
cut out by $t_i=-1$ and $z_i=0$. 
Consider the functors
\[
\QCoh( Y_{I \cup \{i\}, J} ) \xleftarrow{\kappa_i^*} \QCoh(Y_{I, J}) \xrightarrow{\iota_{i\times}} \QCoh( Y_{I, J\setminus\{i\} })
\]
Here \( \iota_{i \times} \coloneqq  \iota_{i*}[-1] \) 
is the left adjoint to \(\iota_i^*\)  (see~\cite[Proposition 4.5]{FSZ}).

\begin{lemma}\label{lem:acyclic-blowup}
There is an equivalence
\begin{equation}\label{eq:acyclic-microlocal-colimit}
 \QCoh\bigl(\cA(\widetilde B )\bigr)
 \simeq\operatorname*{colim}_{\cP}^{\operatorname{Pr}^L}
 \bigl(\QCoh(Y_{I,J}), \iota_{i\times}, \kappa_i^*\bigr).
\end{equation}
\end{lemma}

\begin{proof}

For an independent $S\subset M$, put $O_S=\{z_l\ne0\ (l\in M\setminus S)\}\subset Y_\Sigma$. Then \[O_S \cong \prod_{i\in S} \bC_{z_i} \times (\bC^*)^{m+f-|S|} .\]

Since $E_{ji}=0$ for $i,j\in S$, we have \(c_i = E e_i \in \bZ^{M \sqcup F \setminus S}\). Since $E$ is unimodular we can complete \( \{c_i\}_{i\in S}\) to a basis of \( \bZ^{M \sqcup F \setminus S} \). Hence (recall \(t_i = z^{c_i}\)) 
\[
 O_S\cong\prod_{i\in S}\bigl(\bC_{z_i}\times\bC^*_{t_i}\bigr)\times(\bC^*)^{m+f-2|S|}.
\]
If \(J \not \subset S\), then  $Y_{I,J}\cap O_S=\emptyset$. For \( I \subset J \subset S\),
\[
Y_{I, J} \cap O_S \cong \big\{ z_i = 0 \ (i \in J), \quad t_i = -1\  (i\in I) \big\} \cong \prod_{i\in S\setminus J} \bC_{z_i} \times \prod_{i\in S\setminus I} \bC_{t_i}^* \times (\bC^*)^{m+f-2|S|}.
\]
Therefore
\begin{align*}
    &\colim_{ \cP} \QCoh( Y_{I, J} \cap O_S) \\
    & \simeq \bigotimes_{i \in S} \left(  \colim (\QCoh( \bC_{z_i} \times \bC_{t_i}^* ) \xleftarrow{z_i = 0} \QCoh(\bC_{t_i}^*) \xrightarrow{t_i= -1} \QCoh(\mathrm{pt}) ) \right) \otimes \QCoh((\bC^*)^{m+f-2|S|} )\\
    &\simeq \bigotimes_{i\in S} \QCoh\left( \Bl_{(0, -1)} (\bC_{z_i} \times \bC_{t_i}^*) \setminus \widetilde{ \{ z_i = 0\} } \right) \otimes \QCoh\left( (\bC^*)^{m+f - 2|S|} \right) \\
     &\simeq \QCoh\left(\prod_{i\in S}\{z_iy_i=1+t_i\}\times(\bC^*)^{m+f-2|S|}\right)
    \simeq\QCoh(U_S).
\end{align*}
The second equivalence is \cite[Theorem 4.13]{GLe}, and the last is Lemma~\ref{lem:acyclic-node-charts}(ii).

By \cite[Corollary~5.5.3.4 and Theorem~5.5.3.18]{Lurie-HTT},
\[
\colim_{\cP} \bigl(\QCoh(Y_{I,J}), \iota_{i \times}, \kappa_{i}^*\bigr)
    \simeq \lim_{\cP^{\mathrm{op}}}
    \left(
    \QCoh(Y_{I,J}), \  \iota_i^*, \kappa_{i*}
    \right)
\]
Note that \(\{O_S\}\) (\(S\) independent) forms an open cover of \( Y_\Sigma\). 
By Zariski descent,
\begin{align*}
    \lim_{\cP^{\mathrm{op}}}
    \left(
    \QCoh(Y_{I,J}), \  \iota_i^*, \kappa_{i*}
    \right) 
    &\simeq
    \lim_{\cP^{\mathrm{op}}} 
    \lim_{S}
    \left( \QCoh(O_S \cap Y_{I, J} ) \right) \\
    &\simeq\lim_S \lim_{\cP^{\mathrm{op}}} 
    \left( \QCoh(O_S \cap Y_{I, J} ) \right) \\
    &\simeq \lim_S \QCoh(U_S) \\
    &\simeq \QCoh\bigl(\cA(\widetilde B)\bigr).
\end{align*}
\end{proof}

Now we can prove Theorem~\ref{thm:acyclic-skeletal-hms} by matching the two colimit diagrams. 
\begin{proof}[Proof of Theorem~\ref{thm:acyclic-skeletal-hms}]
We use Lemma~\ref{lem:acyclic-blowup} for $\widetilde B^\vee$ with completion $E^\vee=-E^{\T}$, retaining the notation $Y_\Sigma$ and $Y_{I,J}$.
Write $\c z_j$ for the B-side coordinates, so
\[
 \c t_i=(\c z)^{E^\vee e_i}=\prod_{j\in M\sqcup F}(\c z_j)^{-E_{ij}}.
\]

We first identify the toric part.
On the initial torus, shift the base coordinates of~\eqref{eq:identify_with_torus} by $\pi\mathbf 1_M$:
\[
 \psi=-E^{\T}\theta+\pi\mathbf1_M,\qquad p=\widetilde Du,\qquad \ltor=p^{\T}\dd\psi .
\]
Then $\psi_i=\pi-\eta_i$ and $p_i=d_iu_i$ for $i\in M$, and
\[
 \Lambda_I\cap(\bC^*)^{m+f}=\bigl\{\psi_i=0,\ p_i\le0\ (i\in I),\qquad p_j=0\ (j\notin I)\bigr\}.
\]
So $\Lambda\cap(\bC^*)^{m+f}$ is the Lagrangian $\bigcup_S\sigma_S^\perp\times(-\sigma_S)$ of the fan $\Sigma$. By the coherent--constructible--correspondence(CCC)~\cite{FLTZ,Kuwagaki},
\[
\mu\Sh\bigl(\Lambda\cap(\bC^*)^{m+f}\bigr)\simeq\QCoh(Y_\Sigma).
\]
Here $\c z_j$ is the holonomy of a local system along the $\psi_j$-circle.
For $i\in I$, the circle $\{|z_i|=1\}\subset\Lambda_I$, with the other $z_j$ fixed, has class $-E^{\T}e_i=E^\vee e_i$, so the holonomy along it is $\c t_i$.
The Lie group spin structure on this circle does not extend over the disk $\{|z_i|\le1\}$, as after~\eqref{eq:node-A-diagram}.
So a local system on the disk has $\c t_i=-1$.

Next we identify the vertices and the edges. Fix $I\subseteq J$.
The characters $\c t_j$, $j\in J$, involve only unit coordinates on $Y_{\emptyset,J}$.
Since $E^\vee$ is unimodular, they extend to a basis of the lattice of these unit coordinates.
Choose such a basis, and put $Z_J=\{\c t_j=1\ (j\in J)\}\subset Y_{\emptyset,J}$.
Then
\[
 \begin{aligned}
 Y_{\emptyset,J}&\cong\prod_{j\in J}\bC^*_{\c t_j}\times Z_J,\\
 Y_{I,J}&\cong\prod_{i\in I}\{\c t_i=-1\}
 \times\prod_{j\in J\setminus I}\bC^*_{\c t_j}\times Z_J.
 \end{aligned}
\]
Let $D_i=\{|z_i|<\delta+\varepsilon\}$ and $\mathrm{Ann}_i=\{\delta-\varepsilon<|z_i|<\delta+\varepsilon\}$.
Then
\[
 V_{I,J}\cong
 \prod_{i\in I}D_i\times\prod_{j\in J\setminus I}\mathrm{Ann}_j\times\Lambda_{Z_J},
\]
where $\Lambda_{Z_J}$ is the FLTZ skeleton of $Z_J$.
Thus K\"unneth and CCC for $Z_J$ give
\[
 \begin{aligned}
 \mu\Sh(V_{I,J})
 &\simeq
 \left(\bigotimes_{i\in I}\QCoh(\{\c t_i=-1\})\right)
 \otimes\left(\bigotimes_{j\in J\setminus I}\QCoh(\bC^*_{\c t_j})\right)
 \otimes\QCoh(Z_J)\\
 &\simeq\QCoh(Y_{I,J}).
 \end{aligned}
\]
Under these equivalences, the restriction along $f_i$ is $ \iota_i^*$, the pullback to $\{\c z_i=0\}$.
The restriction along $h_i$ restricts local systems from the disk $\{|z_i|<\delta+\varepsilon\}$ to the annulus, so it is $ \kappa_{i*}$.
Taking left adjoints, $(f_i)_!=\iota_{i*}[-1]=\iota_{i\times}$, since $\{\c z_i=0\}$ has trivial normal bundle, and $(h_i)_!= \kappa_i^*$.
These identifications are compatible with composition, and we omit the routine check.
So the colimit diagrams of Lemmas~\ref{lem:acyclic-microlocal} and~\ref{lem:acyclic-blowup} agree, and $\mu\Sh(\Lambda)\simeq\QCoh(\cA(\widetilde B^\vee))$. Proposition~\ref{prop:finite-atlas} with $T=E^\vee$ gives $\cA(\widetilde B^\vee)\cong\cX(\widetilde B^\vee)$. Taking compact objects proves the statement.
\end{proof}

\begin{remark}
    If we did not allow mutations in the definition of sink-recurrence, the result would be the definition of acyclic seeds.  From this point of view, the reason we can calculate the skeleton in the acyclic case but not the sink-recurrent case is that we do not know how to calculate the result of mutations on skeleta (except in real dimension 4, for which see \cite{STW}).  
\end{remark}

\section{Comparison with the affine Liouville structure}\label{app:affine}

In \cite[Theorem 1.1]{HS26}, we showed that arbitrary symplectomorphisms of Liouville manifolds -- no control at infinity -- give equivalences of Fukaya categories.  
In this appendix,
we use a criterion of \cite{HS26} to construct such symplectomorphisms relating the Liouville structures we have used throughout this article to the Liouville structure given by a primitive
for the K\"ahler form, hence showing that they determine the same Fukaya category.

The criterion from \cite{HS26} concerns the existence of an embedding in which the Poisson brackets for $\Omega$ are quadratic.  We obtain such an embedding here using the theta functions of \cite{GHKK}.

\begin{lemma}\label{lem:theta-filtration}
Suppose that $\cA(\widetilde B)$ has large cluster complex \cite[Definition~8.27]{GHKK} and that $\sA(\widetilde B)=\sU(\widetilde B)$.
Let a torus $T$ act on $\cA(\widetilde B)$ such that all cluster variables are weight vectors. Then there is a filtration $(F_n)_{n\ge0}$ of $\sA(\widetilde B)$ by $T$-invariant subspaces such that
\[
 F_0=\bC,\qquad \dim F_n<\infty,\qquad \bigcup_nF_n=\sA(\widetilde B),\qquad F_nF_l\subset F_{n+l},
\]
and the Rees algebra $\bigoplus_nF_nt^n$ is finitely generated.
\end{lemma}

\begin{proof}
We filter the principal-coefficient variety by theta functions, and then set the coefficients to $1$.

Let $V$ be the variety $\cA_{\mathrm{prin}}$ of \cite{GHKK} for $\widetilde B$.
It has an extra invertible coordinate $y_i$ for each $i\in M\sqcup F$.
By \cite[Propositions~8.28 and~8.20]{GHKK}, $V$ and $V^\vee$ have enough global monomials, and the theta functions $\vartheta_q$, $q\in V^\vee(\bZ^T)$, form a basis of $\cO(V)$.
By \cite[Lemma~8.19]{GHKK}, there are global monomials $\vartheta_{p_1},\dots,\vartheta_{p_r}$ on $V^\vee$ such that, with $\varphi_i=\vartheta_{p_i}^T$, the set
\[
 \Xi=\{q\in V^\vee(\bR^T)\mid\varphi_i(q)\ge-1\ \text{for all }i\}
\]
is bounded.
As in the proof of \cite[Theorem~8.23]{GHKK}, we may take each $p_i$ in the interior of a cluster cone $\mathcal C^+_{\mathbf s_i}$.
Then $\vartheta_{p_i}$ is a character of the seed torus of $\mathbf s_i$, so $\varphi_i$ is linear in the coordinates of $\mathbf s_i$.
Each $\varphi_i$ is decreasing \cite[Proposition~8.14]{GHKK}.
So each set $\{\varphi_i\ge-1\}$ is a positive polytope \cite[Lemma~8.7]{GHKK}, and so is their intersection $\Xi$.
For $n\ge0$, put
\[
 \widetilde F_n=\bigoplus_{q\in n\Xi(\bZ)}\bC\,\vartheta_q .
\]
The set $\Xi$ is bounded, and it contains a neighborhood of $0$, since $\varphi_i(0)=0>-1$.
So $\widetilde F_0=\bC$, $\dim\widetilde F_n<\infty$ and $\bigcup_n\widetilde F_n=\cO(V)$.
Since $\Xi$ is a positive polytope, $\widetilde F_n\widetilde F_l\subset\widetilde F_{n+l}$, and $\bigoplus_n\widetilde F_nt^n$ is finitely generated \cite[Theorem~8.23]{GHKK}.

Next we describe $\widetilde F_n$ by valuations.
For a Laurent polynomial $g=\sum_ua_uz^u$ in the coordinates of $\mathbf s_i$, put $v_i(g)=\min_{a_u\ne0}\varphi_i(u)$.
The function $\varphi_i$ strictly increases at every bend of a broken line \cite[Lemma~8.17]{GHKK}.
The broken line without bends contributes $z^q$.
So in the coordinates of $\mathbf s_i$,
\[
 \vartheta_q=z^q+\sum_{\varphi_i(u)>\varphi_i(q)}a_{q,u}z^u .
\]
Distinct $q$ give distinct leading monomials, so $v_i\bigl(\sum_qc_q\vartheta_q\bigr)=\min_{c_q\ne0}\varphi_i(q)$.
Hence
\begin{equation}\label{eq:theta-valuation}
 \widetilde F_n=\{g\in\cO(V)\mid v_i(g)\ge-n\ \text{for all }i\}.
\end{equation}

Next we set the coefficients to $1$.
On each seed torus of $V$, the locus $\{y=1\}$ is the seed torus of $\cA(\widetilde B)$, and the gluing maps restrict to those of $\cA(\widetilde B)$.
Restriction to this locus is a homomorphism $\mathrm{ev}_1\colon\cO(V)\to\sU(\widetilde B)$.
By assumption, $\sU(\widetilde B)=\sA(\widetilde B)$.
The map $\mathrm{ev}_1$ is onto, because its image contains the cluster variables and the inverses of the frozen variables.
Put $F_n=\mathrm{ev}_1(\widetilde F_n)$.
Then $F_0=\bC$, $\dim F_n<\infty$, $\bigcup_nF_n=\sA(\widetilde B)$ and $F_nF_l\subset F_{n+l}$.
The Rees algebra $\bigoplus_nF_nt^n$ is a quotient of $\bigoplus_n\widetilde F_nt^n$, hence finitely generated.

It remains to show that each $F_n$ is $T$-invariant.
Every exchange relation $x_kx_k'=\mathfrak m_++\mathfrak m_-$ of $\cA(\widetilde B)$ is $T$-homogeneous, so $\mathfrak m_+$ and $\mathfrak m_-$ have the same weight.
In $V$ the relation reads $x_kx_k'=y^{e_+}\mathfrak m_++y^{e_-}\mathfrak m_-$ for some exponents $e_\pm$.
So the $T$-action extends to $V$, with every $y_i$ of weight zero.
It acts diagonally on the monomials of every seed, so it preserves every $v_i$, and every $\widetilde F_n$ by~\eqref{eq:theta-valuation}.
Since $\mathrm{ev}_1$ is $T$-equivariant, each $F_n$ is $T$-invariant.
\end{proof}

\begin{lemma}\label{lem:quadratic-growth}
In the setting of Lemma~\ref{lem:theta-filtration}, let $(F_n)_{n\ge0}$ be the filtration constructed in its proof, and let $\{\cdot,\cdot\}$ be a Poisson bracket on $\sA(\widetilde B)$ which is log-canonical in every seed.
Then $\{F_n,F_l\}\subset F_{n+l}$.
Moreover, there is a closed embedding $\iota=(\iota_1,\dots,\iota_N)\colon\cA(\widetilde B)\hookrightarrow\bC^N$ whose coordinates are $T$-weight vectors and satisfy
\[
 \{\iota_\mu,\iota_\nu\}\in\operatorname{span}_{\bC}\{1,\iota_\rho,\iota_\rho \iota_\sigma\mid 1\le\rho,\sigma\le N\}\qquad(1\le\mu,\nu\le N).
\]
\end{lemma}

\begin{proof}
We keep $V$, $\mathbf s_i$, $v_i$, $\widetilde F_n$, $\mathfrak m_\pm$, $e_\pm$ and $\mathrm{ev}_1$ from the proof of Lemma~\ref{lem:theta-filtration}.
We extend the bracket to $V$, show that it respects $\widetilde F_\bullet$, and then take a Veronese subalgebra.

We first extend the bracket to $V$.
On the initial seed torus of $V$, we keep the brackets of the $x_i$, and make every $y_i$ a Casimir.
Let $x_kx_k'=\mathfrak m_++\mathfrak m_-$ be an exchange relation of $\cA(\widetilde B)$.
For $j\ne k$, write $\{x_k,x_j\}=c_{kj}x_kx_j$ and $\{x_k',x_j\}=c_{kj}'x_k'x_j$ in $\sA(\widetilde B)$.
The monomials $\mathfrak m_\pm$ do not involve $x_k$, so $\{\mathfrak m_\pm,x_j\}=\gamma_\pm\mathfrak m_\pm x_j$ for constants $\gamma_\pm$.
By the Leibniz rule,
\[
 \gamma_+\mathfrak m_+x_j+\gamma_-\mathfrak m_-x_j=\{x_kx_k',x_j\}=(c_{kj}+c_{kj}')(\mathfrak m_++\mathfrak m_-)x_j .
\]
So $\gamma_+=\gamma_-=c_{kj}+c_{kj}'$, since $\mathfrak m_+\ne\mathfrak m_-$, or $\mathfrak m_+=\mathfrak m_-=1$ and $\gamma_\pm=0$.
In $V$, $x_kx_k'=y^{e_+}\mathfrak m_++y^{e_-}\mathfrak m_-$ and the $y_i$ are Casimirs, so $\{x_kx_k',x_j\}=\gamma_+x_kx_k'x_j$.
Hence $\{x_k',x_j\}=(\gamma_+-c_{kj})x_k'x_j=c_{kj}'x_k'x_j$ in $V$ as well.
By induction on mutations, the bracket is log-canonical in every seed of $V$.
So it preserves the Laurent polynomial ring of every seed of $V$, and hence their intersection $\cO(V)$.

Next we show $\{F_n,F_l\}\subset F_{n+l}$.
In the seed $\mathbf s_i$, the bracket reads $\{z^u,z^w\}=\Pi_i(u,w)\,z^{u+w}$ for a bilinear form $\Pi_i$.
Since $\varphi_i(u+w)=\varphi_i(u)+\varphi_i(w)$, this gives $v_i(\{g,g'\})\ge v_i(g)+v_i(g')$.
So $\{\widetilde F_n,\widetilde F_l\}\subset\widetilde F_{n+l}$ by~\eqref{eq:theta-valuation}.
The map $\mathrm{ev}_1$ is Poisson, because the $y_i$ are Casimirs.
Hence $\{F_n,F_l\}\subset F_{n+l}$.

It remains to construct $\iota$.
For sufficiently divisible $d>0$, the Veronese algebra $\bigoplus_{n\ge0}F_{nd}t^{nd}$ is generated by $F_dt^d$.
Thus $F_{nd}=(F_d)^n$, where $(F_d)^n$ is the linear span of $n$-fold products.
Since $F_d$ is $T$-invariant and contains $F_0=\bC$, it has a $T$-weight basis $1,\iota_1,\dots,\iota_N$.
Since $\sA(\widetilde B)=\bigcup_nF_{nd}=\bC[\iota_1,\dots,\iota_N]$, the map $\iota$ is a closed embedding.
Finally, $\{\iota_\mu,\iota_\nu\}\in F_{2d}=(F_d)^2=\operatorname{span}_{\bC}\{1,\iota_\rho,\iota_\rho\iota_\sigma\}$.
\end{proof}

For sink recurrent $\widetilde B$, we have $\sA(\widetilde B)=\sU(\widetilde B)$  \cite[Proposition~3.7]{GLSB}.  We also have:

\begin{lemma}\label{lem:sink-recurrent-large}
If $B$ is sink recurrent,  $\cA(\widetilde B)$ has large cluster complex.
\end{lemma}

\begin{proof}
The mutable matrices with a reddening sequence form a mutation-closed class.
Indeed, if $(k_1,\dots,k_r)$ is a reddening sequence of $B$ and $k\in M$, then so is $(k,k,k_1,\dots,k_r)$, and \cite[Theorem~2.2.4]{BHIT} gives a reddening sequence of $\mu_kB$ from it.
The class contains the zero matrices.
Let $k$ be a sink of $B$, and let $\mathbf k$ be a reddening sequence of $B|_{M\setminus\{k\}}$.
Since $B_{jk}\ge0$ for all $j$, the sequence $(k,\mathbf k)$ is a reddening sequence of $B$ \cite[Theorem~4.5 and Remark~4.6]{CL}.
So every sink-recurrent $B$ has a reddening sequence, by Definition~\ref{def:sink-recurrence}.
At its end all $c$-vectors are nonpositive.
So the last cluster chamber is the negative chamber \cite[Lemma~5.12]{GHKK}, and $\cA(\widetilde B)$ has large cluster complex \cite[Lemma~8.29]{GHKK}.
\end{proof}

We turn to the comparison of Liouville structures.  For a smooth affine variety $X$ with complex structure $I$ and closed embedding $\iota=(\iota_1,\dots,\iota_N)\colon X\hookrightarrow\bC^N$, write $\iota_\nu=x_\nu+iy_\nu$ and put
\[
 R=1+\tfrac14\abs{\iota}^2,\qquad
 \laff=-\dd R\circ I=\tfrac12\sum_{\nu=1}^N\bigl(x_\nu\,\dd y_\nu-y_\nu\,\dd x_\nu\bigr),\qquad
 g=\dd\laff(\cdot,I\cdot).
\]
The metric $g$ is induced by the Euclidean metric, and $Z_{\laff}R\le R-1$.
The Liouville structures obtained from different algebraic closed embeddings are deformation equivalent~\cite[Lemma~5.18]{McLean-growth}.

In the unimodular case, we have:

\begin{proposition}\label{prop:affine unimodular case} 
Let $\widetilde B=\binom BC$ be unimodular, with $B$ sink recurrent, and let $\lad$ be the Liouville form of Proposition~\ref{prop:lambda_ad}.
Then there is a closed embedding $\iota\colon\cA(\widetilde B)\hookrightarrow\bC^N$, a diffeomorphism $\phi$ of $\cA(\widetilde B)$ isotopic to the identity, and a function $f$ such that
$\phi^*\lad=\laff+\dd f
$.
\end{proposition}
\begin{proof}
Let $\widetilde B$ be unimodular, with completion $(E,\widetilde D)$.
By the proof of Proposition~\ref{prop:Omega} and Lemma~\ref{lem:completion-mutation}, $\Omega$ has the form~\eqref{eq:forms} in every seed, for the mutated completion.
So the Poisson bracket of $\Omega$ is log-canonical in every seed.
Lemma~\ref{lem:quadratic-growth}, with $T$ trivial, gives $\iota$ whose coordinate brackets are polynomials of degree at most two in $\iota$.
Since $\dd\lad=\Im\Omega$ by~\ref{P1}, \cite[Corollary~1.6]{HS26} gives $\phi$ and $f$ with $\phi^*\lad=\laff+\dd f$.
\end{proof}
\begin{corollary}
$\Fuk\bigl(\cA(\widetilde B),\lad\bigr)\simeq\Fuk\bigl(\cA(\widetilde B),\laff\bigr)$.
\end{corollary}

In general:

\begin{proposition}\label{prop:affine} 
Let $\widetilde B=\binom BC$ have full rank, with $B$ sink recurrent, and let $\alpha$ be the Liouville form of~\ref{Q3}.
There are a closed embedding $\iota\colon\cA(\widetilde B)\hookrightarrow\bC^N$, a diffeomorphism $\phi$ of $\cA(\widetilde B)$ isotopic to the identity, and a function $f$ such that
$
 \phi^*\alpha=\laff+\dd f
$.

\end{proposition}

\begin{proof}
Use $P$, $\lambda_P$, $\omega_P$, $S$ and $\varrho_{s,q}$ of Section~\ref{sec:fullrank-geometry}, and identify $\cA(\widetilde B)=F^{-1}(1)$ by Lemma~\ref{lem:phase-cover}.
Apply Lemma~\ref{lem:quadratic-growth} to $\cA_{\mathrm{prin}}(B)$, with the bracket of $\Omega$ for the completion~\eqref{eq:principal-completion} and the torus $H$ acting by~\eqref{eq:h-action}.
Let $\iota_0$ be the resulting embedding.
Put $\iota=(\iota_0,U,U^{-1})$ on $P$, and restrict it to $\cA(\widetilde B)$.

Since $\Im\Omega(v,Iv)=0$, the definition of $\lambda_P$ gives 
\[
 \dd\alpha(v,Iv)=\sum_{a\in F}\Bigl(\bigl(\dd\log|U_a|(v)\bigr)^2+\bigl(\dd\arg U_a(v)\bigr)^2\Bigr)\ge0 .
\]
Thus $\dd\alpha$ is weakly tamed by $I$.

For the growth estimate, let $\Pi_P=\omega_P^{-1}$ and let $\Phi$ be the diffeomorphism of Lemma~\ref{lem:phase-cover}(ii).
Inverting $\Phi^*\pi^*\omega_P=\dd\alpha+\dd p\wedge\dd q$ at $(x,0,0)$ and pushing forward by $\pi\Phi$ gives
\[
 \Pi_P=(\dd\alpha)^{-1}\pm\sum_{i=1}^m\xi_i\wedge\eta_i,\qquad
 \xi_i=\partial_{q_i}\varrho_{0,q}(x)\big|_{q=0},\qquad
 \eta_i=\partial_{p_i}\varrho_{-(S^{\T})^{-1}p,0}(x)\big|_{p=0}.
\]
The coordinate brackets of $\iota_0$ have degree at most two, and the fields $\partial_{\log|U_a|}$ and $\partial_{\arg U_a}$ act diagonally on $U_a^{\pm1}$.
Hence $\iota_*\Pi_P$ has coefficients of degree at most two in the real and imaginary parts of $\iota$, as in \cite[proof of Corollary~1.6]{HS26}.
By~\eqref{eq:phase-coordinates} and the choice of $\iota_0$, the vectors $\dd\iota(\xi_i)$ and $\dd\iota(\eta_i)$ are linear in $\iota$.
So the correction term has the same quadratic bound.
Since $g$ is the induced Euclidean metric,
\[
 \norm{(\dd\alpha)^{-1}}_g\le C\bigl(1+\abs{\iota}^2\bigr)\le C'R .
\]
Apply \cite[Theorem~1.4]{HS26} with $J=I$ and growth function $r\mapsto r$ to obtain $\phi^*\alpha=\laff+\dd f$.

In both cases $\phi$ is a symplectomorphism from $\dd\laff$ to $\dd\alpha$, respectively to $\dd\lad$.
The Liouville form $\laff$ is Weinstein, since $R$ is Morse.
The forms $\alpha$ and $\lad$ are Morse--Bott Weinstein by~\ref{Q3} and Proposition~\ref{prop:lambda_ad}.
\end{proof}
\begin{corollary}
$\Fuk\bigl(\cA(\widetilde B),\alpha\bigr)\simeq\Fuk\bigl(\cA(\widetilde B),\laff\bigr).$
\end{corollary}

\section{On geometricity of  $\mathbb{G}_m$ actions}

The following is likely well known, but we did not find a reference.  

\begin{lemma}\label{lem:gm-action-geometric}
Let $X=\Spec A$ be smooth affine over $\bC$.
Every algebraic $\mathbb G_m$-action on $\QCoh(X)$ is induced by an algebraic action on $X$.
\end{lemma}

\begin{proof}
Consider the coaction $
 \rho:\QCoh(X)\to\QCoh(\mathbb G_m\times X)$.
Let $p:\mathbb G_m\times X\to X$ be the projection.
The $\bC[t^{\pm1}]$-linear extension of $\rho$ is an autoequivalence $\Phi$ of $\QCoh(\mathbb G_m\times X)$, characterized by $\Phi(p^*F)=\rho(F)$.  By derived Morita theory \cite[Theorem~8.9]{Toen-Morita} this is given by an invertible $A[t^{\pm1}]$-bimodule over $\bC[t^{\pm1}]$; by the description of the derived Picard group (see e.g. \cite[Corollary~5.7]{Negron-DPic}), we find:
\[
 \Phi(-)\simeq L[n]\otimes f^*(-),
 \qquad f\in\operatorname{Aut}_{\mathbb G_m}(\mathbb G_m\times X),
\]
where $L$ is a line bundle and $n$ is locally constant on $\mathbb G_m\times X$.
The unit gives $n|_{\{1\}\times X}=0$, $L|_{\{1\}\times X}\simeq\cO_X$, and $f_1=\id_X$.
Every connected component meets $\{1\}\times X$, so $n=0$.
Since $A$ is regular, pullback gives $
 \operatorname{Pic}(A)\xrightarrow{\sim}\operatorname{Pic}(A[t^{\pm1}])$
by \cite[Remark~1.5.2]{Weibel-Pic}.
Thus $L$ is trivial; choose a trivialization agreeing with the unit at $t=1$.

The multiplication coherence now identifies $f_s^*f_t^*$ with $f_{st}^*$.
Comparing their maps on functions gives $f_t\circ f_s=f_{st}$.
Since $st=ts$, the map $a(t,x)=f_t(x)$ is an algebraic $\mathbb G_m$-action on $X$.
Relative to the canonical composition of pullbacks, its categorical multiplication differs by a unit
\[
 c\in A[s^{\pm1},t^{\pm1}]^\times,
 \qquad c(1,t)=c(s,1)=1.
\]
Each connected component of $X$ is integral.
On each component, $c=b s^r t^q$ with $b$ a unit and $r,q\in\bZ$.
The two normalizations force $b=1$ and $r=q=0$.
Hence the chosen trivialization respects multiplication.
There are no further coherence obstructions by \cite[Corollary~8.13]{Toen-Morita}: the autoequivalence space is $1$-truncated.
\end{proof}

\section{Sink recurrence for $\mathrm{SL}_2$ marked surfaces}\label{app:surfaces}

Let $S$ be a compact connected oriented surface with boundary, and let $\mathbb M\subset\partial S$ be a finite set, with at least two elements, and meeting every boundary component.  Fix a triangulation \(T\) whose vertices are in this set.

In \cite{FG06}, the corresponding seed for $\cA_{\mathrm{SL}_2, S}$ is the following: 
\begin{align}\label{eq:surface-B}
M = \{\text{Internal edges of } T\}, \qquad F = \{ \text{Boundary edges of } T\}, \notag
\\
 \widetilde B_T=\binom BC,\qquad
 (\widetilde B_T)_{ij}=\#\{\Delta: j\text{ follows }i\}-\#\{\Delta: i\text{ follows }j\},
\end{align}
where $i\in M\sqcup F$, $j\in M$, $\Delta$ runs over the triangles of $T$ with sides $i$ and $j$, and ``follows'' refers to the clockwise order of the sides of $\Delta$.
This is the seed of the moduli space $\mathcal A_{\mathrm{SL}_2,S}$ of Fock and Goncharov \cite[Definition~1.11 and Section~10.2]{FG06}, with $m=2$.
In their notation, $J^2_T=M$ \cite[Figure~1.6]{FG06}, the frozen set is $F$, and $(\widetilde B_T)_{ij}=\varepsilon_{ji}$.
Flips of $T$ are mutations of $\widetilde B_T$ \cite[Section~10.3 and Proposition~10.1]{FG06}. Any two ideal triangulations are related by flips. %

Lam and Speyer have shown that this $\widetilde{B}_T$ is ``Louise'' \cite[Theorem~1.3]{Lam-Speyer} using an argument of Muller \cite[Theorems~10.5--10.6]{M}. Here we show:

\begin{proposition}\label{prop:surface}
The matrix $\widetilde B_T$ has full rank, and $B$ is sink recurrent.
\end{proposition}

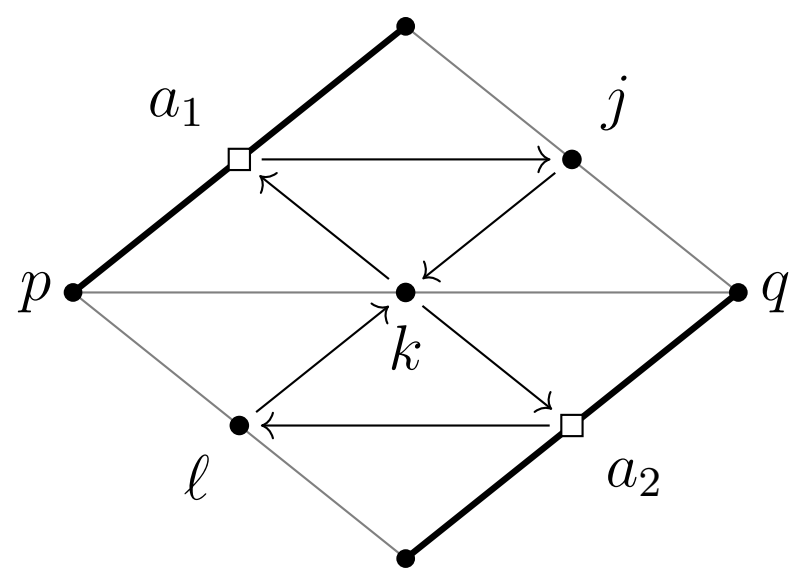
\begin{figure}
\centering
\begin{tikzpicture}[scale=1.1,
 mut/.style={circle,fill,inner sep=1.3pt},
 fro/.style={rectangle,draw,fill=white,inner sep=2pt},
 arr/.style={->,shorten >=2pt,shorten <=2pt}]
 \coordinate (p) at (0,0);
 \coordinate (r) at (2,1.6);
 \coordinate (q) at (4,0);
 \coordinate (s) at (2,-1.6);
 \draw[gray] (p)--(q) (r)--(q) (s)--(p);
 \draw[very thick] (p)--(r) (q)--(s);
 \foreach \v in {p,q,r,s} \fill (\v) circle (1.6pt);
 \node[left] at (p) {$p$};
 \node[right] at (q) {$q$};
 \node[fro,label=above left:{$a_1$}] (A1) at (1,0.8) {};
 \node[mut,label=above right:{$j$}] (J) at (3,0.8) {};
 \node[fro,label=below right:{$a_2$}] (A2) at (3,-0.8) {};
 \node[mut,label=below left:{$\ell$}] (L) at (1,-0.8) {};
 \node[mut,label=below:{$k$}] (K) at (2,0) {};
 \draw[arr] (A1)--(J);
 \draw[arr] (J)--(K);
 \draw[arr] (K)--(A1);
 \draw[arr] (K)--(A2);
 \draw[arr] (A2)--(L);
 \draw[arr] (L)--(K);
\end{tikzpicture}
\caption{The two triangles at the arc $k=pq$. An arrow $i\to j$ joins two sides of a triangle when $j$ follows $i$ in the clockwise order. The thick sides $a_1,a_2$ lie on $\partial S$, and squares mark frozen vertices.}
\label{fig:band}
\end{figure}

\begin{proof}
The matrix $\widetilde B_T$ is the exchange matrix of a quantum seed \cite[Proposition~7.8]{Muller-skein}, so it has full rank \cite[Proposition~7.2]{Muller-skein}.

We prove sink recurrence by induction on $\abs M$, for all $(S,\mathbb M)$ as above at once.
Sink recurrence is mutation invariant, so we may choose any triangulation $T$.
If $k$ is a source of $B$, then $k$ is a sink of $\mu_kB$, and $\mu_kB$ agrees with $B$ on $M\setminus\{k\}$.
So in Definition~\ref{def:sink-recurrence} a source $k$ may replace a sink.

If $S$ is a disk, take a fan triangulation.
Its mutable matrix is acyclic.

Otherwise, \cite[Lemma~4.16, case (2a) of the proof]{Ladkani} gives $T$ and $k\in M$ with the following properties.
The arc $k$ joins two distinct marked points, and $S$ cut along $k$ is connected.
The two triangles at $k$ are $\{a_1,j,k\}$ and $\{a_2,\ell,k\}$ with $a_1,a_2\in F$ and $j,\ell\in M$, possibly $j=\ell$.
They form a quadrilateral with diagonal $k$ and opposite sides $a_1,a_2$, see Figure~\ref{fig:band}.
Its sides in clockwise order are $a_1,j,a_2,\ell$ or $a_1,\ell,a_2,j$.
In the first case $k$ follows $j$ and $\ell$ in their triangles, so \eqref{eq:surface-B} gives $(\widetilde B_T)_{jk}>0$ and $(\widetilde B_T)_{\ell k}>0$.
In the second case $j$ and $\ell$ follow $k$, and both entries are negative.
So $k$ is a sink or a source of $B$, and $N(k)=\{j,\ell\}$.

Let $S'$ be $S$ cut along $k$, and let $S''$ be $S'$ minus the two triangles at $k$.
Turning an arc into a boundary segment deletes its vertex from $B$.
So the triangulations $T\setminus\{k\}$ of $S'$ and $T\setminus\{k,j,\ell\}$ of $S''$ have the mutable matrices
\[
 B|_{M\setminus\{k\}},\qquad B|_{M\setminus(\{k\}\sqcup N(k))} .
\]
The surface $S'$ is connected by the choice of $k$.
In $S'$, the triangle with sides $a_1,j,k$ meets the rest only along $j$, and the other triangle only along $\ell$.
So $S''$ is connected, or empty if $j=\ell$.
Both surfaces contain the two endpoints of $k$ as distinct marked points, and every boundary component meets the marked points.
They have fewer internal arcs than $S$.
By induction both restrictions are sink recurrent, and so is $B$.
\end{proof}

\bibliographystyle{plain}
\bibliography{refs}

\end{document}

%% file: macros.tex
\newcommand{\cP}{\mathcal{P}}
\newcommand{\bC}{\mathbb{C}}

\newcommand{\cF}{\mathcal{F}}

\newcommand{\Bl}{\mathrm{Bl}}

\DeclareMathOperator{\Spec}{Spec}

\newcommand{\cN}{\mathcal{N}}

\newcommand{\cA}{\mathcal{A}}
\newcommand{\sA}{\mathscr{A}}
\newcommand{\sU}{\mathscr{U}}
\newcommand{\sC}{\mathscr{C}}
\newcommand{\cK}{\mathcal{K}}
\newcommand{\cQ}{\mathcal{Q}}

\newcommand{\bR}{\mathbb{R}}

\newcommand{\bZ}{\mathbb{Z}}

\newcommand{\cO}{\mathcal{O}}

\DeclareMathOperator{\Coh}{Coh}
\DeclareMathOperator{\QCoh}{QCoh}
\DeclareMathOperator{\Sh}{Sh}
\DeclareMathOperator{\Loc}{Loc}

\DeclareMathOperator{\Hom}{Hom}

\newcommand{\cX}{\mathcal{X}}

\newcommand{\cC}{\mathcal{C}}
\newcommand{\bQ}{\mathbb{Q}}

\newcommand{\Perf}{\mathrm{Perf}}

\newcommand{\red}[1]{\textcolor{red}{#1}}

\let\Im\relax
\DeclareMathOperator{\Im}{Im}
\let\Re\relax
\DeclareMathOperator{\Re}{Re}

\newcommand{\colim}{\operatorname*{colim}}
\DeclareMathOperator{\Log}{Log}

\newcommand{\R}{\mathbb R}
\newcommand{\C}{\mathbb C}
\newcommand{\id}{\mathrm{id}}
\newcommand{\abs}[1]{\lvert #1\rvert}

\newcommand{\norm}[1]{\lVert #1\rVert}

\DeclareMathOperator{\coker}{coker}

\let\c\relax
\newcommand{\c}{\check}

\DeclareMathOperator{\Crit}{Crit}

\DeclareMathOperator{\Skel}{Skel}
\DeclareMathOperator{\Hess}{Hess}
\DeclareMathOperator{\supp}{supp}
\DeclareMathOperator{\rank}{rank}
\DeclareMathOperator{\diag}{diag}
\DeclareMathOperator{\Mat}{Mat}

\DeclareMathOperator{\ind}{ind}

\DeclareMathOperator{\Zero}{Zero}
\DeclareMathOperator{\fib}{fib}

\DeclareMathOperator{\codim}{codim}

\newcommand{\lad}{\lambda_{\mathrm{ad}}}
\newcommand{\laff}{\lambda_{\mathrm{aff}}}
\newcommand{\ltor}{\lambda_{\mathrm{tor}}}
\newcommand{\lN}{\lambda_{\cN}}
\newcommand{\T}{\mathsf T}
\newcommand{\capp}{\mathrm{cap}}

\newcommand{\dd}{\mathrm{d}}